\documentclass[UTF17]{amsart}
\usepackage{amssymb,amsmath,amsthm,mathrsfs,multirow,xcolor,framed,url,enumitem}
\usepackage{float}
\usepackage{amsaddr}
\usepackage{lipsum} 
\usepackage{diagbox} 
\usepackage{tabularx}
\usepackage{longtable} 
\usepackage{threeparttablex} 
\usepackage{booktabs} 
\usepackage[title]{appendix}
\usepackage{tikz}
\usepackage[all]{xy}
\usepackage{caption}

\usepackage[ 
bookmarks=true,
bookmarksnumbered=true，
bookmarksopen=false,
hypertexnames=false,
breaklinks=true,
colorlinks=true,
pagebackref=true,
hyperindex=true,
plainpages=false
pdfauthor={Luo, Wang, Zhou },
pdftitle={Ternary},
pdfsubject={ },
pdfkeywords={}
]
{hyperref}

\newtheorem{theorem}{Theorem}[section]
\newtheorem{lemma}[theorem]{Lemma}
\newtheorem{cor}[theorem]{Corollary}

\newtheorem{prop}[theorem]{Proposition}

\theoremstyle{definition}
\newtheorem{definition}[theorem]{Definition}

\usepackage{tikz}
\usetikzlibrary{arrows.meta}
\theoremstyle{remark}
\newtheorem{remark}{Remark}  

\numberwithin{equation}{section}

\usepackage{hyperref}
\newcommand\nutwid{\overset{\text{\lower3pt\hbox{$\sim$}}}\nu}

\newcommand{\tr}{\operatorname{tr}}
\newcommand{\n}{\operatorname{n}}

\newcommand{\discrd}{\operatorname{discrd}}

\newcommand{\Aut}{\operatorname{Aut}}

\allowdisplaybreaks

\begin{document}
	
	\title[Ternary quadratic genera and generalized Hurwitz class numbers]{Ternary quadratic genera and generalized Hurwitz class numbers via quaternion orders}

	\author{Yifan Luo$^{1}$, Tianze Wang$^{1,2}$ and Haigang Zhou$^{3}$}
	
	\address{$^{1}$Institute of Mathematics, Henan Academy of Sciences, Zhengzhou, China, 450046}
	\email{kotori@hnas.ac.cn}
	
	\address{$^{2}$School of Mathematics and Statistics, North China University of Water Resources\\ and Electric Power, Zhengzhou, China, 450046}
	\email{wtz@ncwu.edu.cn}
	
	\address{$^{3}$School of Mathematical Sciences, Tongji University, Shanghai, China, 200092}
	\email{haigangz@tongji.edu.cn}
	
\begin{abstract}
We classify the genera of primitive positive definite ternary quadratic
forms in the extremal cases of Lehman's level and discriminant conditions.
The genera are parametrized by their level, discriminant, and genus
conductor, which also distinguishes local forms having the same isotropy.
Using the correspondence with quaternion orders of level $(N_1,N_2)$,
we express their weighted representation numbers and class numbers in
terms of generalized modified Hurwitz class numbers.
The formulas allow arbitrary admissible prime power exponents and
include the additional  genera obtained by Watson's transformation.
As an application, we determine all $255$ one-class genera in these
families and give an explicit representation formula for each.
Of these formulas, $110$ concern nonsquarefree $N_1N_2$.
\end{abstract}

\keywords{ternary quadratic forms, quaternion algebras, Siegel--Weil formula, representation numbers}
\subjclass[2020]{11E20, 11E25, 11R52, 11E41}
\maketitle

\section{Introduction}\label{Sec:1}

Ternary quadratic forms, especially positive definite ones, occupy a distinguished position in number theory. They are the simplest genuine higher-dimensional forms beyond the binary case and arise naturally in the theory of quadratic fields, the geometry of numbers, modular forms, and local--global questions for rational representations. At the same time, their representation theory exhibits deep phenomena that do not occur for binary forms, including spinor exceptions and failures of the local--global principle for integral representations.

The classification of ternary quadratic forms over the integers has close connections to the arithmetic of quaternion algebras and spinor genera. Over $\mathbb{Q}$, the Hasse--Minkowski theorem reduces global equivalence to a collection of local invariants. Over $\mathbb Z$, however, two forms can be equivalent over
$\mathbb Z_p$ for every prime $p$ and over $\mathbb R$,
yet inequivalent over $\mathbb Z$.
 This phenomenon is reflected in the theory of genera and spinor genera, developed in the work of Brandt, Eichler, Kneser, and Schulze-Pillot. One of the most fundamental classification theories is based on a set of invariants introduced by Eichler~\cite{Eic74} and reformulated by Shimura~\cite{Shi06} in a more natural way over algebraic number fields.

A fundamental step toward the systematic classification of ternary quadratic forms was made by Lehman~\cite{Leh92}. For a primitive positive definite ternary quadratic form \(f\), Lehman introduced the level \(N_f\), together with the discriminant \(d_f\), as a basic invariant, and gave a complete formulation of local invariants that determine the genus of \(f\). This framework clarifies the structure of genera of ternary forms and connects the theory with quaternion algebras and modular forms of weight \(3/2\).
Write
\begin{equation}\label{eq:Nd-factorization}
	N_f = 2^{n_0} \prod_{i=1}^k p_i^{n_i}, \qquad
	d_f = 2^{d_0} \prod_{i=1}^k p_i^{d_i},
\end{equation}
for the prime factorizations of \( N_f \) and \( d_f \), where $p_i\geq 3$ are distinct primes for $1\leq i\leq k$.
Then Lehman's work~\cite[Theorem 2]{Leh92} shows that $n_0\geq2$ and that the relations between the exponents $n_i$ and $d_i$ can be formulated as a combination of the conditions: (1) \( n_i \leq d_i \leq 2n_i \) for $1\leq i\leq k$; (2) \( d_0 = n_0-2 \), \( d_0 = 2n_0 \), or \( n_0 \leq d_0 \leq 2n_0-2 \); and (3) additional parity constraints (see Lemma~\ref{the:N&d} below). Taken together, these constraints describe the precise admissible region for the level \(N_f\) and the discriminant \(d_f\).

The squarefree-level results of~\cite{LZ24} suggest the study of arbitrary
prime power exponents. We do this in this paper for the extremal cases of Lehman's
conditions, namely
\begin{itemize}
	\item \( d_i = n_i \) or \( d_i = 2n_i \) for $1\leq i\leq k$, and
	\item \( d_0 = n_0-2 \) or \( d_0 = 2n_0 \).
\end{itemize}
Under these conditions and the parity restriction in
Lemma~\ref{the:N&d}, we have $N_f=4N$ with $N$ not a perfect square.
The admissible discriminants are precisely
\begin{equation}
 d_f=16N^2/\delta,\quad \delta\parallel16N.
\end{equation}
Here $N_1\parallel N$ means $N_1\mid N$ and $(N_1,N/N_1)=1$.
We also consider the genera of level $2N$ and discriminant
$N^2/\delta'$ when $2\parallel N$ and $\delta'\parallel N/2$.
These additional genera are obtained by Watson's transformation and
need not satisfy the extremal condition at $2$.
The corresponding squarefree cases recover the formulas of~\cite{LZ24}.

Throughout this paper, $N_1$ always denotes positive divisors of $N$ of form $N_1=s^2N_0$ with $N_1\parallel N$, where $s$ is a positive integer and $N_0$ is squarefree with $e(N_0)$ odd, i.e.\ $\mu(N_0)=-1$. That is, we will always assume
\begin{equation}\label{equ1.2}
	N=N_1N_2,\quad N_1=s^2N_0,\quad \gcd(N_1,N_2)=1,
\end{equation}
with $\mu(N_0)=-1$, and $N_2$ and $s$ positive integers.
Here $e(N_0)$ denotes the number of distinct prime factors of $N_0$; and $\mu(\cdot)$ is the M\"{o}bius function. 

With these in place, our first main result can be stated as the following theorem, which gives a complete classification and characterization for the genera of primitive positive definite ternary quadratic forms in all extremal cases of Lehman's Theorem~\cite[Theorem 2]{Leh92}.

\begin{theorem}\label{the:genera}
	Let $N$ be a positive integer that is not a perfect square, and let $\delta $ be a positive divisor of $16N$ with $\delta \parallel 16N$. Let $N_1$ denote positive divisors of $N$ as in ~\eqref{equ1.2}. Then for each $\delta $ there are exactly $2^{e(N)-1}$ genera $G_{4N,\,16N^2/\delta ,\,N_1}$ in the set of primitive positive definite ternary quadratic forms $f$ of level $N_f=4N$ and discriminant $d_f=16N^2/\delta $. Here $G_{4N,\,16N^2/\delta ,\,N_1}$ denotes the genus of level $4N$, discriminant $16N^2/\delta $, and genus conductor $N_1$ (see Definitions {\rm ~\ref{def:genus}} and {\rm ~\ref{def:genus111}}). Moreover, if in addition $2\parallel N$, and $\delta '$ is an odd positive integer with $\delta '\parallel N/2$, then for each $\delta '$ there are exactly $2^{e(N)-1}$ genera $G_{2N,\,N^2/\delta ',\,N_1}$ in the set of primitive positive definite forms $f$ of level $N_f=2N$ and discriminant $d_f=N^2/\delta '$.
\end{theorem}

\begin{remark}
The genus conductor records the genus symbols that are not determined by level and discriminant. A form in a genus of genus conductor $N_1$ is anisotropic
exactly at primes dividing $N_0=N_1/s^2$. The conductor $N_1$ also distinguishes genera at places of even valuations with same isotropy. 
For the base genus $G_{4N,\,16N^2,\,N_1}$ (i.e., $\delta=1$) corresponding to an order $\mathcal O$ of level $(N_1,N_2)$, Proposition~\ref{prop:Eichlersym} gives
\[
 N_1=\prod_{\substack{p\mid N\\(\mathcal O/p)=-1}}p^{v_p(N)}.
\]
This description uses the base order (the order corresponding to a base genus); its Eichler symbols can change
after partial reciprocity of orders and quadratic forms.
\end{remark}

\vskip 0.2cm
Parallel to classification stands the problem of representation numbers, or more generally the arithmetic of quadratic forms. This problem traces back to Fermat, Euler, Lagrange, and Legendre, and reached a foundational milestone in Gauss's \emph{Disquisitiones Arithmeticae}. A fundamental question in this area is to understand which integers are represented by a given form $f$, and to relate the corresponding representation numbers to arithmetic invariants attached to the form. For binary quadratic forms, such a connection is expressed elegantly by the Hurwitz class number formula, which links representation numbers to class numbers of imaginary quadratic fields. For ternary forms the analogous theory is more complex, but also richer, owing to the presence of spinor exceptions and the failure of the local--global principle for integral representations. Gauss established an explicit relation between the number of representations of $n$ as a sum of three squares and class numbers of binary quadratic forms. The celebrated results of Duke and Schulze-Pillot on the representation of integers by ternary quadratic forms, building on the modularity of theta series, have provided profound answers linking analytic methods to algebraic ones.

The even Clifford construction gives a bijection between equivalence classes of positive definite integral ternary quadratic forms and isomorphism classes of orders in definite rational quaternion algebras; see \cite[Chapter~22]{Voi21}. This correspondence is compatible with localization. Consequently, the class number of a ternary genus equals the type number of the corresponding genus of orders.

This correspondence is not only theoretical but also practical for applications in arithmetic. The type number is intimately connected to the counting of elements with given trace and reduced norm, and such counting corresponds precisely to the representation number of a ternary quadratic form. Thus the arithmetic of representation numbers is encoded in the distribution of elements within quaternion orders.

The modularity of the generating series attached to quaternion orders is exploited in the foundational works of Gross~\cite{Gro87}, Boylan--Skoruppa--Zhou~\cite{BSZ19}, and Li--Skoruppa--Zhou~\cite{LSZ22}. Gross showed that, for a positive definite quaternion algebra ramified at a prime \( p \), the weighted sum of the theta series corresponding to maximal orders is an Eisenstein series of weight \( 3/2 \) whose Fourier coefficients are modified Hurwitz class numbers. This result was extended to arbitrary squarefree numbers by Boylan, Skoruppa, and the third author~\cite{BSZ19}, and to Eichler orders with squarefree level by Li, Skoruppa, and the third author~\cite{LSZ22}. In each case the key insight is that a suitable weighted average of the representation numbers over a genus (or over a set of conjugate orders) yields modular (Jacobi) forms whose Fourier coefficients are given by modified Hurwitz class numbers.

The modified Hurwitz class numbers connect weighted genus representation numbers with type numbers of quaternion orders.
The squarefree case was treated in \cite{BSZ19,LSZ22,LZ24}.
In \cite{LZ26}, this method was extended to the family in which
every prime divisor of $N_1$ has odd valuation, with $N_2$
arbitrary and coprime to $N_1$.

In this paper, we allow positive even valuations in $N_1$ as well.
They correspond to residually inert orders in the split local
quaternion algebra. Together with Lehman's correspondences,
this extension covers all the extremal genera considered in
Theorem~\ref{the:genera}. The generalized quantities
$H^{(N_1,N_2;\,\delta)}$ evaluate their weighted representation
numbers, and the primitive-normalizer count then gives their
class numbers. The order--form correspondence used in this
argument is that established in \cite{LWZ26}.

We can now state our second and third main results, Theorems~\ref{the:repnum} and~\ref{the:classnum},
 which give explicit and uniform formulas, respectively, for the weighted number of representations and for
  the number of equivalence classes in a genus, valid for the general $N$ considered in this paper.

\begin{theorem}\label{the:repnum}
	Under the notations of Theorem ~{\rm \ref{the:genera}}, together with the generalization of 
	the modified Hurwitz class number $H^{(N_1,\,N/N_1;\,\delta )}(D)$ in Definition~ {\rm \ref{def:HCNd}}, we have for any integer $n\geq 0$
	\begin{equation*}
		\sum_{f\in G_{4N,\,16N^2/\delta ,\,N_1}}\frac{R_{f}(n)}{|\Aut(f)|}=2^{-e(N)-1}H^{(N_1,\,N/N_1;\,\delta )}(n),
	\end{equation*}
	where $R_f(n)$ denotes the number of representations of $n$ by $f$, $|\Aut(f)|$ is the number of automorphisms of $f$, 
	and the sum is over a complete set of representatives of the equivalence classes of genus $G_{4N,\,16N^2/\delta ,\,N_1}$.
	Moreover, if in addition $2\parallel N$, then
	\begin{equation*}
		\sum_{f\in G_{2N,\,N^2/\delta ',\,N_1}}\frac{R_{f}(n)}{|\Aut(f)|}=2^{-e(N)-1}H^{(N_1,\,N/N_1;\,\delta ')}(4n).
	\end{equation*}
\end{theorem}

\begin{remark}
	In the two formulas of this theorem there is a clear correspondence between the triple of the subscript of $G$ on the left,
	 which records the invariants of the genus $G$, and the triple of the superscript of $H$ on the right, 
	 which records the parameters of the Hurwitz class number $H$. Moreover, as we shall show in Remark~\ref{remark4.4}
	  under Definition~\ref{def:HCNd}, the generalized modified Hurwitz class numbers $H^{(N_1,N/N_1;\,\delta )}(n)$ 
	  used here are compatible with the existing ones (see, e.g.,~\cite{LZ24}).
\end{remark}

\begin{theorem}\label{the:classnum}
	Under the notations of Theorem ~{\rm \ref{the:repnum}}, let $|G_{4N,\,16N^2/\delta ,\,N_1}|$ denote the class number of genus $G_{4N,\,16N^2/\delta ,\,N_1}$. Then for any $\delta \parallel 16N$,
	\begin{equation*}
		|G_{4N,\,16N^2/\delta ,\,N_1}|=2^{-e(N)-1}\sum_{\substack{t\parallel N\\1\leq t\leq3}}\sum_{\substack{t\mid r\\r^2\leq4t}}H^{(N_1,\,N/N_1)}(4t-r^2)+2^{-e(N)-1}\sum_{\substack{t\parallel N\\t\geq4}}H^{(N_1,\,N/N_1;\,16^{\delta_{2\mid t}}t)}(4^{\delta_{2\nmid t}}),
	\end{equation*}
	where $\delta_{2\nmid t}=1$ if $2\nmid t$ and $0$ otherwise, $\delta_{2\mid t}=1$ if $2\mid t$ and $0$ otherwise. Equivalently,
	\begin{equation*}
		|G_{4N,\,16N^2/\delta ,\,N_1}|=2^{-e(N)-1}\sum_{\substack{t\parallel N\\t>0}}\sum_{\substack{t\mid r\\r^2\leq4t\\(t,r)\neq(4,\pm4)}}H^{(N_1,\,N/N_1)}(4t-r^2)\prod_{p\mid t}\frac{1-\left(\frac{\Delta(-4t)}{p}\right)/p}{B_p(t;\,N_1,N/N_1)C_p(t)},
	\end{equation*}
	where $\Delta(-4t)=-t_0$ if $-4t=-t_0t_1^2$ with $-t_0$ a fundamental discriminant, and
	\begin{align*}
		B_p(t;\,N_1,\,N/N_1)&=
		\begin{cases}
			(p-1)p^{v_p(t)/2-1}, & \text{if } v_p(t)\equiv0\!\!\!\pmod2,\ p\mid N_1,\\
			(p+1)p^{v_p(t)/2-1}, & \text{if } v_p(t)\equiv0\!\!\!\pmod2,\ p\mid (N/N_1),\\
			p^{(v_p(t)-1)/2},  & \text{if } v_p(t)\equiv1\!\!\!\pmod2,
		\end{cases}\\
		C_p(t)&=
		\begin{cases}
			2, & \text{if } p=2,\ 4\mid t,\ \text{and } \Delta(-4t)\equiv5\!\!\!\pmod8,\\
			1, & \text{otherwise},
		\end{cases}
	\end{align*}
	with $v_p(t)$ to be the valuation of $t$ at $p$.
	Moreover, if in addition $2\parallel N$, then for any $\delta '\parallel N/2$,
	\begin{equation*}
		|G_{2N,\,N^2/\delta ',\,N_1}|=|G_{4N,\,16N^2/\delta ,\,N_1}|.
	\end{equation*}
\end{theorem}

\begin{remark}
	The advantage of the first formula for $|G_{4N,\,16N^2/\delta ,\,N_1}|$ in Theorem~\ref{the:classnum} 
	is that, combined with Theorem~\ref{the:repnum}, it clearly reveals the connection between class numbers
	 and representation numbers, whereas the second, product form is more convenient for explicit computation. 
	 By Theorem~\ref{the:classnum}, it is clear that the class number is independent of $\delta$.
\end{remark}

The one-class genera of primitive positive definite ternary forms have already
been classified; see~\cite{JKI97} and
\cite[Theorem~25.4.6]{Voi21}. The application of our results identifies which
of these genera belong to the present families and supplies uniform
Hurwitz class number formulas for the weighted representation numbers.
We can give a list of all the pairs $(N_1,N/N_1)$ with class number $1$ as follows.

\begin{cor}
There are exactly $39$ pairs $(N_1,N/N_1)$ with class number $1$:
$(2,1)$, $(2,3)$, $(2,5)$, $(2,7)$, $(2,9)$, $(2,11)$, $(2,15)$, $(2,23)$, $(3,1)$, $(3,2)$, $(3,4)$, $(3,5)$, $(3,8)$, $(3,11)$, $(5,1)$, $(5,2)$, $(5,4)$, $(7,1)$, $(7,3)$, $(8,1)$, $(8,5)$, $(12,1)$, $(12,5)$, $(12,11)$, $(13,1)$, $(18,1)$, $(18,7)$, $(20,1)$, $(28,1)$, $(28,3)$, $(30,1)$, $(42,1)$, $(45,1)$, $(52,1)$, $(70,1)$, $(72,1)$, $(78,1)$, $(168,1)$, $(180,1)$, where $N$ and $N_1$ are as in Theorem ~{\rm \ref{the:genera}}. 
\end{cor}

This can be seen as follows. If the genus has one class represented by $f$, then its mass is
$1/|\Aut(f)|\leq1/2$. The mass formula therefore gives
\[
 \prod_{p\mid N_1}\frac{p^{v_p(N_1)-1}(p-1)}2
 \prod_{p\mid N_2}\frac{p^{v_p(N_2)-1}(p+1)}2\leq12.
\]
Every factor on the left is at least $1$, except for the factor $1/2$ when
$2\parallel N_1$. Thus every individual factor is at most $24$, so
\[
 \begin{aligned}
 p^{v_p(N_1)-1}(p-1)&\leq48, &&p\mid N_1,\\
 p^{v_p(N_2)-1}(p+1)&\leq48, &&p\mid N_2.
 \end{aligned}
\]
In particular, $p\leq47$, and the exponents are bounded as well.
Enumerating these prime powers, retaining coprime pairs satisfying
\eqref{equ1.2} and the product bound, gives $185$ candidates.
Then exact evaluation of Theorem~\ref{the:classnum} leaves precisely the
$39$ pairs displayed above. The verification program accompanying
the tables records all candidates, their masses, and their class numbers.
For these pairs, the exact divisors of $16N$ give $212$ extremal genera;
the odd exact divisors of $N/2$ with $2\parallel N$ give $43$ additional
genera of level $2N$. Hence the total is $255$. We can give a complete list of all $255$ one-class genera explicitly, 
together with explicit formulas for the representation numbers of the corresponding $255$ ternary quadratic forms, 
in Appendix~\ref{App:A} at the end of this paper.

The material of this paper is organized as follows. 
Section~\ref{Sec:2} introduces the genus conductor and proves the
classification by genus conductors. Section~\ref{Sec:3} treats orders of
level $(N_1,N_2)$, their normalizers, and their masses.
Section~\ref{Sec:Relationships} identifies the relations of ternary quadratic forms
and quaternion orders.
Section~\ref{Sec:4} defines the generalized modified Hurwitz class
numbers, and Section~\ref{Sec:proofs} proves the main formulas.
Section~\ref{Sec:examples} gives three examples of different types, and
Appendix~\ref{App:A} lists the $255$ one-class genera and their
representation formulas.

\section{Genus symbols and genus conductor}\label{Sec:2}

We recall the local invariants of primitive positive definite
ternary quadratic forms and introduce the genus conductor
through their genus symbols.

\subsection{Basics}
Let $f=(a,b,c,r,s,t)$ denote the ternary quadratic form with integer coefficients given by
\begin{equation*}
	f(x,y,z)=ax^2+by^2+cz^2+ryz+sxz+txy.
\end{equation*}
The matrix associated with $f$ is defined as
\begin{equation*}
	M_f=
	\begin{pmatrix}
		2a & t & s\\
		t & 2b & r\\
		s & r & 2c
	\end{pmatrix}.
\end{equation*}
We define the discriminant of $f$ to be
\begin{equation*}
	d_f={\det(M_f)}/{2},
\end{equation*}
and the divisor of $f$ to be
\begin{equation*}
	m_f=\gcd(M_{11},M_{22},M_{33},2M_{23},2M_{13},2M_{12}),
\end{equation*}
where $M_{ij}$ denotes the cofactor of the $(i,j)$-entry of $M_f$, i.e.\ $M_{ij}=(-1)^{i+j}A_{ij}$ with $A_{ij}$ the corresponding minor of $M_f$.
Define
\begin{equation*}
	\widetilde{m}_f=\gcd(a,b,c,r,s,t).
\end{equation*}
We say that $f$ is primitive if $\widetilde{m}_f = 1$.
Define the level of $f$ to be
\begin{equation*}
	N_f=\frac{4d_f}{m_f}.
\end{equation*}
Note that $N_f$ is the smallest positive integer for which $N_fM_f^{-1}$ is an even integral matrix. The reciprocal of primitive $f$ is defined to be the ternary quadratic form
$$
\psi_f=(M_{11}/m_f,M_{22}/m_f,M_{33}/m_f,2M_{23}/m_f,2M_{13}/m_f,2M_{12}/m_f),
$$
so that $M_{\psi_f}=\tfrac{4d_f}{m_f}M_f^{-1}=N_fM_f^{-1}$.

\vskip 0.2cm
The fundamental connection between level and discriminant is given by the following result of Lehman~\cite{Leh92}.
\begin{lemma}[\cite{Leh92}, Theorem 2]\label{the:N&d}
	Let $f$ be a primitive positive definite ternary quadratic form of level $N_f$ and discriminant $d_f$. Let
	\begin{equation}\label{eq:Nf-fact}
		N_f=2^{n_0}p_1^{n_1}\cdots p_k^{n_k}
	\end{equation}
	be the prime factorization of $N_f$. Then $n_0\geq2$ and $d_f$ has the form
	\begin{equation}\label{eq:df-fact}
		d_f=2^{d_0}p_1^{d_1}\cdots p_k^{d_k},
	\end{equation}
	with the following restrictions on the exponents $d_i$ $(\,0\leq i\leq k\,)$:
	\begin{itemize}
		\item[\rm (1)] $d_0=n_0-2$, $d_0=2n_0$, or $n_0\leq d_0\leq 2n_0-2$; and
		\item[\rm (2)] $n_i\leq d_i\leq 2n_i$ for $1\leq i\leq k$.
	\end{itemize}
	Furthermore, if $n_i$ is even for all $0\leq i\leq k$, then either $n_0\leq d_0\leq 2n_0-2$, or $d_i$ is odd for some $1\leq i\leq k$.
\end{lemma}

Two ternary quadratic forms $f_1$ and $f_2$ are said to be equivalent, written $f_1\sim f_2$, if there exists an $M\in {\rm GL}_3(\mathbb{Z})$ (the group of $3\times 3$ integer matrices with $\det(M)=\pm1$) such that $M_{f_1}=MM_{f_2}M^{t}$. Equivalent forms are said to belong to the same class. Two integral forms $f_1$ and $f_2$ are said to be semi-equivalent if they are equivalent over the $p$-adic integers (written $f_1\underset{p}{\sim} f_2$) for every prime $p$, and over the real numbers. Semi-equivalent forms are said to lie in the same genus. Equivalent forms are clearly semi-equivalent.

\vskip 0.2cm
Let $f$ be a ternary quadratic form over $\mathbb{Q}$. We say that $f$ is isotropic at $p$ if $f_p(x,y,z)=0$ for some nonzero $x,y,z\in\mathbb{Q}_p$, and anisotropic otherwise, where $f_p$ denotes the localization of $f$ at $p$. It is well known that $f$ is equivalent over $\mathbb{Q}$ to a diagonal form
\begin{equation*}
	ax^2+by^2+cz^2.
\end{equation*}
The Hasse invariant of $f$ at $p$ is then defined as
\begin{equation*}
	S_p(f)={(a,-1)}_p{(b,-1)}_p{(c,-1)}_p{(a,b)}_p{(b,c)}_p{(c,a)}_p,
\end{equation*}
which depends only on the equivalence class of $f$; here ${(a,b)}_p$ denotes the Hilbert norm residue symbol. 
The modified Hasse symbol $S^*_p$ records whether $f$ is isotropic at $p$:
\begin{equation}\label{(2.3)}
	S^*_p(f)={(-1)}^{\delta_{p,2}}S_p(f)=
	\begin{cases}
		1, & \text{if $f$ is isotropic at $p$},\\
		-1, & \text{if $f$ is anisotropic at $p$}.
	\end{cases}
\end{equation}
Here $\delta_{m,n}=1$ when $m=n$ and $0$ otherwise.
Hilbert reciprocity gives
\[
 S_\infty(f)\prod_p S_p(f)=1.
\]
For positive definite $f$ we have $S_\infty(f)=1$.
Together with \eqref{(2.3)}, this shows that the number of finite
primes at which $f$ is anisotropic is odd.
\subsection{Lehman's correspondences and Watson's transformations}\label{section2.2}
Lehman's correspondence~\cite[Theorem 4]{Leh92} and Watson's transformation~\cite{Wat62} play important roles for us.
Denote by $C(N,d)$ the set of all classes of primitive positive definite ternary quadratic forms of level $N$ and discriminant $d$. By~\cite[Theorem 4]{Leh92}, for a fixed prime $p\nmid 2Nd$ the sets $C(p^gN,p^hd)$ and $C(p^gN,p^{3g-h}d)$ are in bijection via Lehman's correspondence $\phi_p$:
\begin{align*}
	\phi_p :C(p^gN,p^hd) & \to C(p^gN,p^{3g-h}d)\\
	f=(p^{g}a,p^{h-g}b,c,p^{h-g}r,p^{g}s,p^{g}t)& \mapsto (a,p^{2g-h}b,p^gc,p^{g}r,p^gs,p^{2g-h}t),
\end{align*}
i.e.\ $M_f \mapsto PM_{f}P$, where
\begin{equation*}
	P=
	\begin{pmatrix}
		p^{-g/2} & 0& 0\\
		0 & p^{(3g-2h)/2}& 0\\
		0 & 0& p^{g/2}
	\end{pmatrix},
\end{equation*}
with $a, b, c, r, s, t$ integers and $p\nmid ac$. As for $p=2$, if $2\nmid Nd$, then the sets $C(2^{g}N,2^{h}d)$ and $C(2^{g}N,2^{3g-h-2}d)$ are in bijection via Lehman's correspondence $\phi_2$ (see~\cite[Theorem~5]{Leh92}). Moreover, for $2\nmid Nd$, by~\cite[Proposition~2.7]{LZ24} we know that the sets $C(8N,64d)$ and $C(4N,4d)$ are in bijection via Watson's transformation $\lambda_4$:
\begin{align*}
	\lambda_4 :C(8N,64d) & \to C(4N,4d)\\
	f=(8a,8b,c,8r,8s,8t)& \mapsto\lambda_4(f)=(2a,2b,c,4r,4s,2t),
\end{align*}
with $a, b, c, r, s, t$ integers and $2\nmid ac$. Notice that this will serve as a basis for the consideration of the non-boundary case of level $4N$ and discriminant $4d$. For more as for the definition of $\lambda_4$ one is referred to ~\cite[Proposition~2.7]{LZ24}.

As for the basic properties of $\phi_p$ and $\lambda_4$ we have the following lemma (see~\cite[Proposition~2.7]{LZ24}).

\begin{lemma}\label{prop:bijection}
	For any prime $p\geq 2$, we have
	\begin{itemize}
		\item[\rm (1)] $|\Aut(f)| =|\Aut(\phi_p(f))|$ for every form $f$; and
		\item[\rm (2)] $\phi_p$ is genus-preserving, i.e.\ if two classes lie in the same genus, then so do their images under $\phi_p$.
	\end{itemize}
	The same holds true for $\lambda_4$ provided that $8\parallel N_f$ and $64\parallel d_f$, or $4\parallel N_f$ and $16\parallel d_f$; furthermore, under these conditions $\lambda_4$ is bijective.
\end{lemma}

\begin{remark}\label{rem:phiN}
	(1) For a positive integer $\delta $ we define $\phi_\delta  = \phi_{p_1} \circ \cdots \circ \phi_{p_{e(\delta )}}$, where $p_i$ are the distinct prime factors of $\delta $. This is well-defined, since the correspondences $\phi_p$ for distinct primes commute with each other, and the composition operation is associative.
	(2) One can check readily that, for a form $f$ of level $4N$, we have $\phi_{4N}(f)=\psi_f$, the reciprocal of $f$.
	(3) It is clear that if $p\neq2$, then $\phi_p(f(x,y,z))=p^gf(p^{-g}x,p^{g-h}y,z)$ by definition. Hence for any prime $q\neq p$, we have $\phi_p(f)\sim_q p^g f$. This holds true for $p=2$ similarly. 
\end{remark}

\subsection{Genus symbols and genus conductor}

Recall that by ~\cite{LZ26}, for an integer $N$ of ~\eqref{equ1.2} with $e(N_0)=e(N_1)$ odd, and for a form $f$ that is anisotropic exactly at primes $p\mid N_1$, the genus of $f$ can be characterized explicitly as follows: at each place of odd primes $p$, 
\begin{equation*}
	f\underset{p}{\sim}
	\begin{cases}
		-\epsilon x^2-p^{v_p(N)}y^2+\epsilon p^{v_p(N)}z^2, & \text{if $p\mid N_1$,}\\
		-x^2-p^{v_p(N)}yz, & \text{otherwise,}
	\end{cases}
\end{equation*}
and at place $2$, 
\begin{equation*}
	f\underset{2}{\sim}
	\begin{cases}
		3x^2-2^{v_2(4N)}(y^2+z^2+yz), & \text{if $2\mid N_1$,}\\
		-x^2-2^{v_2(4N)}yz, & \text{otherwise.}
	\end{cases}
\end{equation*}
Here and in the sequel, $\epsilon$ is always used to denote a fixed quadratic non-residue modulo an odd prime $p$. Note that this genus is determined completely by the isotropy together with level and discriminant. 

However, for the general $N$ in ~\eqref{equ1.2} of this paper, we only have $e(N_0)\leq e(N_1)$ in general. And this may lead to some uncertainty, i.e., the combination of level, discriminant and isotropy together does not suffice to characterize a genus in general. To see the reason of this, we suppose that not all $v_p(N_1)$ are odd, i.e.,  $e(N_0)<e(N_1)$. Then we can take an even positive integer $u_0$ and a prime $p$ (resp.\ $2$), and consider a genus of $f_1$ satisfying
$$
f_1\underset{p}\sim -\epsilon x^2-p^{u_0}y^2+\epsilon p^{u_0}z^2\ \text{if }\ p\ne 2,\,\,\,\text{and}\quad f_1\underset{2}\sim 3x^2-2^{u_0+2}(y^2+z^2+yz);
$$
or a genus of $f_2$ satisfying
$$
f_2\underset{p}\sim -x^2-p^{u_0}yz\ \text{if }\ p\ne 2,\,\,\,\text{and}\quad f_2\underset{2}\sim -x^2-2^{u_0+2}yz.
$$
It is clear that $v_p(N_{f_1})=v_p(N_{f_2})=u_0$, $v_p(d_{f_1})=v_p(d_{f_2})=2u_0$ (resp.\ $v_2(N_{f_1})=v_2(N_{f_2})=u_0+2$, $v_2(d_{f_1})=v_2(d_{f_2})=2u_0+4$), and $f_1$ and $f_2$ are isotropic at $p$ (resp.\ at $2$).
Yet, as we will see below (under Lemma ~\ref{prop2.7}), $f_1$ and $f_2$ does not lie in the same genus. 
This observation motivates the use of the notion of genus symbol instead of isotropy in general.

\begin{definition}[\cite{Leh92}]\label{def:genussym}
	Let \(f=(a,b,c,r,s,t)\) be primitive and \(\psi\) its reciprocal. Assume, after a suitable normalization, that \(a\) and \(M_{33}/m_f\) are coprime to \(m_f m_\psi\). If \(p\) is an odd prime with \(p \mid m_f\) (resp.\ \(p \mid m_\psi\)), we set
	\[
	\left(\frac{f}{p}\right) = \left(\frac{a}{p}\right) \quad
	\left(\text{resp. } \left(\frac{\psi}{p}\right) = \left(\frac{M_{33}/m_f}{p}\right)\right);
	\]
	and for prime \(2\), if \(16 \mid m_f\) we set
	\[
	\left(\frac{f}{4}\right) = (-1)^{(a-1)/2},
	\]
	if moreover \(32 \mid m_f\) we also set
	\[
	\left(\frac{f}{8}\right) = (-1)^{(a^2-1)/8};
	\]
	the symbols \(\left(\frac{\psi}{4}\right)\) and \(\left(\frac{\psi}{8}\right)\) are defined analogously. These symbols are collectively called the genus symbols of \(f\).
\end{definition}

\begin{remark}\label{remark2.5}
	(1) For primitive $f$ of level $4N$ and discriminant $16N^2$ one can verify that $m_{\psi}=1$, thus the genus symbols are defined only for $f$ but not for $\psi$. (2) For the compatibility of the use of other notations and for convenience later, in the sequel, we will also write $\left(\frac{f}{2}\right)$ (resp. $\left(\frac{\psi}{2}\right)$) for the symbol $\left(\frac{f}{8}\right)$ (resp. $\left(\frac{\psi}{8}\right)$).
\end{remark}

\begin{remark}\label{rem:genussym}
Let $f$ be of level $4N$, $p\mid4N$, and $g=v_p(4N)$.
Then the combination of Remark~\ref{rem:phiN} and the identity
$M_{\psi_f}=4NM_f^{-1}$ gives
\[
 \psi_{\phi_p(f)}\sim_p(4N/p^g)f,
 \qquad \phi_p(f)\sim_q p^g f\quad(q\ne p).
\]
Consequently, whenever the relevant genus symbols are defined,
\[
 \left(\frac{\psi_{\phi_p(f)}}p\right)
 =\left(\frac{4N/p^g}p\right)\left(\frac fp\right),
 \quad
 \left(\frac{\phi_p(f)}q\right)
 =\left(\frac{p^g}q\right)\left(\frac fq\right).
\]
At place $2$, the divisibility requirements in
Definition~\ref{def:genussym} must be imposed; divisibility by $2$
alone does not define a genus symbol.
\end{remark}

The following proposition gives a clear characterization for a form $f$ to be anisotropic at odd primes $p$ using the notion of genus symbols.

\begin{prop}\label{prop:anisotropic}
Let $f$ be a primitive positive definite ternary quadratic form.  Let $p$ be an odd prime. 
Then $f$ is isotropic at $p$ if both $v_p(N_f)$ and $v_p(d_f)$ are even, and it is anisotropic at $p$ if and only if one of the following holds:
	\begin{itemize}
		\item[\rm (1)] $\left(\frac{f}{p}\right)\left(\frac{\psi}{p}\right)=-\left(\frac{-N_f/p^{v_p(N_f)}}{p}\right)$ when $v_p(N_f)$ is even and $v_p(d_f)$ is odd;
		\item[\rm (2)] $\left(\frac{f}{p}\right)=-\left(\frac{-d_f/p^{v_p(d_f)}}{p}\right)$ when $v_p(N_f)$ is odd and $v_p(d_f)$ is even;
		\item[\rm (3)] $\left(\frac{\psi}{p}\right)=-\left(\frac{-m_f/p^{v_p(m_f)}}{p}\right)$ when $v_p(N_f)$ is odd and $v_p(d_f)$ is odd.
	\end{itemize}
\end{prop}
\begin{proof}
Write
\[
 f\sim_p ax^2+bp^uy^2+cp^vz^2,
 \quad 0\leq u\leq v,\quad a,b,c\in\mathbb Z_p^\times.
\]
Then
\[
 v_p(m_f)=u,\quad v_p(N_f)=v,\quad v_p(d_f)=u+v,
 \quad v_p(m_\psi)=v-u.
\]
The Hilbert symbol formula gives
\[
 S_p(f)=
 \begin{cases}
  1,&\text{ if }u,v\text{ are even},\\
  \left(\dfrac{-ac}{p}\right),&\text{ if }u\text{ is odd},\ v\text{ is even},\\
  \left(\dfrac{-bc}{p}\right),&\text{ if }u,v\text{ are odd},\\
  \left(\dfrac{-ab}{p}\right),&\text{ if }u\text{ is even},\ v\text{ is odd}.
 \end{cases}
\]
Whenever the respective genus symbols are defined, we have
\[
 \left(\frac fp\right)=\left(\frac ap\right)\quad(u>0),
 \quad
 \left(\frac\psi p\right)
 =\left(\frac{cN_f/p^v}{p}\right)\quad(u<v).
\]
Also,
\[
 \left(\frac{d_f/p^{u+v}}p\right)=\left(\frac{abc}p\right),
 \quad N_fm_f=4d_f.
\]
In case \textup{(1)}, $u$ is odd and $v$ is even, so $0<u<v$
and both genus symbols are defined.
In case \textup{(2)}, $u,v$ are odd, so $u>0$ and the symbol
of $f$ is defined.
In case \textup{(3)}, $u$ is even and $v$ is odd, so $u<v$
and the symbol of $\psi$ is defined.
Substitution in the displayed expressions for $S_p(f)$ proves
the three criteria, since $f$ is anisotropic at the odd prime $p$
if and only if $S_p(f)=-1$.
\end{proof}

\begin{table}[H]
\caption{Characterization of isotropy by genus symbols}
\centering
\begin{tabular}{c|c|c|c}
\hline
 $v_p(N_f)$ & $v_p(d_f)$ & genus symbols & isotropy \\
\hline
even & even & --- & isotropic \\
even & odd & $\left(\frac{f}{p}\right)\left(\frac{\psi}{p}\right)=-\left(\frac{-N_f/p^{v_p(N_f)}}{p}\right)$ & anisotropic \\
odd & even  &  $\left(\frac{f}{p}\right)=-\left(\frac{-d_f/p^{v_p(d_f)}}{p}\right)$ & anisotropic \\
odd & odd & $\left(\frac{\psi}{p}\right)=-\left(\frac{-m_f/p^{v_p(m_f)}}{p}\right)$ & anisotropic \\
\hline
\end{tabular}
\end{table}

Using genus symbols we have the following characterization of a genus.

\begin{lemma}[\cite{Leh92}, Proposition 4]\label{prop2.7}
	Two primitive positive definite ternary quadratic forms $f_1$ and $f_2$ lie in the same genus 
	if and only if they have the same level, the same discriminant, and the same collection of genus symbols.
\end{lemma}

\begin{remark}
We note that the statement of this lemma is equivalent to that of~\cite[Chapter 15, Theorems 9 and 10]{CS88} for primitive ternary forms.
\end{remark}

We provide two examples to illustrate the significance of this lemma. Let $u,v$ be positive integers, and $p$ an odd prime. Let $f_1$ and $f_2$ be two forms such that
\begin{equation*}
	f_1\underset{p}{\sim}-\epsilon x^2-p^{u}y^2+\epsilon p^{u}z^2,
	\quad
	f_2\underset{p}{\sim}-x^2-p^{v}yz.
\end{equation*}
Then $\left(\frac{f_1}{p}\right)=\left(\frac{-\epsilon}{p}\right)=-\left(\frac{-1}{p}\right)$, whereas $\left(\frac{f_2}{p}\right)=\left(\frac{-1}{p}\right)$. This is compatible with the fact that $f_1$ and $f_2$ lie in different genera. Alternatively, if we take $f_1$ and $f_2$ such that 
\begin{equation*}
	f_1\underset{2}{\sim} 3x^2-2^{u+2}(y^2+z^2+yz),
	\quad
	f_2\underset{2}{\sim} -x^2-2^{v+2}yz,
\end{equation*}
then $\left(\frac{f_1}{8}\right)={(-1)}^{(3^2-1)/8}=-1$, whereas $\left(\frac{f_2}{8}\right)={(-1)}^{({(-1)}^2-1)/8}=1$. This is again compatible with the fact that $f_1$ and $f_2$ lie in different genera.

\vskip 0.2cm
Based on Lemma ~\ref{prop2.7} we have

\begin{prop}\label{def:genus-lemma}
	Let $N$ be a positive integer that is not a perfect square, and let $N_1$ be a positive divisor of $N$ as in~\eqref{equ1.2}. Let $G_{4N,\,16N^2,\,N_1}$ be the genus of form $f$ of level $4N$ and discriminant $16N^2$ whose collection of genus symbols satisfies
	$$
	\left\{p:\,p \mid N,\ \left(\frac{f}{p}\right)=-\left(\frac{-1}{p}\right)\right\}=\left\{p:\,p\mid N_1\right\}.
	$$
	Then, as $N_1$ runs over all positive divisors of $N$ as in~\eqref{equ1.2}, $G_{4N,\,16N^2,\,N_1}$ runs over all $2^{e(N)-1}$ genera of level $4N$ and discriminant $16N^2$. In particular, the genus $G_{4N,\,16N^2,\,N_1}$ is uniquely determined by the triple $\left(4N,\,16N^2,\,N_1\right)$. 
\end{prop}

\begin{proof}
Start with a genus of level $4N$ and discriminant $16N^2$.
Here $m_f=16N$ and $m_\psi=1$. The genus symbol
$\left(\frac f4\right)=-1$ is fixed by the level and discriminant;
the remaining symbols are those at primes dividing $N$.
They determine a unique exact divisor $N_1$ by the condition in the
statement. The local isotropy criteria and Hilbert reciprocity imply
that $N_1=s^2N_0$ with $\mu(N_0)=-1$.
Lemma~\ref{prop2.7} shows that this assignment from genera to admissible
divisors is injective.

By~\cite[Lemma~3]{Leh92}, the number of these genera is $2^{e(N)-1}$.
There are equally many admissible divisors: since $N$ is not a square,
at least one $v_p(N)$ is odd, and toggling the full $p$-part of an exact
divisor reverses the parity of $e(N_0)$.
Therefore exactly half of the $2^{e(N)}$ exact divisors are admissible.
The injection is consequently a bijection, proving both existence
and uniqueness.
\end{proof}

The above result justifies the following definitions.

\begin{definition}\label{def:genus}
	Let $G_{4N,\,16N^2,\,N_1}$ be the uniquely determined genus of level $4N$, 
	discriminant $16N^2$ and the associated positive divisor $N_1$ of $N$ as in Proposition~\ref{def:genus-lemma}.
	 Then we call the third parameter $N_1$ the conductor of genus $G_{4N,\,16N^2,\,N_1}$.
\end{definition}

\begin{definition}\label{def:genus111}
	For positive integers $\delta $ with $\delta \parallel 16N$, let $G_{4N,\,16N^2/\delta ,\,N_1}$
	 be the genus obtained from the base genus $G_{4N,\,16N^2,\,N_1}$ by applying the correspondence $\phi_\delta $
	 (cf.\ Remark~\ref{rem:phiN}). Moreover, if $2\parallel N$, let $G_{2N,\,N^2/\delta ',\,N_1}$
	  be the genus obtained from $G_{4N,\,16N^2,\,N_1}$ by applying $\phi_{\delta '}\circ\lambda_4$,
	   for a positive integer $\delta '$ with $\delta '\parallel N/2$. In each case $N_1$
	    is called the conductor of genus $G_{4N,\,16N^2/\delta ,\,N_1}$ (resp.\ $G_{2N,\,N^2/\delta ',\,N_1}$).
\end{definition}

\begin{prop}\label{jiajia}
Let $\delta\parallel16N$ and $f\in G_{4N,16N^2,N_1}$.
For each prime $p\mid N$, the conductor $N_1$ of
$G_{4N,16N^2/\delta,N_1}$ is characterized by the genus symbols of $\phi_\delta(f)$ as follows:
\[
 \begin{aligned}
 \left(\frac{\psi_{\phi_\delta(f)}}p\right)
   &=-\left(\frac{-16N/\delta}p\right)
     \quad\Longleftrightarrow\quad p\mid N_1
       &&(p\mid\delta),\\
 \left(\frac{\phi_\delta(f)}p\right)
   &=-\left(\frac{-\delta}p\right)
     \quad\Longleftrightarrow\quad p\mid N_1
       &&(p\nmid\delta).
 \end{aligned}
\]
All genus symbols in these two cases are defined. In particular, the
genus symbol for $p=2$ is used here only when $2\mid N$.
\end{prop}
\begin{proof}
At a prime $p\mid N$, the local similarity formulas in
Remark~\ref{rem:genussym} give
\[
 \begin{aligned}
 \left(\frac{\psi_{\phi_\delta(f)}}p\right)
 &=\left(\frac{16N/\delta}p\right)\left(\frac fp\right)
       &&(p\mid\delta),\\
 \left(\frac{\phi_\delta(f)}p\right)
 &=\left(\frac\delta p\right)\left(\frac fp\right)
       &&(p\nmid\delta).
 \end{aligned}
\]
The factors on the right are nonzero.
The claim follows from Proposition~\ref{def:genus-lemma}.
\end{proof}

\begin{remark}\label{remark2.13}
Suppose $2\parallel N$. For
$f\in G_{2N,\,N^2/\delta',\,N_1}$ with $\delta'\parallel N/2$,
Watson's transformation gives
$\lambda_4(f)\sim_p f/4$ at every odd prime $p$.
Since its level is $2N$, reciprocity gives
\[
 \psi_{\lambda_4(f)}\sim_p2\psi_f.
\]
Thus, whenever the symbols are defined,
\[
 \left(\frac{\lambda_4(f)}p\right)=\left(\frac fp\right),
 \quad
 \left(\frac{\psi_{\lambda_4(f)}}p\right)
 =\left(\frac2p\right)\left(\frac{\psi_f}p\right).
\]
These identities together with Proposition~\ref{jiajia} determine the odd part
of the conductor of $G_{2N,\,N^2/\delta',\,N_1}$.
At place $2$, both divisors have valuation $2$, so neither symbol with
denominator $8$ is defined. The form is anisotropic at $2$ precisely
when $2\mid N_1$, which determines the $2$ part of the conductor.
\end{remark}

\begin{remark}\label{rem:genus}
	One can check that $f$ is anisotropic at $p$ if and only if $v_p(N_1)$ is odd. Since the number of primes at which $f$ 
	is anisotropic is always odd, we must have $N_1 = s^2 N_0$ with $\mu(N_0)=-1$. This is compatible with the assumption on $N_0$ given in ~\eqref{equ1.2}.
\end{remark}

	\section{Quaternion orders and type numbers}\label{Sec:3}
	
	In this section we recall some fundamental facts about quaternion algebras, generalize the notion of level 
	attached to an order by giving explicit characterizations of some new kinds of quaternion orders, and establish general type-number and mass formulas for the general $N_1$ and $N_2$ of~\eqref{equ1.2}.
	
	\subsection{Basics}
	For a field $F$ of characteristic $0$ and $a,b\in F^{\times}$, we write $\left(\frac{a,\,b}{F}\right)$ for the $F$-algebra with basis $1,i,j,k=ij=-ji$ 
	and defining relations $i^2=a$, $j^2=b$. In this paper $F$ is always $\mathbb{Q}$, $\mathbb{Q}_\infty=\mathbb{R}$, or a field $\mathbb{Q}_p$ of $p$-adic numbers for a prime $p$.
	 Over $\mathbb{Q}_p$ or $\mathbb{Q}_{\infty}$ there are exactly two isomorphism classes of quaternion algebras: the class of the matrix algebra $M_2(\mathbb{Q}_p)$, 
	 and the class of the unique division algebra (skew field). One has $Q_p=\left(\frac{a,\,b}{\mathbb{Q}_p}\right)\simeq M_2(\mathbb{Q}_p)$ if and only if ${(a,b)}_p=1$. 
	 Over $\mathbb{Q}$, two quaternion algebras $Q=\left(\frac{a,\,b}{\mathbb{Q}}\right)$ and $Q'=\left(\frac{a',b'}{\mathbb{Q}}\right)$ 
	 are isomorphic if and only if ${(a,b)}_p={(a',b')}_p$ for all $p$ including $\infty$. 
	 A quaternion algebra $Q$ over $\mathbb{Q}$ is said to be ramified at $p$ if $\mathbb{Q}_p\otimes Q$ is a division algebra, 
	 and split at $p$ if $\mathbb{Q}_p\otimes Q\simeq M_2(\mathbb{Q}_p)$; if $Q$ ramifies at $\infty$, it is called definite.
	  A quaternion algebra ramifies at only finitely many places, and the number of these is always even. 
	  Conversely, for any set $S$ of an even number of primes, one can find $a,b\in\mathbb{Q}^{\times}$ with ${(a,b)}_p=-1$ exactly for $p\in S$; 
	  the resulting $Q=\left(\frac{a,\,b}{\mathbb{Q}}\right)$ is then, up to isomorphism, the unique quaternion algebra ramified exactly at the primes in $S$.
		
	For an element $\alpha = t + xi + yj + zk \in Q$, its conjugation is defined by the involution
	\[
	\overline{\alpha} = t - xi - yj - zk.
	\]
	This involution satisfies $\overline{\overline{\alpha}} = \alpha$ and $\overline{\alpha\beta} = \overline{\beta} \, \overline{\alpha}$ for all $\alpha, \beta \in Q$. The reduced trace $\tr$ and the reduced norm $\n$ are defined by
	\[
	\tr(\alpha) = \alpha + \overline{\alpha}, \quad \n(\alpha) = \alpha \cdot \overline{\alpha}.
	\]
		Every element $\alpha \in Q$ satisfies its reduced characteristic equation
	\[
	x^2 - \tr(\alpha) x + \n(\alpha) = 0.
	\]
When $Q \simeq M_2(\mathbb{Q})$ we write
	\[
	\alpha = \begin{pmatrix} a & b \\ c & d \end{pmatrix} \in Q,
	\]
	and we have
	\[
	\overline{\alpha} = \begin{pmatrix} d & -b \\ -c & a \end{pmatrix}, \quad \tr(\alpha) = a + d, \quad \n(\alpha) = ad - bc = \det(\alpha).
	\]

	\subsection{Orders of level $(N_1,N_2)$}
	An ideal of a quaternion algebra $Q$ is a full $\mathbb{Z}$-lattice in $Q$, that is, a finitely generated $\mathbb{Z}$-submodule of $Q$ that contains a $\mathbb{Q}$-basis of $Q$. 
	A (quaternion) order $\mathcal{O}$ in $Q$ is an ideal that is also a ring containing $\mathbb{Z}$. 
	Throughout this paper we always use $\mathcal{O}$ to denote orders in $Q$, and $\mathcal{O}_p$ is used to denote the localization of $\mathcal{O}$ over $\mathbb{Z}_p$ for prime $p\geq 2$. 
	It is known that $\tr(\alpha),\n(\alpha)\in \mathbb{Z}$ for $\alpha\in\mathcal{O}$. 
	Let $I$ be the ideal generated by all $\det(\tr(\alpha_i\bar\alpha_j))$, where $\alpha_1,\ldots,\alpha_4\in \mathcal{O}$. 
	Then $I$ is the square of an ideal. 
	The reduced discriminant, $\discrd(\mathcal{O})$, of $\mathcal{O}$ is defined to be the positive generator of the square root of $I$. 	
	
	An order is maximal if it is not properly contained in any other orders. 
	A quaternion division algebra $Q_p$ over $\mathbb{Q}_p$ contains a unique maximal order, 
	namely $R_p+R_p j$ with $R_p=\mathbb{Z}_p+i\mathbb{Z}_p$ for odd $p$ and $R_2=\mathbb{Z}_2+\frac{i+1}{2}\mathbb{Z}_2$ for $p=2$, 
	where $j^2=p$, $ji=-ij$, and $i^2=\epsilon$ with $\left(\frac{\epsilon}{p}\right)=-1$ for odd $p$, while $i^2=-3$ for $p=2$~\cite[Chapter~13]{Voi21}.
 Writing
	\begin{equation*}
		Q_p\simeq\left\{
		\begin{pmatrix}
			\alpha & \beta\\
			p\bar{\beta} & \bar{\alpha}
		\end{pmatrix}
		: \alpha,\beta\in \mathbb{Q}_p+i\mathbb{Q}_p\right\},
	\end{equation*}
	then the unique maximal order of $Q_p$ is
	\begin{equation*}
		\mathcal{O}_p^{(1)}\simeq\left\{
		\begin{pmatrix}
			\alpha & \beta\\
			p\bar{\beta} & \bar{\alpha}
		\end{pmatrix}
		: \alpha,\beta\in R_p\right\}.
	\end{equation*}
	
A straightforward computation yields the following lemma.
	\begin{lemma}[\cite{Piz76b}, Definition 1]\label{lemma3.1}
		Let $p\geq 2$ be a prime. Suppose $Q$ is ramified at $p$. Set for any nonnegative integer $u$,  
		\begin{equation*}
			\mathcal{O}_p^{(2u+1)}=\left\{
			\begin{pmatrix}
				\alpha & p^u\beta\\
				p^{u+1}\bar{\beta} & \bar{\alpha}
			\end{pmatrix}
			: \alpha,\beta\in R_p\right\}.
		\end{equation*}
		Then $\mathcal{O}_p^{(2u+1)}$ is an order in $Q_p$, and $[\mathcal{O}_p^{(1)}:\mathcal{O}_p^{(2u+1)}]=p^{2u}$.
	\end{lemma}
	
	If, alternatively, $Q$ splits at $p$, then the maximal orders in ${Q}_p$ are isomorphic to $M_2(\mathbb{Z}_p)$ over $\mathbb{Z}_p$,
	 and the Eichler orders $\mathcal{O}$ of level $p^n$ in $Q$ are those isomorphic to
	$\begin{pmatrix}
		\mathbb{Z}_p & \mathbb{Z}_p\\
		p^n\mathbb{Z}_p & \mathbb{Z}_p
	\end{pmatrix}$ over $\mathbb{Z}_p$~\cite[Chapter 23.4]{Voi21}.
	
	When $n>0$ is even, we can give a new explicit description for Eichler orders of level $p^n$ in the following Lemma \ref{lemma3.2}. We point out that this description may have potential use for the classification of Eichler orders later, as well as to be a motivation of Lemma \ref{lemma3.3} below.
	
	\begin{lemma}\label{lemma3.2}
		Let $\mathcal{O}$ be an Eichler order of level $p^{2u}$, where $2u$ is a positive even integer. 
Let $\mathcal{O}_p$ denote the localization of $\mathcal{O}$ over $\mathbb{Z}_p$. Then, for $p\neq2$, $\mathcal{O}_p$ is conjugate to
		\begin{equation*}
			\left\{ M\in M_2(\mathbb{Z}_p):M\equiv
			\begin{pmatrix}
				a & b\\
				b & a
			\end{pmatrix}\!\!\!\!
			\pmod{p^{u}},\,\,a, b\in\mathbb{Z}_p\right\},
		\end{equation*}
		and for prime $2$, $\mathcal{O}_2$ is conjugate to
		\begin{equation*}
			\left\{ M\in M_2(\mathbb{Z}_2):M\equiv
			\begin{pmatrix}
				a+b & b\\
				0 & a
			\end{pmatrix}
			\!\!\!\!
			\pmod{2^{u}},\,\,a, b\in\mathbb{Z}_2\right\}.
		\end{equation*}
	\end{lemma}
	
	\begin{proof}
		The lemma follows from the following isomorphic formulations of $\mathcal{O}_p$ and $\mathcal{O}_2$ respectively: 
		\begin{equation*}
			\begin{pmatrix}
				p^{u} & 1\\
				-p^{u} & 1
			\end{pmatrix}
			\begin{pmatrix}
				\mathbb{Z}_p & \mathbb{Z}_p \\
				p^{2u}\mathbb{Z}_p  & \mathbb{Z}_p
			\end{pmatrix}
			\begin{pmatrix}
				p^{u} & 1\\
				-p^{u} & 1
			\end{pmatrix}^{-1},\quad 
            \begin{pmatrix}
				-2^{u} & 1\\
				2^{u} & 0
			\end{pmatrix}
			\begin{pmatrix}
				\mathbb{Z}_2 & \mathbb{Z}_2 \\
				2^{2u}\mathbb{Z}_2  & \mathbb{Z}_2
			\end{pmatrix}
			\begin{pmatrix}
				-2^{u} & 1\\
				2^{u} & 0
			\end{pmatrix}^{-1}.
		\end{equation*}
		\end{proof}
	
The next lemma gives explicit matrix models for split residually
inert orders that occur when $v_p(N_1)$ is positive and even.
These local orders extend the family used in \cite{LZ26}.

	\begin{lemma}\label{lemma3.3}
		Let $2u$ be positive and even. For a prime $p\neq2$ set
		\begin{equation*}
			\mathcal{O}_p^{(2u)}=\left\{ M\in M_2(\mathbb{Z}_p):M\equiv
			\begin{pmatrix}
				a & \epsilon b\\
				b & a
			\end{pmatrix}
			\!\!\!\!
			\pmod{p^{u}},\,\,a, b\in\mathbb{Z}_p\right\},
		\end{equation*}
		and for prime $2$ set
		\begin{equation*}
			\mathcal{O}_2^{(2u)}=\left\{ M\in M_2(\mathbb{Z}_2):M\equiv
			\begin{pmatrix}
				a+b & b\\
				b & a
			\end{pmatrix}
			\!\!\!\!
			\pmod{2^{u}},\,\,a, b\in\mathbb{Z}_2\right\}.
		\end{equation*}
		Then for every $p\geq 2$, $\mathcal{O}_p^{(2u)}$ is an order in $M_2(\mathbb{Q}_p)$ and $[M_2(\mathbb{Z}_p):\mathcal{O}_p^{(2u)}]=p^{2u}$.
	\end{lemma}
	
	\begin{proof}
		It is easy to verify that $\mathcal{O}_p^{(2u)}$ is a ring containing $\mathbb{Z}$. Using the identities $a+k_1p^{u}=(a+\frac{(k_1+k_2)p^{u}}{2})+\frac{k_1-k_2}{2}p^{u}$ and $a+k_2p^{u}=(a+\frac{(k_1+k_2)p^{u}}{2})-\frac{k_1-k_2}{2}p^{u}$ for $k_1,k_2\in\mathbb{Z}_p$, we obtain
		\begin{equation*}
			\mathcal{O}_p^{(2u)}=
			\begin{pmatrix}
				1 & 0\\
				0 & 1
			\end{pmatrix}\mathbb{Z}_p+\begin{pmatrix}
				p^u & 0\\
				0 & -p^u
			\end{pmatrix}\mathbb{Z}_p+\begin{pmatrix}
				0 & \epsilon\\
				1 & 0
			\end{pmatrix}\mathbb{Z}_p+\begin{pmatrix}
				0 & \epsilon p^u\\
				-p^u & 0
			\end{pmatrix}\mathbb{Z}_p,
		\end{equation*}
		and, for the prime $2$,
		\begin{equation*}
			\mathcal{O}_2^{(2u)}=
			\begin{pmatrix}
				1 & 0\\
				0 & 1
			\end{pmatrix}\mathbb{Z}_2+\begin{pmatrix}
				1 & 1\\
				1 & 0
			\end{pmatrix}\mathbb{Z}_2+\begin{pmatrix}
				0 & 2^{u}\\
				0 & 0
			\end{pmatrix}\mathbb{Z}_2+\begin{pmatrix}
				0 & 0\\
				2^{u} & 0
			\end{pmatrix}\mathbb{Z}_2.
		\end{equation*}
		The lemma now follows by direct verification.
	\end{proof}
	
	\begin{remark}\label{rem:Eichlerorder}
		The orders defined in Lemma~\ref{lemma3.3} are not intersections of maximal orders (see~\cite[Chapter 3.1.3]{Tu12}), i.e.\ they are not Eichler orders, since the set
		\begin{equation*}
			\{\alpha\in\mathcal{O}_p^{(2u)}:\det(\alpha)\equiv0\!\!\!\!\pmod{p}\}=\mathcal{O}_p^{(2u)}\cap pM_2(\mathbb{Z}_p)
		\end{equation*}
		is closed under addition, whereas Eichler orders lack this property.
	\end{remark}
	
	With Lemmas~\ref{lemma3.1}, \ref{lemma3.2}, and~\ref{lemma3.3} in place, we can now extend the notion of level for orders to all admissible cases.
	 This extension is new and crucial for our discussion in the sequel. 
	
	\begin{definition}\label{def3.6}
		Let $N_0$, $N_1$ and $N_2$ be as in ~\eqref{equ1.2}. Let $Q$ be a definite quaternion algebra ramified exactly at finite primes $p\mid N_0$, and let $\mathcal{O}$ be an order in $Q$. If
		\begin{enumerate}
			\item for $p^{2u+1}\parallel N_1$, $\mathcal{O}_{p}$ is isomorphic to $\left\{
			\begin{pmatrix}
				\alpha & p^u\beta\\
				p^{u+1}\bar{\beta} & \bar{\alpha}
			\end{pmatrix}
			: \alpha,\beta\in R_p\right\}$ over $\mathbb{Z}_{p}$,
			
			\item for $p^{2u}\parallel N_1$ with $p\neq 2$, $\mathcal{O}_{p}$ is isomorphic to $\left\{ M\in M_2(\mathbb{Z}_p):M\equiv
			\begin{pmatrix}
				a & \epsilon b\\
				b & a
			\end{pmatrix}
			\pmod{p^u}\right\}$ over $\mathbb{Z}_{p}$,
			
			\item for $2^{2u}\parallel N_1$, $\mathcal{O}_{2}$ is isomorphic to $\left\{ M\in M_2(\mathbb{Z}_2):M\equiv
			\begin{pmatrix}
				a+b & b\\
				b & a
			\end{pmatrix}
			\pmod{2^u}\right\}$ over $\mathbb{Z}_{2}$, and
			
			\item for $p\nmid N_1$, $\mathcal{O}_{p}$ is isomorphic to $\begin{pmatrix}
				\mathbb{Z}_p & \mathbb{Z}_p\\
				N_2\mathbb{Z}_p & \mathbb{Z}_p
			\end{pmatrix}$ over $\mathbb{Z}_p$,
		\end{enumerate}
		then we say that $\mathcal{O}$ is an order of level $(N_1,N_2)$. When $N_2=1$, we simply say that $\mathcal{O}$ is an order of level $N_1$.
	\end{definition}

\subsection{Relations with Eichler symbols}
The generalization of level $(N_1,N_2)$ of orders $\mathcal{O}$ in Definition \ref{def3.6} has a close relation to Eichler symbols.
 To see this we first recall some basic relevant concepts; for more details one is referred to~\cite[Chapter 24.3]{Voi21}.
Let  $\operatorname{rad}(\mathcal{O})$ be the Jacobson radical of order $\mathcal{O}$, i.e., the intersection of all maximal left ideals of $\mathcal{O}$. 
Then the $\mathbb{F}_p$-algebra $\mathcal{O}_p/\operatorname{rad}(\mathcal{O}_p)$ is semisimple and has a standard involution, 
where $\mathbb{F}_p$ is the finite field $\{0,1,\ldots,p-1\}$. Thus there exist three
possibilities for $\mathcal{O}_p/\operatorname{rad}(\mathcal{O}_p)$ up to isomorphisms: it is either $\mathbb{F}_p$,
 or a separable quadratic $\mathbb{F}_p$-algebra, or a quaternion algebra over $\mathbb{F}_p$.
\begin{definition}
Let $p\geq 2$ be a prime, then the Eichler symbol $\left(\frac{\mathcal{O}}{p}\right)$ is defined as
\begin{equation*}
\left(\frac{\mathcal{O}}{p}\right)=
\begin{cases}
*, &\text{if $\mathcal{O}_p/\operatorname{rad}(\mathcal{O}_p)$ is a quaternion algebra over $\mathbb{F}_p$;}\\
1, &\text{if $\mathcal{O}_p/\operatorname{rad}(\mathcal{O}_p)\simeq \mathbb{F}_p\times\mathbb{F}_p$; in this case $\mathcal{O}_p$ is said residually split;}\\
0, &\text{if $\mathcal{O}_p/\operatorname{rad}(\mathcal{O}_p)\simeq \mathbb{F}_p$; in this case $\mathcal{O}_p$ is said residually ramified;}\\
-1, &\text{if $\mathcal{O}_p/\operatorname{rad}(\mathcal{O}_p)$ is a quadratic field extension of $\mathbb{F}_p$; in this case $\mathcal{O}_p$ is said residually inert.}
\end{cases}
\end{equation*}
\end{definition}
\noindent Here $*$ is an indication for the split maximal case, whose semisimple quotient is $M_2(\mathbb F_p)$. 
Eichler symbol is in fact determined by the Kronecker symbol of the discriminant quadratic form $\tr(\alpha)^2-4\n(\alpha)$.
More accurately, we have~\cite[24.3.3]{Voi21}
\begin{itemize}
\item $\left(\frac{\mathcal{O}}{p}\right)=*$ if and only if $\left\{\left(\frac{\tr(\alpha)^2-4\n(\alpha)}{p}\right):\,\alpha\in\mathcal{O}\right\}=\{-1, 0, 1\}$;
\item $\left(\frac{\mathcal{O}}{p}\right)=1$ (resp. $-1$ or $0$) if and only if $\left\{\left(\frac{\tr(\alpha)^2-4\n(\alpha)}{p}\right):\,\alpha\in\mathcal{O}\right\}=\{0, 1\}$ (resp. $\{0,-1\}$ or $\{0\}$).
\end{itemize}

\begin{prop}\label{prop:Eichlersym}
Let $O$ be an order of level $(N_1,N_2)$. Then $O_p$ is residually inert for $p\mid N_1$, residually split for $p\mid N_2$, 
and $O_p/\operatorname{rad}(O_p)$ is a quaternion algebra over $\mathbb F_p$ for $p\nmid N_1N_2$.
\end{prop}

\begin{proof}
We use the local orders in Definition \ref{def3.6}. Suppose first that $p^{2u+1}\parallel N_1$. Write
\[
 x=\begin{pmatrix}\alpha&p^u\beta\\p^{u+1}\bar\beta&\bar\alpha\end{pmatrix}\in O_p,
 \quad \alpha,\beta\in R_p.
\]
The reduction of this order modulo its Jacobson radical is $R_p/pR_p\simeq\mathbb F_{p^2}$: the map sending $x$ to $\alpha\bmod pR_p$ is a surjective ring homomorphism, 
and its kernel is the radical. Thus $O_p$ is residually inert. For odd $p$, this can also be seen from the discriminant calculation
\[
 \operatorname{tr}(x)^2-4\operatorname n(x)
 = (\alpha-\bar\alpha)^2+4p^{2u+1}\beta\bar\beta,
\]
whose quadratic character values modulo $p$ are precisely $\{0,-1\}$.

Next suppose $p^{2u}\parallel N_1$ with $u\geq1$. If $p$ is odd, reduction of the matrices
\[
 x\equiv\begin{pmatrix}a&\epsilon b\\ b&a\end{pmatrix}\pmod{p^u}
\]
identifies the semisimple quotient with $\mathbb F_p[t]/(t^2-\epsilon)\simeq\mathbb F_{p^2}$, since $\epsilon$ is a nonsquare modulo $p$. Equivalently,
\[
 \operatorname{tr}(x)^2-4\operatorname n(x)\equiv4\epsilon b^2\pmod p,
\]
so its quadratic character values are $\{0,-1\}$. For $p=2$, the matrices modulo $2$ have the form
\[
 \begin{pmatrix}a+b&b\\b&a\end{pmatrix}
 =aI+bT,\quad T=\begin{pmatrix}1&1\\1&0\end{pmatrix},\quad T^2+T+I=0.
\]
Since $t^2+t+1$ is irreducible over $\mathbb F_2$, the semisimple quotient is $\mathbb F_4$. Hence $O_p$ is residually inert in every case $p\mid N_1$.

If $p\mid N_2$, the local order is an Eichler order of positive level; its reduction modulo $p$ is triangular, with semisimple quotient $\mathbb F_p\times\mathbb F_p$. Thus $O_p$ is residually split. Finally, if $p\nmid N_1N_2$, then $O_p\simeq M_2(\mathbb Z_p)$ and
\[
 O_p/\operatorname{rad}(O_p)\simeq M_2(\mathbb F_p),
\]
which is a split quaternion algebra over $\mathbb F_p$.
\end{proof}

\subsection{Type numbers}
Two orders in $Q$ are of the same type if they are conjugate by an element of $Q^\times$, or equivalently, if they are isomorphic as
$\mathbb Z$-algebras. The genus of an order $\mathcal O$ consists of the orders in $Q$ that are locally isomorphic to $\mathcal O$ at every
finite prime. Its type set is the set of isomorphism classes in this genus. We write $T_{N_1,\,N_2}$ for the cardinality when $\mathcal O$
has level $(N_1,N_2)$.

We will use the normalizer-counting method of \cite[Section~2]{LZ26}.
In that paper every prime divisor of $N_1$ has odd valuation.
The additional local orders in the present paper occur at
$p^{2u}\parallel N_1$, where $u\geq1$ and $Q_p$ is split.
We first determine their normalizers. The conversion of the resulting
count to representation numbers is carried out in
Section~\ref{sec:6.3}.

Let $N=N_1N_2$ and define
\[
 N(\mathcal O_p)\!:\,=\{x\in Q_p^\times:x^{-1}\mathcal O_px=\mathcal O_p\},
 \quad
 N_Q(\mathcal O)\!:\,=\{x\in Q^\times:x^{-1}\mathcal Ox=\mathcal O\}.
\]
The matrix representatives below refer to the local models in
Definition~\ref{def3.6}.

\begin{lemma}\label{lemmafornormalizer}
If $p\nmid N$, then
$N(\mathcal O_p)=\mathbb Q_p^\times\mathcal O_p^\times$.
If $p\mid N$, then
\[
 N(\mathcal O_p)=\mathbb Q_p^\times\mathcal O_p^\times
 \cup\mathbb Q_p^\times w_p\mathcal O_p^\times,
\]
where 
\[
w_p=
\begin{cases}
 \begin{pmatrix}0&p^u\\p^{u+1}&0\end{pmatrix},
     &\text{if }\,\, p^{2u+1}\parallel N_1,\quad u\geq0,\\[6pt]
 p^u\begin{pmatrix}0&\epsilon_p\\-1&0\end{pmatrix},
     &\text{if }\,\, p^{2u}\parallel N_1,\quad p\ne2,\quad u\geq1,\\[6pt]
 2^u\begin{pmatrix}0&1\\-1&0\end{pmatrix},
     &\text{if }\,\, 2^{2u}\parallel N_1,\quad u\geq1,\\[6pt]
 \begin{pmatrix}0&1\\p^v&0\end{pmatrix},
     &\text{if }\,\, p^v\parallel N_2,\quad v\geq1.
\end{cases}
\]
Here $\epsilon_p$ is a quadratic nonresidue in $\mathbb Z_p^\times$.
In each case $w_p\in\mathcal O_p$ and
$v_p(\n(w_p))=v_p(N)$.
\end{lemma}

\begin{proof}
The maximal, division, and Eichler cases follow from
\cite[Proposition~2.18]{Piz77}.
Suppose $p^{2u}\parallel N_1$ with $u\geq1$, $p\geq 2$, and put
\[
 T=
 \begin{cases}
 \begin{pmatrix}0&\epsilon_p\\1&0\end{pmatrix},&\text{if }\,\,p\ne2,\\[5pt]
 \begin{pmatrix}1&1\\1&0\end{pmatrix},&\text{if }\,\,p=2.
 \end{cases}
\]
Then $\mathcal O_p=\mathbb Z_p[T]+p^uM_2(\mathbb Z_p)$.
The matrix $p^{-u}w_p$ has unit determinant and conjugates $T$
to $-T$ if $p\ne2$, and to $I-T$ if $p=2$.
Thus the asserted union is contained in the normalizer.

Conversely, after scalar multiplication we may take
$g\in N(\mathcal O_p)$ with
$g\in M_2(\mathbb Z_p)\setminus pM_2(\mathbb Z_p)$.
Write
\[
 gTg^{-1}\equiv aI+bT\!\!\pmod{p^u}.
\]
For odd $p$, comparison of traces and determinants gives
$a\equiv0$ and $b\equiv\pm1\pmod{p^u}$.
For $p=2$, it gives
\[
 2a+b\equiv1,\quad a^2+ab-b^2\equiv-1\!\!\pmod{2^u},
\]
and hence $5a(a-1)\equiv0\pmod{2^u}$.
Thus $(a,b)\equiv(0,1)$ or $(1,-1)\pmod{2^u}$.
After replacing $g$ by $(p^{-u}w_p)^{-1}g$ if necessary,
we therefore have $gT\equiv Tg\pmod{p^u}$.
Comparison of the matrix entries gives $g\in\mathcal O_p$.
Its reduction modulo $p$ is a nonzero element of the field
generated by $T$ over $\mathbb F_p$, so its determinant is a unit.
Hence $g\in\mathcal O_p^\times$, proving the reverse inclusion.
The two cosets are distinct because $p^{-u}w_p$ has unit
determinant and does not belong to $\mathcal O_p$.
\end{proof}

A nonzero element $x\in\mathcal O$ is called primitive if
$x\notin p\mathcal O$ for every prime $p$. Define
\begin{equation}\label{rho'}
 \rho'_{\mathcal O}(n,r)\!:\,=
 \#\{x\in\mathcal O:x\text{ is primitive},\ x\in N_Q(\mathcal O),
              \ \n(x)=n,\ \tr(x)=r\}.
\end{equation}
Each class in $N_Q(\mathcal O)/\mathbb Q^\times$ has exactly two
primitive integral representatives, differing by sign.
By Skolem--Noether theorem this quotient is $\Aut(\mathcal O)$.
Lemma~\ref{lemmaforprincipalideal} below will show that a primitive
normalizing element satisfies $n\parallel N$ and $n\mid r$.
Consequently, for representatives
$\mathcal O_1,\ldots,\mathcal O_{T_{N_1,\,N_2}}$ of the type set,
\begin{equation}\label{eq:type-normalizer-count}
 T_{N_1,\,N_2}=
 \frac12\sum_{n\parallel N}
 \sum_{\substack{n\mid r\\r^2\leq4n}}
 \sum_{\mu=1}^{T_{N_1,\,N_2}}
 \frac{\rho'_{\mathcal O_\mu}(n,r)}{|\Aut(\mathcal O_\mu)|}.
\end{equation}
This is the counting identity of \cite[Theorem~2.7]{LZ26} in the
present setting. Its proof uses the normalizer description and
primitive representatives, rather than a restriction on the parity
of the exponents in $N_1$.

	\subsection{The mass formula}
	
Let $I_1,\ldots,I_h$ represent the classes of locally principal
fractional left $\mathcal O$-ideals, and put
$\mathcal O_i=\mathcal O_{\operatorname R}(I_i)$.
Here $h=h_{N_1,\,N_2}$, and ideal classes are taken modulo right
multiplication by $Q^\times$. Choose the indexing so that
$\mathcal O_1,\ldots,\mathcal O_{T_{N_1,\,N_2}}$ represent the distinct types.

\begin{prop}\label{prop:mass3}
Writing $N=N_1N_2$ and $e=e(N)$, we have
\[
 \sum_{i=1}^{h}\frac1{|\mathcal O_i^\times|}
 =\frac N{24}\prod_{p\mid N_1}\left(1-\frac1p\right)
                  \prod_{p\mid N_2}\left(1+\frac1p\right),
\]
and
\begin{equation}\label{gongshi3.2}
 \sum_{\mu=1}^{T_{N_1,\,N_2}}\frac1{|\Aut(\mathcal O_\mu)|}
 =2^{1-e}\sum_{i=1}^{h}\frac1{|\mathcal O_i^\times|}
 =\frac{2^{-e}N}{12}
   \prod_{p\mid N_1}\left(1-\frac1p\right)
   \prod_{p\mid N_2}\left(1+\frac1p\right).
\end{equation}
\end{prop}

\begin{proof}
We first compute the ideal-class mass, using the mass comparison
formula \cite[Proposition~26.6.4]{Voi21}.
At a new prime $p^{2u}\parallel N_1$, write
$\mathcal O_p=\mathbb Z_p[T]+p^uM_2(\mathbb Z_p)$.
Reduction modulo $p^u$ gives, for every prime $p$, including $2$,
\[
 [M_2(\mathbb Z_p)^\times:\mathcal O_p^\times]
 =\frac{|\operatorname{GL}_2(\mathbb Z/p^u\mathbb Z)|}
             {|(\mathbb Z_p[T]/p^u\mathbb Z_p[T])^\times|}
 =\frac{(p^2-1)(p^2-p)p^{4u-4}}{(p^2-1)p^{2u-2}}=p^{2u-1}(p-1).
\]
At $p^{2u+1}\parallel N_1$, the index in the unit group of the
maximal division order is $p^{2u}$; at $p^v\parallel N_2$, the
corresponding split index is $p^{v-1}(p+1)$.
Multiplying these indices by the maximal-order mass
$\frac1{24}\prod_{p\mid N_0}(p-1)$ proves the first formula.

For \eqref{gongshi3.2}, we need to explain the passage from ideal classes to
types. Let $m_\mu$ be the size of the fiber above
$[\mathcal O_\mu]$ of the map $[I]\mapsto[\mathcal O_{\operatorname R}(I)]$.
The group of locally principal two-sided fractional
$\mathcal O_\mu$-ideals modulo rational scalars has order $2^e$,
by Lemma~\ref{lemmafornormalizer}. Its orbit on ideal classes has
size $m_\mu$, and its stabilizer is the group of principal
two-sided ideals modulo rational scalars. Together with the exact
sequence
\[
 1\longrightarrow\mathcal O_\mu^\times/\{\pm1\}
 \longrightarrow\Aut(\mathcal O_\mu)
 \longrightarrow
 \frac{\{\text{principal two-sided }\mathcal O_\mu\text{-ideals}\}}
      {\mathbb Q^\times}
 \longrightarrow1,
\]
this gives
\[
 |\Aut(\mathcal O_\mu)|
   =\frac{2^{e-1}|\mathcal O_\mu^\times|}{m_\mu}.
\]
These are the orbit and stabilizer calculations of
\cite[Theorem~2.7]{LZ26}; the local quotient still has order two
at every prime dividing $N$. Summing over the fibers yields
\[
 \sum_{i=1}^{h}\frac1{|\mathcal O_i^\times|}
 =\sum_{\mu=1}^{T_{N_1,\,N_2}}\frac{m_\mu}{|\mathcal O_\mu^\times|}
 =2^{e-1}\sum_{\mu=1}^{T_{N_1,\,N_2}}\frac1{|\Aut(\mathcal O_\mu)|},
\]
which proves the second formula.
\end{proof}

\section{Relations of orders and ternary quadratic forms}\label{Sec:Relationships}

We apply the order--form correspondence of
\cite[Theorem~1.1]{LWZ26} to the local orders of
Definition~\ref{def3.6}. The correspondence and the norm--trace
identity are taken from that paper. The local invariants and the
identification of the genera $G_{4N,\,16N^2/\delta,\,N_1}$ are established
here for the present family.

\subsection{Orders of level $(N_1,N_2)$ and ternary quadratic forms}
Write $C_0$ for the even Clifford construction from ternary forms
to quaternion orders, and $C_0^{-1}(\mathcal O)=f_{\mathcal O}$
for the associated integral ternary form; see
\cite[Chapter~22]{Voi21}. In particular,
$d_{f_{\mathcal O}}=\discrd(\mathcal O)$.
Write
\[
 Q^0\!:\,=\{x\in Q:\tr(x)=0\},\quad
 S\!:\,=\mathbb Z+2\mathcal O,\quad S^0\!:\,=S\cap Q^0,
\]
and let $f_{S^0}$ be the reduced norm on $S^0$.
We use $\widetilde f_{\mathcal O}$ with the normalization of
\cite[Definition~2]{LWZ26}; by its Theorem~1.1,
$\widetilde f_{\mathcal O}=\psi_{f_{\mathcal O}}$.
For an odd Gorenstein order this is the primitive part of $f_{S^0}$.
An order is odd if its trace ideal is $\mathbb Z$.
Each order of level $(N_1,N_2)$ is odd: the local order at $2$
contain an element of trace $1$ (in the division case use
$(1+i)/2$ with $i^2=-3$), and $\tr(1)=2$ is a unit at every
odd prime. Hence the global trace ideal is $\mathbb Z$.

The integral map
\[
 \mathcal O/\mathbb Z\longrightarrow S^0,
 \quad [\alpha]\longmapsto2\alpha-\tr(\alpha)
\]
is an isomorphism and identifies the norm on $S^0$ with
$[\alpha]\mapsto4\n(\alpha)-\tr(\alpha)^2$.
Write $\rho_{\mathcal O}(n,r)$ for the number of elements of
$\mathcal O$ with reduced norm $n$ and reduced trace $r$.
Thus \cite[Lemma~4.1]{LWZ26} gives
\begin{equation}\label{eq:integral-norm-trace}
 \rho_{\mathcal O}(n,r)=R_{f_{S^0}}(4n-r^2)
 \quad(n\geq0,\ r\in\mathbb Z).
\end{equation}
We next identify this form and its reciprocal in the present setting.
\begin{prop}\label{prop3.9}
Let $\mathcal O$ be an order of level $(N_1,N_2)$ and put $N=N_1N_2$.
Then $\mathcal O$ is Gorenstein, $m_{f_{\mathcal O}}=1$, and
\[
 f_{S^0}\in G_{4N,\,16N^2,\,N_1},\quad
 f_{\mathcal O}\in G_{4N,\,N,\,N_1},\quad
 \widetilde f_{\mathcal O}=f_{S^0}=\psi_{f_{\mathcal O}}.
\]
The equalities of forms here and below mean integral equivalence.
\end{prop}

\begin{proof}
The local models give $\discrd(\mathcal O)=N$. At odd primes,
direct computation of the norm on $S^0_p=\mathcal O_p^0$ gives
\[
 f_{S^0}\underset p\sim
 \begin{cases}
  -\epsilon_p x^2-p^{v_p(N_1)}y^2
                  +\epsilon_p p^{v_p(N_1)}z^2,&\text{ if }\,\, p\mid N_1,\\
  -x^2-p^{v_p(N_2)}yz,&\text{ if }\,\, p\nmid N_1.
 \end{cases}
\]
At $2$ the corresponding forms are
\[
 f_{S^0}\underset2\sim
 \begin{cases}
  3x^2-2^{v_2(N_1)+2}(y^2+yz+z^2),&\text{ if }\,\, 2\mid N_1,\\
  -x^2-2^{v_2(N_2)+2}yz,&\text{ if }\,\, 2\nmid N_1.
 \end{cases}
\]
For example, in the residually inert case write
$\mathcal O_2=\mathbb Z_2[T]+2^uM_2(\mathbb Z_2)$ with $T=\left(\begin{smallmatrix}1&1\\1&0\end{smallmatrix}\right)$.
The generators $2T-I$, $2^{u+1}E_{12}$, $2^{u+1}E_{21}$ of $S^0_2$
give
\[
 -5x^2-2^{u+2}x(y+z)-2^{2u+2}yz.
\]
Completing the square in $x$ leaves an anisotropic unimodular binary
norm form multiplied by $2^{2u+2}$. The surjectivity of the norm
$(\mathbb Z_2[T])^\times\to\mathbb Z_2^\times$, and the fact that $-5/3$ is a
$2$-adic square, give the displayed inert normal form.

In particular, $f_{S^0}$ is primitive.
If $m$ is the Gorenstein conductor of $\mathcal O$, then
\cite[Theorem~1.1 and Remark~1]{LWZ26} give
$\widetilde m_{f_{S^0}}=m_{f_{\mathcal O}}$
and $m_{f_{\mathcal O}}=m^2m_{f_{\operatorname{Gor}(\mathcal O)}}$.
Thus $m_{f_{\mathcal O}}=1$ and $m=1$.
The same theorem now gives
$f_{S^0}=\widetilde f_{\mathcal O}=\psi_{f_{\mathcal O}}$.
Since $d_{f_{\mathcal O}}=N$ and $m_{f_{\mathcal O}}=1$, reciprocity gives
\[
 N_{f_{\mathcal O}}=N_{f_{S^0}}=4N,
 \quad d_{f_{S^0}}=16N^2.
\]
The local forms above determine the genus conductor as $N_1$.
This proves both genus assertions.
\end{proof}

\begin{cor}\label{cor:tqf}
Let $\mathcal O_1,\ldots,\mathcal O_{T_{N_1,\,N_2}}$ represent the types of level
$(N_1,N_2)$, where $N=N_1N_2$. Then
\[
 T_{N_1,\,N_2}=|G_{4N,\,16N^2,\,N_1}|=|G_{4N,\,N,\,N_1}|,
\]
and
\[
 \sum_{\mu=1}^{T_{N_1,\,N_2}}\frac{\rho_{\mathcal O_\mu}(n,r)}
                            {|\Aut(\mathcal O_\mu)|}
 =2\sum_{f\in G_{4N,\,16N^2,\,N_1}}
                  \frac{R_f(4n-r^2)}{|\Aut(f)|}.
\]
\end{cor}

\begin{proof}
The genus compatibility and bijectivity of the correspondence in
\cite[Theorem~1.1]{LWZ26}, together with
Proposition~\ref{prop3.9}, give the two bijections on classes.
Moreover,
\[
 |\Aut(f_{S^0})|=|\Aut(f_{\mathcal O})|
                   =2|\Aut(\mathcal O)|.
\]
Combining this identity with \eqref{eq:integral-norm-trace} proves
the weighted formula. The factor $2$ comes from using the full
$\operatorname{GL}_3(\mathbb Z)$-automorphism group of a ternary form.
\end{proof}

\begin{prop}\label{the:odd}
Let $\mathcal O$ be an odd Gorenstein order and $S=\mathbb Z+2\mathcal O$.
The following are equivalent:
\begin{enumerate}[label=\textup{(\arabic*)}]
 \item $\mathcal O$ has level $(N_1,N_2)$ for some admissible coprime
       positive integers $N_1,N_2$;
 \item $m_{f_{\mathcal O}}=1$;
 \item $f_{S^0}$ is primitive, or equivalently $f_{\mathcal O^0}$ is primitive;
 \item $\mathcal O_p$ is not residually ramified at any prime;
 \item $N_{f_{\mathcal O}}=4d_{f_{\mathcal O}}$, or equivalently
       $N_{\widetilde f_{\mathcal O}}^2=d_{\widetilde f_{\mathcal O}}$.
\end{enumerate}
\end{prop}

\begin{proof}
By \cite[Corollary~4.2]{LWZ26},
$\widetilde m_{f_{S^0}}=m_{f_{\mathcal O}}$.
The lattices $S^0$ and $\mathcal O^0$ agree locally at odd primes.
They both represent an odd integer at $2$: if $\tr(\alpha)$ is odd,
then $2\alpha-\tr(\alpha)$ lies in both lattices and has odd norm.
Their norm forms therefore have the same content. This proves
\textup{(2)}$\iff$\textup{(3)}.
The equivalence \textup{(2)}$\iff$\textup{(4)} follows from
\cite[Corollary~1.4]{LWZ26}, including its separate assertion at $2$.

For a primitive $f_{\mathcal O}$, the identity
$N_{f_{\mathcal O}}=4d_{f_{\mathcal O}}/m_{f_{\mathcal O}}$
gives the first equivalence in \textup{(5)}.
Also $\widetilde f_{\mathcal O}=\psi_{f_{\mathcal O}}$ and the
discriminant and level formulas for reciprocity give
\[
 \frac{d_{\widetilde f_{\mathcal O}}}
      {N_{\widetilde f_{\mathcal O}}^2}
 =\frac1{m_{f_{\mathcal O}}}.
\]
Thus \textup{(2)}$\iff$\textup{(5)}.
Proposition~\ref{prop3.9} proves \textup{(1)}$\Rightarrow$\textup{(2)}.
Conversely, if \textup{(2)} holds and $N=d_{f_{\mathcal O}}$, then
$\psi_{f_{\mathcal O}}$ is a primitive form of level $4N$ and
discriminant $16N^2$. Proposition~\ref{def:genus-lemma} places it in
$G_{4N,\,16N^2,\,N_1}$ for an admissible exact divisor $N_1$ of $N$.
The bijection in Corollary~\ref{cor:tqf} and its injectivity identify
$\mathcal O$ with an order of level $(N_1,N/N_1)$.
\end{proof}

\begin{cor}
Let $\mathcal O$ be a Gorenstein order. At a division prime $p$,
the residually inert local orders are precisely
$\mathcal O_p^{(2u+1)}$, $u\geq0$; all other local orders are
residually ramified. At a split prime, the nonmaximal residually
inert local orders are $\mathcal O_p^{(2u)}$, $u\geq1$,
and the residually split local orders are the Eichler orders of
positive level. The symbol $*$ occurs precisely for the split
maximal order. The remaining local orders are residually ramified.
\end{cor}
\begin{proof}
This is the local classification of non-residually-ramified
Gorenstein orders; see \cite[Sections~24.3--24.4]{Voi21}.
The explicit representatives and their symbols are given in
Definition~\ref{def3.6} and Proposition~\ref{prop:Eichlersym}.
\end{proof}

The radical idealizer is used locally. For a $\mathbb Z_p$-order put
\[
 \mathcal O_p^\natural
 =\mathcal O_{\operatorname L}(\operatorname{rad}(\mathcal O_p))
   \cap\mathcal O_{\operatorname R}(\operatorname{rad}(\mathcal O_p)).
\]
\begin{lemma}[\cite{Voi21}, 24.4.5 and Proposition~24.4.7]
\label{the:rad}
If $\mathcal O_p$ is Eichler of level $p^{n}$ with $n\geq2$,
then $\mathcal O_p^\natural$ is Eichler of level $p^{n-2}$.
If $\mathcal O_p$ is residually inert and nonmaximal, its radical
idealizer is its unique minimal superorder, which is residually
inert or maximal. In both cases
$[\mathcal O_p^\natural:\mathcal O_p]=p^2$.
\end{lemma}

Fix such a prime $p$. In the following proposition,
$\mathcal O^\natural$ denotes the global overorder characterized by
\[
 (\mathcal O^\natural)_p=\mathcal O_p^\natural,
 \quad (\mathcal O^\natural)_q=\mathcal O_q\quad(q\ne p).
\]
It exists uniquely by the local--global correspondence for lattices.
We use the normalization of Watson's transformation in
\cite[Section~2.5]{LWZ26}. For an integral quadratic lattice
$L$ with norm form $f$, put
\[
 \Lambda_m(L)\!:\,=
 \{v\in L:f(v+w)-f(w)\in m\mathbb Z
                    \text{ for every }w\in L\},
 \quad
 \lambda_m(f)\!:\,=\frac1m f\big|_{\Lambda_m(L)}.
\]
Here the restriction is expressed in an integral basis of
$\Lambda_m(L)$. 
\begin{prop}\label{Prop4.8}
Let $\mathcal O$ be an order of level $(N_1,N_2)$, and suppose that
$p^2\mid N_1$ or $p^2\mid N_2$. Put
$S^\natural=\mathbb Z+2\mathcal O^\natural$.
For odd $p$, we have
\[
 \frac1p\lambda_p(f_{S^0})
       \sim f_{(S^\natural)^0},
 \quad
 \frac1p\lambda_p(f_{\mathcal O^0})
       \sim f_{(\mathcal O^\natural)^0}.
\]
For $p=2$, we have
\[
 \lambda_4(f_{S^0})
       \sim f_{(S^\natural)^0},
 \quad
 \lambda_4(f_{\mathcal O^0})
       \sim f_{(\mathcal O^\natural)^0}.
\]
Moreover,
\[
 T_{N_1,N_2}\geq
 \begin{cases}
 T_{N_1/p^2,\,N_2},&\text{ if }\,\,p^2\mid N_1,\\
 T_{N_1,\,N_2/p^2},&\text{ if }\,\,p^2\mid N_2.
 \end{cases}
\]
The corresponding inequalities for the genera
$G_{4N,\,16N^2,\,N_1}$ follow from Corollary ~{\rm\ref{cor:tqf}}.
\end{prop}

\begin{proof}
The local bases in Definition~\ref{def3.6} and the definition of
$\Lambda_m$ give the following equalities of lattices in $Q^0$:
\[
 \begin{aligned}
 (S^\natural)^0&=p^{-1}\Lambda_p(S^0),&
 (\mathcal O^\natural)^0&=p^{-1}\Lambda_p(\mathcal O^0)
       &&(p\ne2),\\
 (S^\natural)^0&=\tfrac12\Lambda_4(S^0),&
 (\mathcal O^\natural)^0&=\tfrac12\Lambda_4(\mathcal O^0)
       &&(p=2).
 \end{aligned}
\]
At $p$, these follow by imposing the defining congruences on the
local bases; at every other prime the corresponding lattices
are unchanged. The local--global correspondence for lattices
therefore gives the displayed global equalities.
Taking reduced norms proves the asserted form identities.

Watson's transformation is surjective
on the corresponding genus classes.
For odd $p$, the additional factor $p^{-1}$ does not affect
surjectivity on equivalence classes.
The target genera correspond to level $(N_1/p^2,N_2)$ or
$(N_1,N_2/p^2)$, respectively.
Corollary~\ref{cor:tqf} now gives the type number inequalities.
\end{proof}

\subsection{Orders and ternary quadratic forms for extremal cases}

For a Gorenstein order define its reciprocal, up to isomorphism, by
\[
 \psi(\mathcal O)=C_0(\psi_{f_{\mathcal O}}).
\]
By \cite[Corollary~1.5]{LWZ26}, this is an involution on isomorphism
classes and induces a bijection between the corresponding type sets.
It satisfies
\[
 f_{\psi(\mathcal O)}=\widetilde f_{\mathcal O},\quad
 \widetilde f_{\psi(\mathcal O)}=f_{\mathcal O}.
\]
We write $\psi(\mathcal O_p)$ for its local isomorphism class.
The following proposition records the explicit local models needed
for the extremal genera.

\begin{prop}\label{prop:phiO}
Let $p$ be an odd prime and let $\epsilon$ be a quadratic
nonresidue modulo $p$. For $v\geq0$, we have
\begin{equation*}
\psi\left(\begin{pmatrix}
\mathbb{Z}_p & \mathbb{Z}_p\\
p^v\mathbb{Z}_p & \mathbb{Z}_p
\end{pmatrix}\right)\simeq\begin{pmatrix}
a+p^v\mathbb{Z}_p & \mathbb{Z}_p\\
p^v\mathbb{Z}_p & a+p^v\mathbb{Z}_p
\end{pmatrix},~
\psi\left(\begin{pmatrix}
\mathbb{Z}_2 & \mathbb{Z}_2\\
2^v\mathbb{Z}_2 & \mathbb{Z}_2
\end{pmatrix}\right)\simeq\begin{pmatrix}
b+2^{v+2}\mathbb{Z}_2 & \mathbb{Z}_2\\
2^{v+2}\mathbb{Z}_2 & b+2^{v+2}\mathbb{Z}_2
\end{pmatrix},
\end{equation*}  
where $a\in\mathbb{Z}_p$, $b\in\mathbb{Z}_2$. 
Moreover, we have
\begin{equation*}
\psi(\mathcal{O}^{(2u+1)}_p)\simeq\left\{\begin{pmatrix}
\alpha & p^u\beta\\
p^{u+1}\bar{\beta} & \bar{\alpha}
\end{pmatrix}
: \alpha\in \mathbb{Z}_p+p^{2u+1}\sqrt{\epsilon}\mathbb{Z}_p,\,\beta\in \mathbb{Z}_p+\sqrt{\epsilon}\mathbb{Z}_p\right\},
\end{equation*} 
and
\begin{equation*}
\psi(\mathcal{O}^{(2u)}_p)\simeq\left\{ 
			\begin{pmatrix}
				a+p^uc & \epsilon (p^{2u}b-p^ud)\\
				p^{2u}b+p^ud & a-p^uc 
			\end{pmatrix}
			:a,b,c,d\in\mathbb{Z}_p\right\}.
\end{equation*} 
Further, for $p=2$, we have
\begin{equation*}
\psi(\mathcal{O}^{(2u+1)}_2)\simeq\left\{\begin{pmatrix}
\alpha & 2^{u+1}\beta\\
2^{u+2}\bar{\beta} & \bar{\alpha}
\end{pmatrix}
: \alpha\in\mathbb{Z}_2+2^{2u+2}\sqrt{-3}\mathbb{Z}_2,\,\beta\in \mathbb{Z}_2+\frac{1+\sqrt{-3}}{2}\mathbb{Z}_2\right\},  
\end{equation*} 
and
\begin{equation*}
\psi(\mathcal{O}^{(2u)}_2)\simeq\left\{ 
			\begin{pmatrix}
				a+2^{u+1}b & 2^{u+1}c\\
				2^{u+1}d & a-2^{u+1}b 
			\end{pmatrix}
			:a,b,c,d\in\mathbb{Z}_2,\,b+c+d\equiv0\!\!\!\pmod{2^u}\right\}.
\end{equation*} 
\end{prop}

\begin{proof}
Apply the even Clifford construction to the local forms for
$f_{S^0}$ in Proposition~\ref{prop3.9}.
For example, at an odd division prime take
$i^2=\epsilon$, $j^2=p$, and $ji=-ij$.
The resulting order has basis
\[
 1,\quad p^{2u+1}i,\quad p^uj,\quad p^uij,
\]
which gives the displayed matrix model. At a dyadic division prime,
take $i^2=-3$, $j^2=2$, and $\omega=(1+i)/2$.
A basis of the reciprocal order is
\[
 1,\quad2^{2u+2}i,\quad2^{u+1}j,\quad2^{u+1}\omega j.
\]
The split formulas follow by the same Clifford multiplication in
$M_2(\mathbb Q_p)$. Their reduced discriminant exponents are
$2v$ and $4u$ at odd primes, and $2v+4$ and $4u+4$ at $2$,
respectively. In the division cases the exponents are $4u+2$
at odd primes and $4u+6$ at $2$, as also follows from the displayed
bases. These agree with the discriminants prescribed by reciprocity.
\end{proof}

\begin{definition}\label{partial}
Let $\mathcal O$ be an order of level $(N_1,N_2)$, put $N=N_1N_2$, and let
$\delta\parallel16N$. An order $\mathcal O'$ has level
$(N_1,N_2;\,\delta)$ if
\[
 \mathcal O'_p\simeq
 \begin{cases}
 \psi(\mathcal O_p),&\text{ if }\,\,p\mid\delta,\\
 \mathcal O_p,&\text{ if }\,\,p\nmid\delta.
 \end{cases}
\]
\end{definition}
These local conditions are realized by a global order: apply
Lehman's correspondence $\phi_\delta$ to $f_{\mathcal O}$ and then
apply $C_0$. The local compatibility established in
\cite[Corollary~1.5]{LWZ26} identifies its local orders with the
ones in the definition. It also gives a bijection of the two type
sets. This assertion concerns local orders; the associated globally
normalized ternary forms are in general only unit-similar away
from the selected primes, as explained in Remark~\ref{rem:phiN}.

Every order $\mathcal O'$ of level $(N_1,N_2;\,\delta)$ is Gorenstein.
Put $S'=\mathbb Z+2\mathcal O'$. The discriminant formulas and the
correspondence give
\[
 d_{f_{\mathcal O'}}=\discrd(\mathcal O')=N\delta,
 \quad N_{f_{\mathcal O'}}=4N,
 \quad m_{f_{\mathcal O'}}=\delta,
\]
and
\[
 \widetilde f_{\mathcal O'}=f_{(S')^0}/\delta,
 \quad N_{\widetilde f_{\mathcal O'}}=4N,
 \quad d_{\widetilde f_{\mathcal O'}}=16N^2/\delta,
 \quad m_{\widetilde f_{\mathcal O'}}=16N/\delta.
\]
The normalization by $\delta$ remains valid when $\mathcal O'$ is
even: in that case $(S')^0=2(\mathcal O')^0$ and
$m_{f_{\mathcal O'}}=4\widetilde m_{f_{(\mathcal O')^0}}$.
By Definition~\ref{def:genus111},
\[
 f_{\mathcal O'}\in G_{4N,\,N\delta,\,N_1},\quad
 \widetilde f_{\mathcal O'}\in G_{4N,\,16N^2/\delta,\,N_1}.
\]

\vskip 0.2cm
An order $\mathcal{O}$ is a Bass order if and only if every order $\mathcal{O}'\supset \mathcal{O}$ is a Gorenstein order.
 Pacetti and Sirolli~\cite{PS14} gave a tree of Bass orders on a totally definite quaternion algebra $Q$ over a totally real number field $K$, 
 and classified all the Bass orders into the so-called classes A1, A2, B, and C.

\vskip 0.2cm
Now we can establish a concise table, which gives an illustration for the relations of the classifications of orders of level $(N_1,N_2;\,\delta )$
 and ternary quadratic forms for extremal cases. The relations with the  classification of Bass orders by Pacetti and Sirolli~\cite{PS14} are also indicated.
\begin{table}[H]
\caption{Orders of level $(N_1,N_2;\,\delta )$ and ternary quadratic forms of the extremal cases}
\centering
\begin{tabular}{c|c|c|c}
\hline
Odd primes & Orders of level $(N_1,N_2;\,\delta )$ & Classes of $\mathcal{O}_p$ & Genus symbols\\
\hline
$p\mid N_1$ and $p\nmid\delta $ & residually inert & A2 &$\left(\frac{\psi_{f_{\mathcal{O}}}}{p}\right)=-\left(\frac{-\delta}{p}\right)$, $\left(\frac{{f_{\mathcal{O}}}}{p}\right)=\times$ \\
$p\mid N_2$ and $p\nmid\delta $ & residually split & A1 &$\left(\frac{\psi_{f_{\mathcal{O}}}}{p}\right)=\left(\frac{-\delta}{p}\right)$, $\left(\frac{{f_{\mathcal{O}}}}{p}\right)=\times$ \\
$p\mid\delta $ & residually ramified & B (if $p\parallel\delta )$  &$\left(\frac{\psi_{f_{\mathcal{O}}}}{p}\right)=\times$, $\left(\frac{{f_{\mathcal{O}}}}{p}\right)$ makes sense\\
$p\nmid N_1N_2$ & $\left(\frac{\mathcal{O}}{p}\right)=*$ & A1 &$\left(\frac{\psi_{f_{\mathcal{O}}}}{p}\right)=\times$, $\left(\frac{{f_{\mathcal{O}}}}{p}\right)=\times$\\
\hline
\end{tabular}
\end{table}
\noindent Here $\times$ is used to indicate that the relevant genus symbols do not make sense.
For an odd prime $p$, if $v_p(\delta)>1$, then $\mathcal O_p$ is not a Bass order.
The table concerns odd primes, and for $p=2$ one can give similar assertions. The corresponding genus symbols for
$f_{S^0}/\delta$ follow from Remark~\ref{rem:genussym}.

\section{Generalization of the modified Hurwitz class numbers}
\label{Sec:4}

We extend the modified Hurwitz class numbers of
\cite[Definition~4.2]{LZ26} by allowing positive even valuations
in $N_1$. We give explicit local factors for $\delta=1$ and
recursive formulas for general $\delta$.

\begin{definition}\label{def:HCN11}
Let $N_1,N_2$ be coprime positive integers. For a
negative discriminant $-D$, let $\mathfrak{f}=\mathfrak{f}_{N_1,N_2}$ be the largest
positive integer, whose prime factors divide $N_1N_2$, 
such that $\mathfrak{f}^2\mid D$ and $-D/\mathfrak{f}^2$ is a discriminant.
Put, for $p\mid N_1N_2$, 
\[
 \mathfrak{f}_p\!:\,=p^{v_p(\mathfrak{f})},\quad
 \alpha_p\!:\,=v_p(N_1N_2),\quad
 \chi_p\!:\,=\left(\frac{-D/\mathfrak{f}^2}{p}\right).
\]
With the usual Hurwitz class number $H(\cdot)$, define
\begin{equation}\label{Hn1n2}
 H^{(N_1,N_2;\,1)}(D)\!:\,=H^{(N_1,N_2)}(D)\!:\,
 =H(D/\mathfrak{f}^2)\prod_{p\mid N_1N_2}A_p(D;N_1,N_2),
\end{equation}
where, with the column conditions $p\mid N_1$ and $p\mid N_2$, as well as the row conditions with respect to the valuations,
the local factors $A_p(D;N_1,N_2)$ are given in Table~\ref{tab:local_factors} below.

\begin{table}[H]
\centering
\caption{Local factors $A_p(D;N_1,N_2)$}
\label{tab:local_factors}
\small
\setlength{\tabcolsep}{4pt}
\renewcommand{\arraystretch}{1.35}
\begin{tabular}{l|cc}
\toprule
 & $p\mid N_1$ & $p\mid N_2$\\
\midrule
$v_p(p\mathfrak{f}^2)<\alpha_p$, $p\mid D/\mathfrak{f}^2$
 & $0$ & $0$\\
$v_p(p\mathfrak{f}^2)<\alpha_p$, $p\nmid D/\mathfrak{f}^2$
 & $\mathfrak{f}_p^2(1-\chi_p)$
 & $\mathfrak{f}_p^2(1+\chi_p)$\\
$v_p(p\mathfrak{f}^2)\geq \alpha_p$, $\alpha_p$ odd
 & $p^{\alpha_p-1}(1-\chi_p)$
 & $\dfrac{2p^{(\alpha_p-1)/2}\mathfrak{f}_p(p-\chi_p)
              -(p+1)p^{\alpha_p-1}(1-\chi_p)}{p-1}$\\
$v_p(p\mathfrak{f}^2)\geq \alpha_p$, $\alpha_p$ even
 & $p^{\alpha_p/2-1}\mathfrak{f}_p(p-\chi_p)-p^{\alpha_p-1}(1-\chi_p)$
 & $\dfrac{p+1}{p-1}
       \bigl(p^{\alpha_p/2-1}\mathfrak{f}_p(p-\chi_p)
                 -p^{\alpha_p-1}(1-\chi_p)\bigr)$\\
\bottomrule
\end{tabular}
\end{table}
For positive $D\equiv1,2\pmod4$, set
$H^{(N_1,N_2)}(D)=0$. At zero, set
\[
 H^{(N_1,\,N_2)}(0)\!:\,
 =(-1)^{e(N_0)+1}\frac{N_1N_2}{12}
   \prod_{p\mid N_1}\left(1-\frac1p\right)
   \prod_{p\mid N_2}\left(1+\frac1p\right),
\]
where $N_0$ is squarefree such that $N_1=s^2N_0$.
\end{definition}

\begin{definition}\label{def:HCNd}
Let $\delta\parallel16N_1N_2$ and $D\geq0$, with the convention
$v_p(0)=+\infty$. For primes $p\nmid\delta$, define recursively
as follows.
For $p^u\parallel N_1$ with $u\geq1$, define, if $v_p(D)<u$ (the low-valuation branch), 
\[
\begin{split}
 H^{(N_1,N_2;16^{\delta_{p,2}}p^u\delta)}(D)\!:\,
 ={}&\frac{1+(-1)^{v_p(D)}}{2p}
 \Bigl[(p+1)H^{(N_1/p^u,\,N_2;\,\delta)}
                         (4^{\delta_{p,2}}p^uD)\\
 &\hspace{30mm}
       -H^{(N_1/p^u,\,pN_2;\,\delta)}
                         (4^{\delta_{p,2}}p^uD)\Bigr],
\end{split}
\]
and, if $v_p(D)\geq u$ (the high-valuation branch), define
\[
 H^{(N_1,N_2;16^{\delta_{p,2}}p^u\delta)}(D)\!:\,
 =
 \begin{cases}
 p^{u-1}(p-1)H^{(N_1/p^u,\,N_2;\,\delta)}
             (4^{\delta_{p,2}}D/p^u),
       &\text{ if }\,\,u\text{ is even},\\[3pt]
 p^{u-1}H^{(N_1/p^{u-1},\,N_2;\,\delta)}
             (4^{\delta_{p,2}}D/p^u),
       &\text{ if }\,\,u\text{ is odd}.
 \end{cases}
\]

For $p^v\parallel N_2$ with $v\geq1$, define, if $v_p(D)<v$ (the low-valuation branch), 
\[
\begin{split}
 H^{(N_1,N_2;\,16^{\delta_{p,2}}p^v\delta)}(D)\!:\,
 ={}&\frac{v_p(D)+1}{p}
 \Bigl[(p+1)H^{(N_1,N_2/p^v;\,\delta)}
                         (4^{\delta_{p,2}}p^vD)\\
 &\hspace{26mm}
       -H^{(N_1,N_2/p^{v-1};\,\delta)}
                         (4^{\delta_{p,2}}p^vD)\Bigr],
\end{split}
\]
and, if $v_p(D)\geq v$ (the high-valuation branch), define
\[
\begin{split}
 H^{(N_1,N_2;\,16^{\delta_{p,2}}p^v\delta)}(D)\!:\,
 =p^{v-1}\Bigl[
 &(v+1)(p+1)H^{(N_1,N_2/p^v;\,\delta)}
                         (4^{\delta_{p,2}}D/p^v)\\
 &-vH^{(N_1,N_2/p^{v-1};\,\delta)}
                         (4^{\delta_{p,2}}D/p^v)
 \Bigr].
\end{split}
\]

If $2\nmid N_1N_2\delta$, set
\[
 H^{(N_1,N_2;\,16\delta)}(D)\!:\,
 =H^{(N_1,N_2;\,\delta)}(4D).
\]
\end{definition}

\begin{remark}\label{remark11}
The recursion in Definition~\ref{def:HCNd} is independent of the order in which the prime power factors of $\delta$ are selected.
 To make this precise, regard the values $H^{(N_1,\,N_2;\,\delta)}(D)$, with the auxiliary parameter pairs allowed in Definition~\ref{def:HCNd}, 
 as entries of a formal family. For a prime $p$ not yet selected, let $T_p$ denote the entire linear substitution prescribed by Definition~\ref{def:HCNd}: 
 it includes the change of $\delta$, the relevant $p$-primary part of $(N_1,N_2)$, the represented argument, 
 and every coefficient in the applicable valuation branch. We claim that
\[
 T_pT_q=T_qT_p\quad(p\ne q)
\]
whenever both sides are defined.

To see this, we fix a branch at each of $p$ and $q$. Each branch is a finite linear combination of terms obtained by {\rm (i)} 
modifying only the corresponding prime-primary part of $(N_1,N_2)$, and {\rm (ii)} replacing $D$ by $c_rD$, 
where $c_r$ is a rational power of that prime, multiplied by $4$ if $r=2$ as specified in Definition~\ref{def:HCNd}. 
These operations for distinct primes commute. Their scalar coefficients depend only on the unchanged local valuation and parameter data. 
On expanding either composite substitution, corresponding terms therefore have the same coefficient, the same final parameter pair, 
the same selected divisor, and the same final represented argument. This proves the claimed equality, including the branches containing two summands.

At $D=0$ we use $v_p(0)=+\infty$, so the high-valuation branches apply at every prime and the same argument remains valid.
 Thus every ordering of the distinct prime power selections yields the same value of $H^{(N_1,N_2;\delta)}(D)$.
\end{remark}

\begin{remark}\label{remark4.4}
For $D>0$ and $v_p(N_1N_2)=1$, the definitions give
\[
 H^{(N_1,N_2;\,16^{\delta_{p,2}}p)}(D)
 =H^{(N_1,N_2)}(4^{\delta_{p,2}}pD).
\]
Consequently, if $\delta\parallel N_1N_2$ is odd and squarefree,
then
\[
 H^{(N_1,N_2;\,\delta)}(D)=H^{(N_1,N_2)}(\delta D).
\]
This recovers the squarefree definition in \cite{LZ24}.
\end{remark}	
	\begin{remark}
Definition~\ref{def:HCNd} avoids a longer closed formula with multiple
valuation cases and suffices for exact computation.
For example, with $N_1=180=2^2\cdot3^2\cdot5$, $N_2=1$, and
$\delta=576=16\cdot2^2\cdot3^2$, it gives
\begin{equation*}
			H^{(180,\,1;\,576)}(D)=
			\begin{cases}
				\dfrac{(1+(-1)^{v_2(D)})(1+(-1)^{v_3(D)})}{24} & \\
				\quad \times \bigl(12H^{(5,1)}(144D)-3H^{(5,3)}(144D) & \\
				\quad -4H^{(5,2)}(144D)+H^{(5,6)}(144D)\bigr), & \text{if $v_2(D),v_3(D)<2,$} \\[15pt]
				\dfrac{1+(-1)^{v_2(D)}}{2}\bigl(9H^{(5,1)}(16D/9)-3H^{(5,2)}(16D/9)\bigr), & \text{if $v_2(D)<2,v_3(D)\geq2,$} \\[10pt]
				\dfrac{1+(-1)^{v_3(D)}}{3}\bigl(4H^{(5,1)}(9D)-H^{(5,3)}(9D)\bigr), & \text{if $v_2(D)\geq2,v_3(D)<2,$} \\[10pt]
				12H^{(5,1)}(D/9), & \text{if $v_2(D),v_3(D)\geq2.$}
			\end{cases}
		\end{equation*}
	\end{remark}
	\section{Proof of Theorems}\label{Sec:proofs}
	
	In this section we complete the proofs of Theorems~\ref{the:genera},~\ref{the:repnum} and~\ref{the:classnum}.
	We use the notations $N$, $N_1$ of Theorem~\ref{the:genera}, and we set $N_2=N/N_1$.
	\subsection{Proof of Theorem~\ref{the:genera}}
Proposition~\ref{def:genus-lemma} proves the assertion for $\delta=1$.
Lehman's correspondences are bijections on classes and on genera by
Lemma~\ref{prop:bijection}; applying $\phi_\delta$ gives precisely the
genera in Definition~\ref{def:genus111}.
Proposition~\ref{jiajia} recovers their conductors intrinsically.
When $2\parallel N$, the same argument applies to
$\phi_{\delta'}\circ\lambda_4$, using the bijection in
Lemma~\ref{prop:bijection} and Remark~\ref{remark2.13}.
This proves both assertions.

\subsection{Proof of Theorem~\ref{the:repnum}}

We follow the proof of \cite[Theorem~4.3]{LZ26}, given in its
Section~4.2, and retain only the additional local calculations.
For a positive definite ternary genus $G$ and $n>0$, with discriminant $d_G$, the
Siegel--Weil formula reads
\begin{equation}\label{gongshi5.1}
 \sum_{f\in G}\frac{R_f(n)}{|\Aut(f)|}
 =4\pi M(G)\sqrt{\frac{n}{d_G}}\prod_p d_{G,p}(n),
 \quad
 M(G)=\sum_{f\in G}\frac1{|\Aut(f)|}.
\end{equation}
Here, the sums are over class representatives, and the density
normalization is, for a representative $f\in G$,
\[
 d_{G,\,p}(n)=d_{f,\,p}(n)
 =\lim_{t\to\infty}p^{-2t}
 \#\{(x,y,z)\in(\mathbb{Z}/p^t\mathbb{Z})^3:
          f(x,y,z)\equiv n\!\!\!\pmod{p^t}\}.
\]
For the local forms of Proposition~\ref{prop3.9}, the odd-prime
formulas with $v_p(N_1)$ odd and the formulas at $p\mid N_2$
are \cite[Propositions~4.4 and~4.5]{LZ26}; their dyadic
counterparts are \cite[Propositions~4.6 and~4.7]{LZ26}.
We record the additional even-exponent and reciprocal
formulas needed here. Throughout the next three lemmas,
$k\geq0$, $p\nmid m$, and $\epsilon_p$ is a nonsquare unit
when $p$ is odd.

\begin{lemma}\label{lem:density-even-inert}
Let $u_0>0$ be even. For odd prime $p$,
\[
 d_{-\epsilon_p x^2-p^{u_0}y^2+\epsilon_p p^{u_0}z^2,\,p}(p^km)
 =
 \begin{cases}
  p^{k/2}\left(1-\left(\dfrac{-m}{p}\right)\right),
                         &\text{ if }\,\,k<u_0,\ k\text{ even},\\
 0,                     &\text{ if }\,\,k<u_0,\ k\text{ odd},\\
 p^{u_0/2}d_{-x^2-yz,\,p}(p^{k-u_0}m),
                         &\text{ if }\,\,k\geq u_0.
 \end{cases}
\]
At $p=2$,
\[
 d_{3x^2-2^{u_0+2}(y^2+yz+z^2),\,2}(2^km)
 =
 \begin{cases}
 2^{k/2+2},
       &\text{ if }\,\,k<u_0,\ k\text{ even},\ m\equiv3\!\!\!\pmod8,\\
        2^{u_0/2}d_{-x^2-4yz,\,2}(2^{k-u_0}m),
       &\text{ if }\,\,k\geq u_0,\\
       0,    &\text{ otherwise}.
          \end{cases}
\]
\end{lemma}

\begin{proof}
When $k\geq2$, the congruence forces $p\mid x$.
Substituting $x=px'$ and dividing by $p^2$ lowers
$k$ and $u_0$ by two and multiplies the density by $p$.
For $k<u_0$, iteration leaves a unary square congruence:
an odd valuation $k$ is impossible, whereas an even valuation
gives the first branch.
At $2$, an odd square is $1$ modulo $8$, and each soluble
unit-square congruence modulo $2^t$, $t\geq3$, has four roots.
This gives the stated factor $2^{k/2+2}$.

For $k\geq u_0$, iterate $u_0/2$ times.
At odd $p$, the remaining unimodular ternary form is
equivalent to $-x^2-yz$. At $2$, the identity
\[
 3x^2-4(y^2+yz+z^2)
 =-(x+2y+2z)^2+4(x+y)(x+z)
\]
gives an integral change of variables of determinant $3$,
hence equivalence to $-x^2-4yz$ over $\mathbb Z_2$.
\end{proof}

\begin{lemma}\label{lem:density-reciprocal-inert}
Let $u\geq1$. For odd $p$,
\[
 d_{-\epsilon_p p^ux^2-y^2+\epsilon_p z^2,\,p}(p^km)
 =
 \begin{cases}
 \dfrac{1+(-1)^k}{2}\left(1+\dfrac1p\right),
       &\text{ if }\,\,k<u,\\[3pt]
 d_{-x^2-yz,\,p}(p^{k-u}m),
       &\text{ if }\,\,k\geq u,\ u\text{ even},\\[3pt]
 \dfrac1p
 d_{-\epsilon_p x^2-py^2+\epsilon_p pz^2,\,p}(p^{k-u}m),
       &\text{ if }\,\,k\geq u,\ u\text{ odd}.
 \end{cases}
\]
At $p=2$,
\[
 d_{3\cdot2^ux^2-(y^2+yz+z^2),\,2}(2^km)
 =
 \begin{cases}
 \dfrac{3(1+(-1)^k)}4,
       &\text{ if }\,\,k<u,\\[3pt]
 d_{-x^2-yz,\,2}(2^{k-u}m),
       &\text{ if }\,\,k\geq u,\ u\text{ even},\\[3pt]
 \dfrac14
 d_{3x^2-8(y^2+yz+z^2),\,2}(2^{k-u+2}m),
       &\text{ if }\,\,k\geq u,\ u\text{ odd}.
 \end{cases}
\]
\end{lemma}

\begin{proof}
The binary forms $-y^2+\epsilon_p z^2$ and $y^2+yz+z^2$
are unramified norm forms. A unit has binary density
$1+1/p$, and an odd valuation is impossible.
If $u,k\geq2$, the congruence forces $p\mid y,z$;
substituting $y=py'$, $z=pz'$ lowers $u,k$ by two
without changing the ternary density.
This proves the branches $k<u$ and reduces $k\geq u$
to $u=0$ or $u=1$.

At $u=0$, the remaining form is equivalent to $-x^2-yz$.
For $p=2$ this follows from
\[
 3x^2-(y^2+yz+z^2)
 =-(x+y)^2+(2x-z)(2x+y+z),
\]
whose change of variables has determinant $3$.
At odd $p$ and $u=1$, substituting $y=py'$, $z=pz'$
once gives the factor $1/p$ in the last branch.
At $2$ the corresponding reduction is
\[
\begin{split}
 d_{6x^2-(y^2+yz+z^2),\,2}(2n)
 &=\frac12d_{3x^2-2(y^2+yz+z^2),\,2}(n)\\
 &=\frac14d_{3x^2-8(y^2+yz+z^2),\,2}(4n).
\end{split}
\]
The last equality follows by substituting $x=2x'$
in the density on the right.
These normalized congruence counts also explain
the dyadic factor $1/4$.
\end{proof}

\begin{lemma}\label{lem:density-reciprocal-split}
For every prime $p$ including $2$, and positive integer $v\geq1$,
\[
 d_{-p^vx^2-yz,\,p}(p^km)
 =
 \begin{cases}
 (k+1)(1-p^{-1}),
       &\text{ if }\,\,k<v,\\
 d_{-x^2-yz,\,p}(p^{k-v}m)+v(1-p^{-1}),
       &\text{ if }\,\,k\geq v.
 \end{cases}
\]
Moreover, for every $n>0$,
\begin{equation}\label{eq:density-base-identity}
 p\,d_{-x^2-yz,\,p}(n)-d_{-x^2-pyz,\,p}(n)=p-1.
\end{equation}
\end{lemma}

\begin{proof}
If $k<v$, each possible valuation $0,\ldots,k$ of $y$
contributes $1-p^{-1}$ to the density.
If $k\geq v$, the valuations $0,\ldots,v-1$ contribute
$v(1-p^{-1})$; the remaining solutions have $p^v\mid y$,
and substitution and division by $p^v$ give
$d_{-x^2-yz,\,p}(p^{k-v}m)$.
Finally, in $-x^2-yz$, the solutions with $p\nmid y$
contribute $1-p^{-1}$, while those with $p\mid y$
contribute $p^{-1}d_{-x^2-pyz,\,p}(n)$.
This proves \eqref{eq:density-base-identity} for every $p$.
\end{proof}

\begin{proof}[Proof of Theorem~{\upshape\ref{the:repnum}}]
By Proposition~\ref{prop:mass3}, Corollary~\ref{cor:tqf}, and Lemma~\ref{prop:bijection}, the mass of every genus under consideration is
\begin{equation}\label{eq:mass}
 M\bigl(G_{4N,\,16N^2/\delta,\,N_1}\bigr)
 =2^{-e(N)-1}\frac{N}{12}
 \prod_{p\mid N_1}(1-p^{-1})\prod_{p\mid N_2}(1+p^{-1}).
\end{equation}
Since every positive definite form represents zero only by the zero vector, its representation number at zero is one. 
Thus \eqref{eq:mass}, the parity assumption on $e(N_0)$, and Remark~\ref{remark11} prove the first assertion for $n=0$. We assume $n>0$ from now on.

\medskip
We first consider the initial case of $\delta=1$.
For $n\equiv1,2\pmod4$, the two sides vanish by the dyadic condition and Definition~\ref{def:HCN11}. 
Suppose therefore that $n\equiv0,3\pmod4$. We use the classical Hurwitz class number identity, in the density normalization of ~\eqref{gongshi5.1},
\begin{equation}\label{eq:hurwitz}
 H(D)=\frac{\pi\sqrt D}{12}\,d_{-x^2-4yz,\,2}(D)
 \prod_{p\ne2}d_{-x^2-yz,\,p}(D),
 \quad D>0,\quad D\equiv0,3\!\!\pmod4.
\end{equation}
For fundamental discriminants this is the Dirichlet class number formula expressed using the base densities;
 the conductor formula for $H$ and the corresponding density recurrences extend it to every negative discriminant.

Let $\mathfrak{f},\mathfrak{f}_p,\chi_p$ be as in Definition~\ref{def:HCN11}, with $D=n$. 
The local density formulas cited before Lemma~\ref{lem:density-even-inert}, together with that lemma, give, for $p\mid N$,
\begin{equation}\label{eq:localbase}
 d_{G_{4N,\,16N^2,\,N_1},\,p}(n)=
 \begin{cases}
 \displaystyle\frac{2^{\delta_{p,2}}(p+1)}{{\mathfrak{f}}_p(p-\chi_p)}
 A_p(n;N_1,N_2),&p\mid N_1,\\[6pt]
 \displaystyle\frac{2^{\delta_{p,2}}(p-1)}{{\mathfrak{f}}_p(p-\chi_p)}
 A_p(n;N_1,N_2),&p\mid N_2.
 \end{cases}
\end{equation}
The maximality of $\mathfrak{f}$ gives the base density
\[
 \begin{cases}
 d_{-x^2-4yz,\,2}(n/{\mathfrak{f}}^2),&p=2,\\
 d_{-x^2-yz,\,p}(n/{\mathfrak{f}}^2),&p\ne2
 \end{cases}
 =\frac{2^{\delta_{p,2}}(p^2-1)}{p(p-\chi_p)}
 \quad(p\mid N).
\]
Consequently, multiplying \eqref{eq:localbase} by $1-p^{-1}$ when $p\mid N_1$, and by $1+p^{-1}$ when $p\mid N_2$,
 gives in either case the base density at $n/\mathfrak{f}^2$, multiplied by ${\mathfrak{f}}_p^{-1}A_p(n;N_1,N_2)$. 
 At every $p\nmid N$, the local form is the base form and ${\mathfrak{f}}$ is a $p$-adic unit; 
 hence the local density is unchanged when the represented argument is divided by ${\mathfrak{f}}^2$.

Applying ~\eqref{gongshi5.1} and \eqref{eq:mass}, multiplying these local identities over all primes,
 and using ${\mathfrak{f}}=\prod_{p\mid N}{\mathfrak{f}}_p$, we obtain
\begin{align*}
 \sum_{f\in G_{4N,\,16N^2,\,N_1}}\frac{R_f(n)}{|\operatorname{Aut}(f)|}
 &=2^{-e(N)-1}\frac{\pi\sqrt{n/{\mathfrak{f}}^2}}{12}
 d_{-x^2-4yz,\,2}(n/{\mathfrak{f}}^2)
 \prod_{p\ne2}d_{-x^2-yz,\,p}(n/{\mathfrak{f}}^2)
 \prod_{p\mid N}A_p(n;N_1,N_2)\\
 &=2^{-e(N)-1}H(n/{\mathfrak{f}}^2)\prod_{p\mid N}A_p(n;N_1,N_2)\\
 &=2^{-e(N)-1}H^{(N_1,\,N_2)}(n),
\end{align*}
by \eqref{eq:hurwitz} and Definition~\ref{def:HCN11}. This proves the first assertion when $\delta=1$.

\medskip
Turn to the general case. We consider the expression, for arbitrary coprime positive integers $N_1,N_2$,
\begin{equation}\label{eq:formalproduct}
 \frac{\pi\sqrt{\delta n}}{12}
 \prod_{p\mid N_1}(1-p^{-1})
 \prod_{p\mid N_2}(1+p^{-1})
 \prod_q d_{\text{corresponding local form at }q,\,\,q}(n).
\end{equation}
Here ``corresponding local form'' means the form specified as follows:
first use the local forms of Proposition~\ref{prop3.9} for the base form, and at each selected prime $p\mid \delta$ use 
the reciprocal form of Lemma~\ref{lem:density-reciprocal-inert} or Lemma~\ref{lem:density-reciprocal-split};
 if $2\nmid N_1N_2$, the reciprocal dyadic form is taken as $-x^2-yz$. Then, at an unselected prime $p$,
  multiply the base form by $\delta$; at a selected prime $p$, multiply its reciprocal form by $\delta/p^{v_p(\delta)}$.

The preceding calculation, being local, identifies \eqref{eq:formalproduct} for $\delta=1$ with $H^{(N_1,\,N_2)}(n)$ 
whenever the base local forms and the factors of Definition~\ref{def:HCN11} are applicable.  
This identification is used only at the parameter pairs occurring in Definition~\ref{def:HCNd}.

For the admissible parameters in Theorem~\ref{the:repnum}, Proposition~\ref{prop3.9} and Remark~\ref{rem:phiN} 
identify the relevant forms with the localizations of $G_{4N,\,16N^2/\delta,\,N_1}$, 
up to the unit square changes in represented arguments implicit in the density comparison. 
Thus ~\eqref{gongshi5.1} and \eqref{eq:mass} identify \eqref{eq:formalproduct} with
\[
 2^{e(N)+1}\sum_{f\in G_{4N,\,16N^2/\delta,\,N_1}}
 \frac{R_f(n)}{|\operatorname{Aut}(f)|}.
\]
Therefore, it remains to check that \eqref{eq:formalproduct} satisfies all the recursions in Definition~\ref{def:HCNd}.

\medskip
Suppose $p\nmid\delta$ and $p^u\parallel N_1$ or $p^v\parallel N_2$. Selecting $p$ replaces
 $\delta$ by $16^{\delta_{p,2}}p^u\delta$ or $16^{\delta_{p,2}}p^v\delta$, respectively. 
 The arguments on the right hand side of Definition~\ref{def:HCNd} are $4^{\delta_{p,2}}p^un$
  and $4^{\delta_{p,2}}n/p^u$, or their $v$-analogues; 
  division by $p^u$ or $p^v$ is used only in the corresponding high-valuation branch.

We make the comparison at an arbitrary prime $q\ne p$ explicit.  By Remark~\ref{rem:phiN},
applying the $\phi_p$ correspondence multiplies the relevant local form at $q$
by $p^u$ (respectively $p^v$), up to a change of variables in
$\operatorname{GL}_3(\mathbb Z_q)$.  The local forms used in
\eqref{eq:formalproduct} have already absorbed these similarity factors
into their scalar coefficients.  More precisely, at an unselected $q$
the coefficient is the third parameter, while at a previously selected
$q$ it is that parameter divided by its $q$-primary part.  Since
$p\ne q$, the latter divisor is unchanged by the present step.
Thus the ratio of the two scalar coefficients at $q$ is
$16^{\delta_{p,2}}p^u$ (respectively $16^{\delta_{p,2}}p^v$).

For a quadratic form $Q$ over $\mathbb Z_q$, a unit $c$, and an
integer $m$, the definition of the local density gives
\[
 d_{cQ,\,q}(m)=d_{Q,\,q}(m/c),\quad
 d_{Q,\,q}(a^2m)=d_{Q,\,q}(m)
 \quad(a\in\mathbb Z_q^\times).
\]
The first identity follows by dividing the congruence by $c$;
the second follows from the bijection $X\mapsto aX$ modulo every
power of $q$.  Consequently, after transferring the scalar
coefficients to the represented arguments, the low-valuation
comparison at $q$ is between
\[
 \frac{4^{\delta_{p,2}}p^un}{16^{\delta_{p,2}}p^u\delta}
 \quad\text{and}\quad \frac n\delta,
 \quad
 \frac{4^{\delta_{p,2}}p^un/(16^{\delta_{p,2}}p^u\delta)}{n/\delta}
 =4^{-\delta_{p,2}}
 =\bigl(2^{-\delta_{p,2}}\bigr)^2.
\]
In the high-valuation comparison, the analogous ratio is
\[
 \frac{4^{\delta_{p,2}}n/(p^u16^{\delta_{p,2}}p^u\delta)}{n/\delta}
 =4^{-\delta_{p,2}}p^{-2u}
 =\bigl(2^{-\delta_{p,2}}p^{-u}\bigr)^2.
\]
Here $2^{-\delta_{p,2}}$ is a $q$-adic unit: if
$\delta_{p,2}=1$, then $p=2$ and $q\ne2$.  Also $p^{-u}$ is a
$q$-adic unit because $q\ne p$.  The same calculation with $v$ in
place of $u$ applies when $p^v\parallel N_2$.  Thus the local
density factors at every $q\ne p$ agree in each recursion branch.

The square-root prefactors contribute, respectively,
\begin{equation}\label{eq:roots}
 \frac{\sqrt{16^{\delta_{p,2}}p^u\delta n}}
 {\sqrt{\delta\,4^{\delta_{p,2}}p^un}}
 =2^{\delta_{p,2}},\quad
 \frac{\sqrt{16^{\delta_{p,2}}p^u\delta n}}
 {\sqrt{\delta\,4^{\delta_{p,2}}n/p^u}}
 =2^{\delta_{p,2}}p^u,
\end{equation}
with identical formulas for $v$. Therefore the entire recursion reduces to the selected local density at $p$, 
the local factor $1\mp p^{-1}$, and \eqref{eq:roots}. Since $p\nmid\delta$, the number $\delta$ is a $p$-adic unit; 
the density identities below apply to the argument $n/\delta$ and its arbitrary unit part. In particular $v_p(n/\delta)=v_p(n)$.

For odd $p$, Lemma~\ref{lem:density-reciprocal-split}, equation~\eqref{eq:density-base-identity}, gives
\begin{equation}\label{eq:basic}
 p\,d_{-x^2-yz,\,p}(t)-d_{-x^2-pyz,\,p}(t)=p-1,
\end{equation}
and hence
\begin{equation}\label{eq:basic2}
 (p+1)d_{-x^2-yz,\,p}(t)
 -(1+p^{-1})d_{-x^2-pyz,\,p}(t)
 =(p+1)(1-p^{-1}).
\end{equation}

If $p^u\parallel N_1$ and $v_p(n)<u$, Lemma~\ref{lem:density-reciprocal-inert} gives the reciprocal density
\[
 \frac{1+(-1)^{v_p(n)}}2(1+p^{-1}).
\]
Multiplication by $1-p^{-1}$ and use of \eqref{eq:basic2} at $p^u n/\delta$ give the low-valuation inert recursion in Definition~\ref{def:HCNd}, with coefficient $(1+(-1)^{v_p(n)})/(2p)$ and the two terms carrying coefficients $p+1$ and $-1$. If $v_p(n)\ge u$, Lemma~\ref{lem:density-reciprocal-inert} reduces the reciprocal density to the base density at $n/(\delta p^u)$ when $u$ is even, and to $p^{-1}$ times the corresponding exponent-one inert base density when $u$ is odd. Multiplying by $1-p^{-1}$ and the second ratio in \eqref{eq:roots} gives, respectively, $p^{u-1}(p-1)$ and $p^{u-1}$, precisely the two high-valuation inert coefficients of Definition~\ref{def:HCNd}.

If $p^v\parallel N_2$ and $v_p(n)<v$, Lemma~\ref{lem:density-reciprocal-split} gives $(v_p(n)+1)(1-p^{-1})$. Multiplication by $1+p^{-1}$ and \eqref{eq:basic2} yield the low-valuation split recursion, with outer coefficient $(v_p(n)+1)/p$ and inner coefficients $p+1$ and $-1$. For $v_p(n)\ge v$, the same lemma gives
\begin{align*}
 &d_{-x^2-yz,\,p}\bigl(n/(\delta p^v)\bigr)+v(1-p^{-1})\\
 &\quad=(v+1)d_{-x^2-yz,\,p}\bigl(n/(\delta p^v)\bigr)
 -\frac vp\,d_{-x^2-pyz,\,p}\bigl(n/(\delta p^v)\bigr),
\end{align*}
where the equality follows from \eqref{eq:basic}. After multiplication by $1+p^{-1}$ and the second ratio in \eqref{eq:roots}, this gives the high-valuation split recursion with outer factor $p^{v-1}$ and inner coefficients $(v+1)(p+1)$ and $-v$. Thus every odd-prime branch agrees with Definition~\ref{def:HCNd}.

We now make explicit the calculations at $p=2$ that were implicit in the preceding argument.  Write $k=v_2(n)$, and suppose $2\nmid\delta$.  In addition to the dyadic cases of Lemmas~\ref{lem:density-reciprocal-inert} and ~\ref{lem:density-reciprocal-split}, we use
\begin{equation}\label{eq:dyadic}
 d_{-x^2-4yz,\,2}(4t)=2d_{-x^2-yz,\,2}(t),\quad
 d_{-x^2-8yz,\,2}(4t)=2d_{-x^2-2yz,\,2}(t),
\end{equation}
and the specialization of ~\eqref{eq:density-base-identity},
\begin{equation}\label{eq:twoidentity}
 2d_{-x^2-yz,\,2}(t)-d_{-x^2-2yz,\,2}(t)=1.
\end{equation}
The ratio of square-root factors in ~\eqref{gongshi5.1} is $2$ for a low-valuation recursion and $2^{u+1}$ (respectively $2^{v+1}$) for a high-valuation recursion, as follows from \eqref{eq:roots}.  Throughout the computations below, comparisons are made relative to the mass factor of the genus on the left.  In particular, deleting the entire $2$-part of $N_1$ multiplies the mass factor by $2$, while deleting the entire $2$-part of $N_2$ multiplies it by $2/3$.  Retaining a single factor $2$ on the same side leaves that mass factor unchanged.

\smallskip\noindent\emph{Case 1: $2^u\parallel N_1$ and $k<u$.}
Lemma~\ref{lem:density-reciprocal-inert} gives the selected density $3(1+(-1)^k)/4$.  After taking the square-root ratio $2$, the contribution relative to the original mass factor is
\[
 2\cdot\frac{3(1+(-1)^k)}4
 =\frac{3(1+(-1)^k)}2.
\]
The two terms on the right of the low-valuation inert recursion of Definition~\ref{def:HCNd} have coefficients $(1+(-1)^k)/4$ times $3$ and $-1$, respectively.  The corresponding mass-factor ratios are $2$ and $3$; the latter occurs because the factor $2$ is transferred from $N_1$ to $N_2$.  By \eqref{eq:dyadic}, their combined local contribution, relative to the original mass factor, is
\begin{align*}
 &\frac{1+(-1)^k}{4}\bigl[3\cdot2\cdot2d_{-x^2-yz,\,2}(t)
        -3\cdot2d_{-x^2-2yz,\,2}(t)\bigr]\\
 &\qquad=\frac{3(1+(-1)^k)}2
       \bigl[2d_{-x^2-yz,\,2}(t)-d_{-x^2-2yz,\,2}(t)\bigr]
 =\frac{3(1+(-1)^k)}2,
\end{align*}
by \eqref{eq:twoidentity}.  Here $t$ denotes the common adjusted argument at $2$ after the unit-scaling comparison already established above.

\smallskip\noindent\emph{Case 2: $2^u\parallel N_1$, $k\ge u$, and $u$ even.}
Lemma~\ref{lem:density-reciprocal-inert} gives the selected density $d_{-x^2-yz,\,2}(n/(\delta2^u))$.  Thus the left contribution relative to its mass factor is
\[
 2^{u+1}d_{-x^2-yz,\,2}\left(\frac{n}{\delta2^u}\right).
\]
The right side of Definition~\ref{def:HCNd} has coefficient $2^{u-1}$, a mass-factor ratio $2$, and a dyadic density factor $2$ from \eqref{eq:dyadic}.  Their product is $2^{u-1}\cdot2\cdot2=2^{u+1}$, as required.

\smallskip\noindent\emph{Case 3: $2^u\parallel N_1$, $k\ge u$, and $u$ odd.}
Lemma~\ref{lem:density-reciprocal-inert} gives the exact identity
\[
 d_{3\cdot2^ux^2-(y^2+yz+z^2),\,2}(n/\delta)
 =\frac14d_{3x^2-8(y^2+yz+z^2),\,2}
   \left(\frac{4n}{\delta2^u}\right).
\]
The square-root ratio is $2^{u+1}$, hence $2^{u+1}/4=2^{u-1}$.  The right-side parameter $N_1/2^{u-1}$ still contains one factor of $2$, so the mass factor does not change.  This proves precisely the odd-$u$ high-valuation coefficient $2^{u-1}$.

\smallskip\noindent\emph{Case 4: $2^v\parallel N_2$ and $k<v$.}
Lemma~\ref{lem:density-reciprocal-split} gives the selected density $(k+1)/2$.  With the square-root ratio $2$, the contribution relative to the original mass factor is $k+1$.  In Definition~\ref{def:HCNd} the right-side coefficient is $(k+1)/2$ multiplying terms with coefficients $3$ and $-1$.  The mass-factor ratios of those terms are $2/3$ and $1$, respectively.  Using \eqref{eq:dyadic} and \eqref{eq:twoidentity}, their total contribution is
\begin{align*}
 &\frac{k+1}{2}\left[3\cdot\frac23\cdot2d_{-x^2-yz,\,2}(t)
                  -2d_{-x^2-2yz,\,2}(t)\right]\\
 &\qquad=(k+1)\left[2d_{-x^2-yz,\,2}(t)-d_{-x^2-2yz,\,2}(t)\right]
 =k+1.
\end{align*}
This checks both the outer coefficient $(k+1)/2$ and the inner coefficients $3,-1$.

\smallskip\noindent\emph{Case 5: $2^v\parallel N_2$ and $k\ge v$.}
Lemma~\ref{lem:density-reciprocal-split} gives
\[
 d_{-2^vx^2-yz,\,2}(n/\delta)
 =d_{-x^2-yz,\,2}\left(\frac{n}{\delta2^v}\right)+\frac v2.
\]
The left contribution relative to its mass factor is therefore
\[
 2^{v+1}\left[d_{-x^2-yz,\,2}(t)+\frac v2\right],
 \quad t=\frac{n}{\delta2^v}.
\]
On the right, Definition~\ref{def:HCNd} has the outer coefficient $2^{v-1}$ and inner coefficients $3(v+1)$ and $-v$.  The corresponding mass-factor ratios are $2/3$ and $1$.  Thus, using \eqref{eq:dyadic}, the right contribution is
\begin{align*}
 &2^{v-1}\left[3(v+1)\cdot\frac23\cdot2d_{-x^2-yz,\,2}(t)
                   -v\cdot2d_{-x^2-2yz,\,2}(t)\right]\\
 &\qquad=2^v\left[2(v+1)d_{-x^2-yz,\,2}(t)
                   -v d_{-x^2-2yz,\,2}(t)\right]\\
 &\qquad=2^v\left[2d_{-x^2-yz,\,2}(t)+v\right]
 =2^{v+1}\left[d_{-x^2-yz,\,2}(t)+\frac v2\right],
\end{align*}
where the penultimate equality is \eqref{eq:twoidentity}.  This verifies all three coefficients $2^{v-1}$, $3(v+1)$, and $-v$.

\smallskip\noindent\emph{Case 6: $2\nmid N_1N_2\delta$.}
The final dyadic rule in Definition~\ref{def:HCNd} is
\[
 H^{(N_1,\,N_2;\,16\delta)}(n)=H^{(N_1,\,N_2;\,\delta)}(4n).
\]
The square-root ratio is $\sqrt{16\delta n}/\sqrt{4\delta n}=2$, and the first identity in \eqref{eq:dyadic} supplies the matching factor $2$ in the dyadic density. 
 There is no change in the mass factor.  Thus the coefficient of this recursion is $1$.

These six computations verify the dyadic coefficients in Definition~\ref{def:HCNd}, conditional on the local density formulas of Lemmas~\ref{lem:density-reciprocal-inert} and ~\ref{lem:density-reciprocal-split} and the dyadic identities \eqref{eq:dyadic}.  
The comparisons at $q\ne2$ use the adjusted-argument square-unit calculation proved above, not merely the fact that the local similarity factors are units.
  Together with the odd-prime verification, this establishes the required recursions of the local Euler products.  Starting from the initial case $\delta=1$,
   induction on the number of selected prime factors identifies \eqref{eq:formalproduct} with $H^{(N_1,\,N_2;\,\delta)}(n)$. 
   The induction is independent of the order of selection: by Remarks \ref{rem:phiN} and \ref{remark11} the correspondences
    and the defining recursions commute, while earlier selections change the local represented argument at a different prime by a square unit, 
    as established in the preceding local comparison.
Using the identification with the actual genus above, we conclude that
\[
 \sum_{f\in G_{4N,\,16N^2/\delta,\,N_1}}
 \frac{R_f(n)}{|\operatorname{Aut}(f)|}
 =2^{-e(N)-1}H^{(N_1,\,N_2;\,\delta)}(n).
\]
Together with the case $n=0$, this proves the first assertion.

Assume $2\parallel N$ and $\delta'\parallel N/2$. By the change of variables defining $\lambda_4$,
\[
 R_{\lambda_4(f)}(n)=R_f(4n).
\]
Lemma~\ref{prop:bijection} shows that $\lambda_4$ gives a bijection of the relevant classes and preserves their automorphism groups; 
it also commutes with the odd prime correspondences. Hence Definition~\ref{def:genus111} and the first assertion give
\begin{align*}
 \sum_{f\in G_{2N,\,N^2/\delta',\,N_1}}
 \frac{R_f(n)}{|\operatorname{Aut}(f)|}
 &=\sum_{f\in G_{4N,\,16N^2/\delta',\,N_1}}
 \frac{R_f(4n)}{|\operatorname{Aut}(f)|}\\
 &=2^{-e(N)-1}H^{(N_1,\,N_2;\,\delta')}(4n),
\end{align*}
including $n=0$. This completes the proof.
\end{proof}

\subsection{Proof of Theorem~\ref{the:classnum}}\label{sec:6.3}
	
Put $N=N_1N_2$. We start from
\eqref{eq:type-normalizer-count}. The three steps are to characterize
primitive normalizing elements, express their counts as ordinary
representation numbers after Lehman's correspondence, and evaluate
the resulting genus averages. The first two steps extend
\cite[Lemmas~2.6 and~3.13]{LZ26}; the final local comparison extends
\cite[Lemma~5.3]{LZ26}. The normalizer calculation needed for these
steps has already been proved in Lemma~\ref{lemmafornormalizer}.

We write $c_{ij}(\alpha_p)$ for the $(i,j)$-entry of $\alpha_p$ in the
local matrix model of Definition~\ref{def3.6}.
\begin{lemma}\label{lemmaforprincipalideal}
Let $0\ne \alpha\in\mathcal O$, where $\mathcal O$ has level $(N_1,N_2)$.
Then $\alpha$ is primitive and $\mathcal O\alpha$ is a two-sided
$\mathcal O$-ideal if and only if, writing $n=\n(\alpha)$, all the
following conditions hold:
\begin{enumerate}[label=\textup{(\arabic*)}]
 \item $n\parallel N$;
 \item $n\mid\tr(\alpha)$;
 \item if $p^{2u+1}\parallel n$ and $p\mid N_1$, then
       $c_{11}(\alpha_p)\in p^{2u+1}R_p$;
 \item if $p^{2u}\parallel n$, $p\mid N_1$, and $p\ne2$, then
       $c_{12}(\alpha_p)+\epsilon_pc_{21}(\alpha_p)\in p^{2u}\mathbb Z_p$;
 \item if $2^{2u}\parallel n$ and $2\mid N_1$, then
       $c_{11}(\alpha_2)+c_{12}(\alpha_2)+c_{21}(\alpha_2)
                                      \in2^{2u}\mathbb Z_2$;
 \item if $p^v\parallel n$ and $p\mid N_2$, then
       $c_{11}(\alpha_p)\in p^v\mathbb Z_p$.
\end{enumerate}
Equivalently, at every prime $p$,
\[
 \alpha_p\in
 \begin{cases}
 \mathcal O_p^\times, &\text{ if }\,\,p\nmid n,\\
 w_p\mathcal O_p^\times, &\text{ if }\,\,p\mid n,
 \end{cases}
\]
where $w_p$ is as in Lemma~{\rm\ref{lemmafornormalizer}}.
\end{lemma}

\begin{proof}
The ideal $\mathcal O \alpha$ is
two-sided if and only if $\alpha\mathcal O \alpha^{-1}\subseteq\mathcal O$.
Since conjugate orders have the same discriminant, this inclusion is
necessarily an equality.  Thus the left-hand condition of the lemma
means that $\alpha$ is a primitive element of $N_Q(\mathcal O)$.

Membership in the normalizer and primitivity can both be checked locally.
The local forms of $w_p$ in Lemma~{\rm\ref{lemmafornormalizer}} show that $w_p\in\mathcal O_p\setminus p\mathcal O_p$ for
$p\mid N$, and right multiplication by $\mathcal O_p^\times$ preserves
both $\mathcal O_p$ and $p\mathcal O_p$.  In the description of
Lemma~{\rm\ref{lemmafornormalizer}}, integrality therefore forces the scalar factor to have
nonnegative valuation, and primitivity forces that valuation to be zero.
The primitive integral normalizing elements are consequently
$\mathcal O_p^\times$ for $p\nmid N$, and
$\mathcal O_p^\times\cup w_p\mathcal O_p^\times$ for $p\mid N$.
Their norm valuations are respectively $0$ and $v_p(N)$.  It follows that
the left-hand condition is equivalent to $n\parallel N$ together with
\begin{equation}\label{eq:lemma64-local-cosets}
 \alpha_p\in
 \begin{cases}
  \mathcal O_p^\times,&p\nmid n,\\
  w_p\mathcal O_p^\times,&p\mid n.
 \end{cases}
\end{equation}

Assume now that $n\parallel N$.  For $p\nmid n$, the first condition
in \eqref{eq:lemma64-local-cosets} is automatic.  For $p\mid n$, put
$\beta_p=w_p^{-1}\alpha_p$.  Since
\[
 v_p(\mathrm n(\beta_p))=v_p(n)-v_p(\mathrm n(w_p))=0,
\]
we have
\begin{equation}\label{eq:lemma64-integrality}
 \alpha_p\in w_p\mathcal O_p^\times
 \quad\Longleftrightarrow\quad \beta_p\in\mathcal O_p.
\end{equation}
Indeed, an element of $\mathcal O_p$ with unit reduced norm has inverse
$\bar \beta_p/\mathrm n(\beta_p)\in\mathcal O_p$.  We translate
\eqref{eq:lemma64-integrality} using Definition~\ref{def3.6}.

If $p^{2u+1}\parallel n$ and $p\mid N_1$, direct multiplication gives
\[
 \beta_p\in\mathcal O_p
 \quad\Longleftrightarrow\quad
 c_{11}(\alpha_p)\in p^{2u+1}R_p.
\]
This is condition~(3), and it also implies
$\operatorname{tr}(\alpha_p)\in p^{2u+1}\mathbb Z_p$.
In the remaining, split cases, write
$\alpha_p=\left(\begin{smallmatrix}a&b\\c&d\end{smallmatrix}\right)$.
If $p^v\parallel n$ and $p\mid N_2$, then
\[
 \beta_p=\begin{pmatrix}c/p^v&d/p^v\\a&b\end{pmatrix}\in\mathcal O_p
 \quad\Longleftrightarrow\quad a,d\in p^v\mathbb Z_p.
\]
This is precisely condition~(6) together with the local trace condition.

Suppose next that $p^{2u}\parallel n$, $p\mid N_1$, and $p\ne2$.
Membership of $\alpha_p$ in $\mathcal O_p$ gives
$a-d,b-\epsilon_pc\in p^u\mathbb Z_p$.  We claim that
\[
 \beta_p\in\mathcal O_p
 \quad\Longleftrightarrow\quad
 a+d,\ b+\epsilon_pc\in p^{2u}\mathbb Z_p.
\]
The congruences on the right, combined with the preceding ones, imply
$a,b,c,d\in p^u\mathbb Z_p$, since $2$ is a unit.  Hence
\[
 \beta_p=p^{-u}\begin{pmatrix}
  -c&-d\\ \epsilon_p^{-1}a&\epsilon_p^{-1}b
 \end{pmatrix}
\]
is integral, and its defining congruences are exactly
\[
 c_{11}(\beta_p)-c_{22}(\beta_p)
   =-\frac{b+\epsilon_pc}{\epsilon_pp^u}\in p^u\mathbb Z_p,
 \quad
 c_{12}(\beta_p)-\epsilon_pc_{21}(\beta_p)
   =-\frac{a+d}{p^u}\in p^u\mathbb Z_p.
\]
These identities also prove the reverse implication.  Thus this case
gives condition~(4) together with the local trace condition.

Finally, suppose that $2^{2u}\parallel n$ and $2\mid N_1$.
We have $b\equiv c$ and $a-d\equiv b\pmod{2^u}$.  Here the required
equivalence is
\[
 \beta_2\in\mathcal O_2
 \quad\Longleftrightarrow\quad
 a+d,\ a+b+c\in2^{2u}\mathbb Z_2.
\]
Under the congruences on the right,
\[
 d\equiv-a,\quad b\equiv c\equiv2a,\quad
 0\equiv a+b+c\equiv5a\pmod{2^u}.
\]
Since $5$ is a unit, all four entries belong to $2^u\mathbb Z_2$, so
\[
 \beta_2=2^{-u}\begin{pmatrix}-c&-d\\a&b\end{pmatrix}
\]
is integral.  Its defining congruences are
\[
 \begin{aligned}
  c_{12}(\beta_2)-c_{21}(\beta_2)
    &=-\frac{a+d}{2^u}\in2^u\mathbb Z_2,\\
  c_{11}(\beta_2)-c_{22}(\beta_2)-c_{12}(\beta_2)
    &=\frac{(a+d)-(a+b+c)}{2^u}\in2^u\mathbb Z_2.
 \end{aligned}
\]
Again these identities also give the reverse implication.  This is
condition~(5) together with the local trace condition.

Combining the trace conditions over $p\mid n$ gives
$n\mid\operatorname{tr}(\alpha)$.  Thus, under condition~(1), the local
coset description \eqref{eq:lemma64-local-cosets} is equivalent to
conditions~(2)--(6).  By the first part of the proof, this proves both
assertions of the lemma.
\end{proof}

For $n\parallel N$, $n\geq4$, define
\[
\mathcal R_{\mathcal O}(n)\!:\,=\{\alpha\in\mathcal O:
 \alpha\text{ is primitive},\ \alpha\in N_Q(\mathcal O),\ \n(\alpha)=n\}.
\]
For $f=f_{S^0}$ we set $R'_f(4n)=\rho'_{\mathcal O}(n,0)$.
Thus the prime $'$ denotes the full primitive-normalizer restriction
of Lemma~\ref{lemmaforprincipalideal}, including its dyadic condition.

\begin{lemma}\label{lemmaforxingshu1}
Every element of $\mathcal R_{\mathcal O}(n)$ has trace zero, and
\[
 |\mathcal R_{\mathcal O}(n)|
 =R'_{f_{S^0}}(4n)
 =R_{\phi_n(f_{S^0})}(4^{\delta_{2\nmid n}}).
\]
Consequently,
\[
 \sum_{f\in G_{4N,\,16N^2,\,N_1}}
       \frac{R'_f(4n)}{|\Aut(f)|}
 =
 \sum_{f\in G_{4N,\,16N^2,\,N_1}}
       \frac{R_{\phi_n(f)}(4^{\delta_{2\nmid n}})}
            {|\Aut(\phi_n(f))|}.
\]
\end{lemma}

\begin{proof}
Let \(\alpha\in\mathcal R_{\mathcal O}(n)\). Lemma~\ref{lemmaforprincipalideal} gives
\(n\mid\operatorname{tr}(\alpha)\). The definiteness of $Q$ implies 
\(
0\leq \n\left(\alpha-{\operatorname{tr}(\alpha)}/2\right)
=n-{\operatorname{tr}(\alpha)^2}/{4},
\)
i.e., \(|\operatorname{tr}(\alpha)|\leq2\sqrt n\).
For \(n>4\), the inequality \(2\sqrt n<n\) implies
\(\operatorname{tr}(\alpha)=0\). If \(n=4\), a nonzero trace would
necessarily be \(\pm4\). So the preceding inequality by definiteness would give
\(\alpha=\operatorname{tr}(\alpha)/2=\pm2\), which are not primitive. Thus
$$
\mathcal R_{\mathcal O}(n)\subseteq\mathcal O^0.
$$

\medskip
Let \(B:Q^0\times Q^0\longrightarrow\mathbb Q\) be the symmetric
bilinear form associated with the reduced norm $\n$, given by
$$
B(\alpha,\beta)\!:\,=\n(\alpha+\beta)-\n(\alpha)-\n(\beta)
=\operatorname{tr}(\alpha\bar \beta).
$$
In particular, \(B(\alpha,\alpha)=2\n(\alpha)\), and \(B\) is integral valued on
\(S^0\times S^0\). 
Define the sublattice
\[
L_n\!:\,=\{\alpha\in S^0:B(\alpha,\beta)\in 4^{\delta_{2\mid n}}n\mathbb Z\text{ for every }\beta\in S^0\}.
\]
We claim that
\begin{equation}\label{lem65:lehman-lattice}
\phi_n(f_{S^0})
\,\,\,\text{is integrally equivalent to the quadratic form of}\,\,\,
(L_n,\n/(4^{\delta_{2\mid n}}n)).
\end{equation}
This identifies the integral equivalence class associated with the
specific sublattice \(L_n\), which is needed for the representation count.

First consider one prime \(p\mid n\), and write \(p^g\parallel n\).
Since \(n\parallel N\), we have \(g=v_p(N)\). Put
\[
\ell_p\!:\,=
\begin{cases}
p^g,&p\neq2,\\
2^{g+2},&p=2,
\end{cases}
\quad
\kappa_p\!:\,=
\begin{cases}
p^g,&p\neq2,\\
2^{g+1},&p=2.
\end{cases}
\]
In the integral normal form used to define Lehman's correspondence in
Section~\ref{section2.2}, the quadratic form under consideration has a representative
\begin{equation}\label{quadratic form}
f(x,y,z)
=\ell_pax^2+\ell_pby^2+cz^2
+\ell_pryz+\ell_psxz+\ell_ptxy,
\quad p\nmid c.
\end{equation}
Here \(a,b,c,r,s,t\) are integer coefficients. For odd \(p\), the
exponents of \(p\) in the level and discriminant are respectively \(g\) and \(2g\).
For \(p=2\), these exponents are \(g+2\) and \(2(g+2)\).
In these cases, the coordinate formula for Lehman's correspondence is
\begin{equation}\label{fangcheng}
\phi_p(f)(x,y,z)
=\frac{f(x,y,\kappa_pz)}{\ell_p}
=ax^2+by^2+\frac{\kappa_p^2}{\ell_p}cz^2
  +\kappa_pryz+\kappa_psxz+txy.
\end{equation}
In particular, at \(p=2\), the coefficient of \(z^2\) is \(2^gc\)
and the two mixed coefficients involving \(z\) are \(2^{g+1}r\)
and \(2^{g+1}s\).

Let \(B_f(\boldsymbol{v},\boldsymbol{w})=f(\boldsymbol{v}+\boldsymbol{w})-f(\boldsymbol{v})-f(\boldsymbol{w})\) be the symmetric bilinear form associated with \(f\), and let \(e_1=(1,0,0),e_2=(0,1,0),e_3=(0,0,1)\) be the basis used in
~\eqref{quadratic form}.
For \(\boldsymbol{v}=(x,y,z)\), both \(B_f(\boldsymbol{v},e_1)\) and \(B_f(\boldsymbol{v},e_2)\) are
automatically divisible by \(\ell_p\), whereas
\[
B_f(\boldsymbol{v},e_3)=2cz+\ell_p(ry+sx).
\]
It follows that
\begin{equation}\label{lem65:coordinate-condition}
B_f(\boldsymbol{v},\mathbb Z^3)\subseteq \ell_p\mathbb Z
\quad\Longleftrightarrow\quad
\ell_p\mid2cz
\quad\Longleftrightarrow\quad
\kappa_p\mid z.
\end{equation}
The last equivalence uses \(p\nmid c\).
Thus the integral substitution
\[
(x,y,z)\longmapsto(x,y,\kappa_pz)
\]
has image exactly the sublattice defined by the divisibility condition
for \(B_f\) in \eqref{lem65:coordinate-condition}.
Equation~\eqref{fangcheng} identifies the form on this
sublattice, divided by \(\ell_p\), with \(\phi_p(f)\).

Apply these substitutions successively at all primes dividing \(n\),
and regard their images as sublattices of the original rational space
\(Q^0\). Each substitution has determinant a power of its selected prime
and is invertible at every other prime. At a newly selected prime \(p\),
all preceding lattice modifications are therefore invertible, and the
preceding norm-scaling factors are \(p\)-adic units. Consequently, the
new condition at \(p\) is precisely
\[
B(x,L_p)\subseteq \ell_p\mathbb Z_p
\]
for the original symmetric bilinear form \(B\). Previously imposed local conditions
remain unchanged. The final image thus has the same localization as
\(L_n\) at every prime, since
\[
\prod_{p\mid n}\ell_p=4^{\delta_{2\mid n}}n.
\]
Equality of all localizations, inside the same rational space, implies
equality of the two integral lattices. The total norm-scaling factor is
\(4^{\delta_{2\mid n}}n\), proving \eqref{lem65:lehman-lattice}.

In particular, \eqref{lem65:lehman-lattice} yields
\begin{equation}\label{lem65:transformed-count}
R_{\phi_n(f_{S^0})}\bigl(4^{\delta_{2\nmid n}}\bigr)
=
\#\{x\in L_n:\n(x)=4n\},
\end{equation}
because \(4n/(4^{\delta_{2\mid n}}n)=4^{\delta_{2\nmid n}}\).

\medskip
The map \(\alpha\mapsto2\alpha\) gives a bijection
\begin{equation}\label{shuangshe}
\{\alpha\in\mathcal O^0:\n(\alpha)=n\}
\xrightarrow{\ \sim\ }
\{\alpha\in L:\n(\alpha)=4n\}.
\end{equation}
Clearly, the map is well-defined and injective. For surjectivity, let
\(\alpha\in L\) satisfy \(\n(\alpha)=4n\), and write
\[
\alpha=a+2\gamma,\quad a\in\mathbb Z,\quad \gamma\in\mathcal O.
\]
The equation \(\operatorname{tr}(\alpha)=0\) gives
\(a=-\operatorname{tr}(\gamma)\), whence
\[
\n(\alpha)
=a^2+2a\operatorname{tr}(\gamma)+4\n(\gamma)
=-a^2+4\n(\gamma).
\]
Thus \(4\mid a^2\), and \(a\) is even. Therefore
\[
\alpha/2=a/2+\gamma\in\mathcal O^0,
\quad \n(\alpha/2)=n.
\]

It remains to show that, for \(\alpha\in\mathcal O^0\) of norm \(n\),
\begin{equation}\label{lem65:restriction-equivalence}
2\alpha\in L_n
\quad\Longleftrightarrow\quad
\alpha\in\mathcal R_{\mathcal O}(n).
\end{equation}

\medskip
Fix \(\alpha\in\mathcal O^0\) with \(\n(\alpha)=n\).
If \(p\nmid n\), then \(4^{\delta_{2\mid n}}n\) is a \(p\)-adic unit, so the divisibility
condition for \(B\) imposes no restriction. Moreover,
\(\alpha_p^{-1}=\bar\alpha_p/n\in\mathcal O_p\), and hence
\(\alpha_p\in\mathcal O_p^\times\), as required by Lemma~\ref{lemmaforprincipalideal}.

Now let \(p^g\parallel n\). We must check
\begin{equation}\label{lem65:local-pairing}
B(2\alpha_p,L_p)\subseteq \ell_p\mathbb Z_p.
\end{equation}
For a trace zero element \(\beta\), we have
\[
B(2\alpha_p,\beta)=-2\operatorname{tr}(\alpha_p\beta).
\]
Write \(\nu_p=1\) for odd \(p\), and \(\nu_2=2\); then
\(\ell_p=\nu_p^2p^g\). We use the local orders of Definition~\ref{def3.6}.
The matrices \(E_{ij}\) below are the standard matrix units.
We divide the argument into the following four cases. 

\smallskip
\emph{\rm (i) The Eichler case: \(p\mid N_2\).}
Write
\[
\alpha_p=
\begin{pmatrix}a&b\\p^gc&-a\end{pmatrix},
\quad a,b,c\in\mathbb Z_p.
\]
A basis of \(L_p\) is
\[
\begin{pmatrix}1&0\\0&-1\end{pmatrix},
\quad \nu_pE_{12},\quad \nu_pp^gE_{21}.
\]
Evaluating \(B(2\alpha_p,\cdot)\) on these basis vectors yields, respectively,
\[
-4a,\quad -2\nu_pp^gc,\quad -2\nu_pp^gb.
\]
The last two are automatically divisible by \(\ell_p\). Therefore
\eqref{lem65:local-pairing} is equivalent to
\[
a=c_{11}(\alpha_p)\in p^g\mathbb Z_p.
\]
This is condition~(6) of Lemma~\ref{lemmaforprincipalideal}, including \(p=2\).

\smallskip
\emph{\rm (ii) The division case: \(p\mid N_1\) and \(g=2u+1\).}
Let \(K_p/\mathbb Q_p\) be the unramified quadratic extension with
ring of integers \(R_p\). Choose a trace zero unit \(\theta\) such that
\[
\theta^2=
\begin{cases}
\epsilon_p,&p\neq2,\\
-3,&p=2.
\end{cases}
\]
Then
\[
R_p=\mathbb Z_p[\theta]\quad(p\neq2),
\quad
R_2=\mathbb Z_2[(1+\theta)/2],
\]
and the trace zero part of \(R_p\) is \(\mathbb Z_p\theta\).
We may write
\[
\alpha_p=
\begin{pmatrix}
a&p^ub\\p^{u+1}\bar b&\bar a
\end{pmatrix},
\quad a=x\theta,\quad x\in\mathbb Z_p,\quad b\in R_p.
\]
Set
\[
H_\theta=\begin{pmatrix}\theta&0\\0&-\theta\end{pmatrix},
\quad
U(w)=\begin{pmatrix}0&p^uw\\p^{u+1}\bar w&0\end{pmatrix}.
\]
If \(w_1,w_2\) is a \(\mathbb Z_p\)-basis of \(R_p\), then
\(H_\theta,\nu_pU(w_1),\nu_pU(w_2)\) form a basis of \(L_p\).
Direct calculation gives
\begin{align*}
B(2\alpha_p,H_\theta)&=-4\theta^2x,\\
B(2\alpha_p,\nu_pU(w))
&=-2\nu_pp^{2u+1}
  \operatorname{Tr}_{K_p/\mathbb Q_p}(b\bar w).
\end{align*}
The second value is automatically divisible by \(\ell_p\).
Since \(\theta^2\) is a unit, the first value has the required
divisibility if and only if \(x\in p^{2u+1}\mathbb Z_p\), equivalently,
\[
c_{11}(\alpha_p)=a\in p^{2u+1}R_p.
\]
This is condition~(3) of Lemma~\ref{lemmaforprincipalideal}, again including \(p=2\).

\smallskip
\emph{\rm (iii) The split residually inert case at an odd prime:
\(p\mid N_1\), \(p\neq2\), and \(g=2u\).}
Write
\[
\alpha_p=\begin{pmatrix}a&b\\c&-a\end{pmatrix}.
\]
Membership in \(\mathcal O_p^0\) gives
\[
a\in p^u\mathbb Z_p,\quad b-\epsilon_pc\in p^u\mathbb Z_p.
\]
Since \(p\neq2\), we have \(L_p=\mathcal O_p^0\). By Lemma~\ref{lemma3.3},
a basis is
\[
H=\begin{pmatrix}p^u&0\\0&-p^u\end{pmatrix},
\quad
U=\begin{pmatrix}0&\epsilon_p\\1&0\end{pmatrix},
\quad
V=\begin{pmatrix}0&\epsilon_pp^u\\-p^u&0\end{pmatrix}.
\]
Direct evaluation of the bilinear form gives
\begin{align*}
B(2\alpha_p,H)&=-4p^ua,\\
B(2\alpha_p,U)&=-2(b+\epsilon_pc),\\
B(2\alpha_p,V)&=2p^u(b-\epsilon_pc).
\end{align*}
The first and third are already divisible by \(p^{2u}\). Thus
\eqref{lem65:local-pairing} is equivalent to
\[
b+\epsilon_pc\in p^{2u}\mathbb Z_p,
\]
which is condition~(4) of Lemma~\ref{lemmaforprincipalideal}.

\smallskip
\emph{\rm (iv) The split residually inert case at \(2\):
\(2\mid N_1\) and \(g=2u\).}
Write
\[
\alpha_2=\begin{pmatrix}a&b\\c&-a\end{pmatrix}.
\]
The defining congruences for \(\mathcal O_2\), together with trace zero,
give
\[
b\equiv c\equiv2a\pmod{2^u}.
\]
In the original matrix expression, a basis of \(L_2\) is
\[
H=\begin{pmatrix}1&2\\2&-1\end{pmatrix},
\quad
U=\begin{pmatrix}0&2^{u+1}\\0&0\end{pmatrix},
\quad
V=\begin{pmatrix}0&0\\2^{u+1}&0\end{pmatrix}.
\]
Hence
\begin{align*}
B(2\alpha_2,H)&=-4(a+b+c),\\
B(2\alpha_2,U)&=-2^{u+2}c,\\
B(2\alpha_2,V)&=-2^{u+2}b.
\end{align*}
Since \(\ell_2=2^{2u+2}\), the divisibility condition
\eqref{lem65:local-pairing} is equivalent to
\[
a+b+c\in2^{2u}\mathbb Z_2,
\quad b,c\in2^u\mathbb Z_2.
\]
The last two conditions follow from the first one, which is exactly condition~(5) of Lemma~\ref{lemmaforprincipalideal}. Indeed, since 
\(a+b+c\equiv5a\pmod{2^u}\), so the first condition implies
\(5a\in2^u\mathbb Z_2\), hence \(a\in2^u\mathbb Z_2\), and then 
\(b,c\in2^u\mathbb Z_2\). 

The computations above prove that, under the standing assumptions
\(n\parallel N\), \(\n(\alpha)=n\), and \(\operatorname{tr}(\alpha)=0\),
membership \(2\alpha\in L_n\) is equivalent to all the remaining
conditions of Lemma~\ref{lemmaforprincipalideal}. That lemma identifies these conditions with
primitivity and the normalizer condition, proving
\eqref{lem65:restriction-equivalence}.

\medskip
Combining \eqref{shuangshe} and
\eqref{lem65:restriction-equivalence}, we obtain the bijection
\[
\mathcal R_{\mathcal O}(n)
\xrightarrow[\alpha\mapsto2\alpha]{\ \sim\ }
\{x\in L_n:\n(x)=4n\}.
\]
Together with \eqref{lem65:transformed-count}, this gives
\[
|\mathcal R_{\mathcal O}(n)|
=
R_{\phi_n(f_{S^0})}\bigl(4^{\delta_{2\nmid n}}\bigr).
\]
Thus by the definition of the restricted count we get
\[
|\mathcal R_{\mathcal O}(n)|
=
\rho'_{\mathcal O}(n,0)
=
R'_{f_{S^0}}(4n).
\]
This proves the asserted pointwise identity.

Finally, Corollary~\ref{cor:tqf} identifies every class in
\(G_{4N,\,16N^2,\,N_1}\) with the norm form \(f_{S^0}\) of an order
of level \((N_1,N_2)\). By Lemma~\ref{prop:bijection},
\[
|\operatorname{Aut}(f)|
=
|\operatorname{Aut}(\phi_n(f))|.
\]
Dividing the pointwise identities by these equal automorphism orders
and summing over the classes in the genus proves the weighted identity.
\end{proof}

\begin{lemma}\label{lemma5.5}
Let $n\parallel N$ and $n\ge4$. With the notations of Definition~{\rm\ref{def:HCNd}} and Theorem~{\rm\ref{the:classnum}}, one has
\[
 H^{(N_1,\,N_2;\,16^{\delta_{2\mid n}}n)}
       (4^{\delta_{2\nmid n}})
 =
 H^{(N_1,N_2)}(4n)
 \prod_{p\mid n}
 \frac{1-\left(\dfrac{\Delta(-4n)}p\right)/p}
      {B_p(n;N_1,N_2)C_p(n)}.
\]
\end{lemma}

\begin{proof}
Put 
\[
 G_1=G_{4N,\,16N^2,\,N_1},\quad
 G_n=\phi_n(G)=G_{4N,\,16N^2/(16^{\delta_{2\mid n}} n),\,N_1},
\]
and write 
\[
 A_{G_1}(t)
 =\sum_{f\in{G_1}}
   \frac{R_f(t)}{|\operatorname{Aut}(f)|},\quad
A_{G_n}(t)
 =\sum_{f\in{G_n}}
   \frac{R_f(t)}{|\operatorname{Aut}(f)|}.   
\]
The Lehman correspondence preserves automorphism group orders, so $G_1$ and
$G_n$ have the same mass.  Since
$d_{G_1}=16N^2$ and $d_{G_n}=16N^2/(16^{\delta_{2\mid n}} n)$, the ratio of their
archimedean factors in the Siegel-Weil formula ~\eqref{gongshi5.1} is
\[
 \frac{\sqrt{{4^{\delta_{2\nmid n}}}/d_{G_n}}}{\sqrt{4n/d_{G_1}}}
 =\sqrt{\frac{4^{1-{\delta_{2\mid n}}}16^{\delta_{2\mid n}} n}{4n}}=2^{\delta_{2\mid n}}.
\]
We first justify the cancellation of the local factors at primes $\ell\nmid n$, including $\ell=2$ when $n$ is odd.
 The local form of $G_n$ is integrally equivalent over $\mathbb Z_\ell$ to $16^{\delta_{2\mid n}}n$ times that of $G_1$, 
 by the off-prime similarity property of Lehman's correspondence (Remark~\ref{rem:phiN}, applied successively). Consequently, for the density normalization in ~\eqref{gongshi5.1},
\[
 d_{G_n,\,\ell}({4^{\delta_{2\nmid n}}})=d_{G_1,\,\ell}({4^{\delta_{2\nmid n}}}/(16^{\delta_{2\mid n}}n)),
\]
which follows from the defining congruence counts, since $16^{\delta_{2\mid n}}n\in\mathbb Z_\ell^\times$. 
Moreover,
\[
 \frac{{4^{\delta_{2\nmid n}}}/(16^{\delta_{2\mid n}}n)}{4n}=\frac{{4^{\delta_{2\nmid n}}}}{4n\cdot 16^{\delta_{2\mid n}}n}
 =\begin{cases}
  1/n^2,&2\nmid n,\\[2pt]
  1/(64n^2),&2\mid n.
 \end{cases}
\]
Both are squares of $\ell$-adic units whenever $\ell\nmid n$ and, in the second case, $\ell\ne2$ (which is automatic). Hence a unit change of variables gives
\begin{equation}\label{eq:off}
 d_{G_n,\,\ell}({4^{\delta_{2\nmid n}}})=d_{G_1,\,\ell}(4n)\quad(\ell\nmid n).
\end{equation}
In particular, when $n$ is odd and $\ell=2$, the scaling factor is $16^{\delta_{2\mid n}}n=n$, which \emph{is} a $2$-adic unit; the represented arguments differ by the square unit $n^{-2}$. 
Consequently, on setting
\[
 K_p(n)=2^{\delta_{p,2}}
        \frac{d_{G_n,p}({4^{\delta_{2\nmid n}}})}{d_{{G_1},p}(4n)}
 \quad(p\mid n),
\]
formula~\eqref{gongshi5.1} gives
\begin{equation}\label{eq:lemma66-genus-comparison}
 A_{G_n}({4^{\delta_{2\nmid n}}})=A_{G_1}(4n)\prod_{p\mid n}K_p(n).
\end{equation}
The positivity of the denominators follows from the local formulas below.
It remains to show that, with $\chi_p(n)=\left(\frac{\Delta(-4n)}p\right)$, 
\begin{equation}\label{eq:lemma66-local-factor}
 K_p(n)=\frac{1-\chi_p(n)/p}{B_p(n;N_1,N_2)C_p(n)}.
\end{equation}

The local comparisons in ~\cite[Lemma~5.3]{LZ26} prove
\eqref{eq:lemma66-local-factor} when $p\mid N_2$, or when $p\mid N_1$
and $v_p(n)$ is odd, with the archimedean factor already included in $K_2$.
We therefore consider $p\mid N_1$ and $n=p^{2u}m$, where $u\geq1$ and
$p\nmid m$.

If $p$ is odd, Lemmas~\ref{lem:density-even-inert} and~\ref{lem:density-reciprocal-inert} give
\[
 d_{{G_1},\,p}(4n)=p^u(1+\chi_p(n)/p),\quad
 d_{G_n,\,p}({4^{\delta_{2\nmid n}}})=1+p^{-1}.
\]
Here $\chi_p(n)=\left(\frac{-m}p\right)\in\{1,-1\}$; the factor $4$ is a
unit square and does not affect these densities.  Hence
\[
 K_p(n)=\frac{1+p^{-1}}{p^u(1+\chi_p(n)/p)}
       =\frac{1-\chi_p(n)/p}{p^{u-1}(p-1)}.
\]
This is \eqref{eq:lemma66-local-factor}, since
$B_p(n;N_1,N_2)=p^{u-1}(p-1)$ and $C_p(n)=1$.

Suppose now that $p=2$.  After removing the unit contributed by the
other prime factors of $n$, the transformed local equation is
\[
 3\cdot2^{2u}x^2-(y^2+yz+z^2)=m.
\]
Thus Lemma~\ref{lem:density-reciprocal-inert} gives $d_{G_n,\,2}(1)=3/2$.  Moreover, Lemma~\ref{lem:density-even-inert} and
~\cite[Lemmas~B.4--B.5]{LZ26} yield
\[
 d_{{G_1},\,2}(4n)
 =2^u d_{-x^2-4yz,\,2}(4m)
 =2^{u+1}d_{-x^2-yz,\,2}(m).
\]
The resulting densities, together with the corresponding values of
$\chi_2(n)$ and $C_2(n)$, are as follows:
\[
 \renewcommand{\arraystretch}{1.2}
 \begin{array}{c|ccc}
  m & d_{{G_1},\,2}(4n) & \chi_2(n) & C_2(n)\\ \hline
  7\pmod 8 & 3\cdot2^u     & 1  & 1\\
  3\pmod 8 & 2^{u+1}       & -1 & 2\\
  1\pmod 4 & 3\cdot2^{u-1} & 0  & 1
 \end{array}
\]
Since $B_2(n;N_1,N_2)=2^{u-1}$, each row gives
\[
 K_2(n)=\frac{2(3/2)}{d_{{G_1},\,2}(4n)}
       =\frac{1-\chi_2(n)/2}{2^{u-1}C_2(n)},
\]
proving \eqref{eq:lemma66-local-factor} also in this case.

All selected local denominators are positive.  If a common density at
some $\ell\nmid n$ vanishes, both weighted representation numbers vanish,
so \eqref{eq:lemma66-genus-comparison} remains valid.  Finally, substitute
\eqref{eq:lemma66-local-factor} into
\eqref{eq:lemma66-genus-comparison} and use Theorem~\ref{the:repnum}:
\[
 \begin{aligned}
  A_{G_1}(4n)&=2^{-e(N)-1}H^{(N_1,N_2)}(4n),\\
  A_{G_n}({4^{\delta_{2\nmid n}}})&=2^{-e(N)-1}
    H^{(N_1,N_2;\,16^{\delta_{2\mid n}} n)}({4^{\delta_{2\nmid n}}}).
 \end{aligned}
\]
Cancelling the common factor proves the asserted identity.
\end{proof}

\begin{proof}[Proof of Theorem~{\upshape\ref{the:classnum}}]
By Corollary~\ref{cor:tqf} and the class bijections in Definition~\ref{def:genus111} we have
\[
 |G_{4N,16N^2/\delta,N_1}|=|G_{4N,16N^2,N_1}|=T_{N_1,N_2}
 \quad(\delta\parallel16N).
\]
When $2\parallel N$, the same bijections also give
$|G_{2N,N^2/\delta',N_1}|=T_{N_1,N_2}$ for every $\delta'\parallel N/2$.
It therefore suffices to evaluate $T_{N_1,N_2}$.

Let $\mathcal O_1,\ldots,\mathcal O_{T_{N_1,N_2}}$ represent the types of orders of
level $(N_1,N_2)$.  The primitive normalizer counting identity~\eqref{eq:type-normalizer-count} is
\[
 T_{N_1,N_2}=\frac12\sum_{t\parallel N}
       \sum_{\substack{t\mid r\\r^2\leq4t}}
       \sum_{\mu=1}^{T_{N_1,N_2}}
       \frac{\rho'_{\mathcal O_\mu}(t,r)}
            {|\operatorname{Aut}(\mathcal O_\mu)|}.
\]

For $1\le t\le 3$, consider an element
$\alpha\in\mathcal O_\mu$ satisfying
$\operatorname{n}(\alpha)=t$ and
$\operatorname{tr}(\alpha)=r$, where
$t\parallel N$ and $t\mid r$.
The norm and trace conditions imply the local normalizer
conditions in Lemma~\ref{lemmaforprincipalideal}, as in~\cite[Remark~3]{LZ26}.
Thus $\alpha$ normalizes $\mathcal O_\mu$.
Moreover, $\alpha$ is primitive. Indeed, if
$\alpha\in p\mathcal O_\mu$ for some prime $p$, then
$\alpha=p\beta$ for some $\beta\in\mathcal O_\mu$, whence
\[
 t=\operatorname{n}(\alpha)
  =p^2\operatorname{n}(\beta).
\]
Since $\operatorname{n}(\beta)\in\mathbb Z$, this
contradicts the squarefreeness of $t\in\{1,2,3\}$.
Consequently, every element counted by
$\rho_{\mathcal O_\mu}(t,r)$ is a primitive normalizing
element, and hence
\[
 \rho'_{\mathcal O_\mu}(t,r)
 =\rho_{\mathcal O_\mu}(t,r)
 \quad(1\le t\le3).
\]

For $t\geq4$, Lemma~\ref{lemmaforxingshu1} shows that only trace zero contributes and
identifies the restricted count with
$R_{\phi_t(f)}(4^{\delta_{2\nmid t}})$, where $f=f_{S^0}$.
Corollary~\ref{cor:tqf} and Lemma~\ref{lemmaforxingshu1} give
\[
 |\operatorname{Aut}(\phi_t(f))|
 =|\operatorname{Aut}(f)|
 =2|\operatorname{Aut}(\mathcal O)|.
\]
Thus the factor $1/2$ in~\eqref{eq:type-normalizer-count} cancels when the order counts are
converted into weighted representation numbers.  Applying Theorem~\ref{the:repnum}
to $G_{4N,16N^2,N_1}$ and to its transformed genera gives
\begin{equation}\label{eq:theorem13-first-formula}
\begin{aligned}
 T_{N_1,N_2}=2^{-e(N)-1}\Biggl(&
  \sum_{\substack{t\parallel N\\1\leq t\leq3}}
  \sum_{\substack{t\mid r\\r^2\leq4t}}
  H^{(N_1,N_2)}(4t-r^2)\\
 &+\sum_{\substack{t\parallel N\\t\geq4}}
  H^{(N_1,N_2;\,16^{{\delta_{2\mid t}}}t)}(4^{\delta_{2\nmid t}})
 \Biggr),
\end{aligned}
\end{equation}
which proves the first formula.

For the product formula, set
\[
 F(t)=\prod_{p\mid t}
 \frac{1-\left(\frac{\Delta(-4t)}p\right)/p}
      {B_p(t;N_1,N_2)C_p(t)}.
\]
Then Lemma~\ref{lemma5.5} gives, for $t\geq4$,
\[
 H^{(N_1,N_2;\,16^{{\delta_{2\mid t}}}t)}(4^{{\delta_{2\nmid t}}})
 =H^{(N_1,N_2)}(4t)F(t).
\]
Moreover, $F(1)=1$ by the empty-product convention, and $F(2)=F(3)=1$:
for the unique prime $p\mid t$ in these two cases, $v_p(t)=1$,
$\left(\frac{\Delta(-4t)}p\right)=0$, and
$B_p(t;N_1,N_2)=C_p(t)=1$.

If $t>4$, the conditions $t\mid r$ and $r^2\leq4t$ force $r=0$.
For $t=4$, the additional possibilities $r=\pm4$ correspond to the
scalar elements $\pm2$, which are not primitive and must be omitted.
Substituting the preceding identities into
\eqref{eq:theorem13-first-formula} therefore yields
\[
 T_{N_1,N_2}=2^{-e(N)-1}
   \sum_{t\parallel N}
   \sum_{\substack{t\mid r,\ r^2\leq4t\\
                    (t,r)\neq(4,\pm4)}}
   H^{(N_1,N_2)}(4t-r^2)F(t),
\]
which is the asserted product formula.
\end{proof}

\section{Examples}\label{Sec:examples}
In this section, we give three examples to illustrate the applications of our main results, corresponding respectively to the cases of $2\nmid N$, $2\parallel N$, and $4\mid N$.
	
	\subsection{Example for $2\nmid N$}
	Take $N=45=3^2\cdot5$. 
	By Theorem~\ref{the:genera}, the forms of level $180=4\cdot 45$ and discriminant $32400=16\cdot 45^2$ fall into $2^{2-1}=2$ genera, $G_{180,32400,5}$ and $G_{180,32400,45}$.
	By Theorem~\ref{the:classnum},
	\begin{equation*}
		|G_{180,32400,5}|=2,\, |G_{180,32400,45}|=1.
	\end{equation*}
The possible values of $\delta $ are $1,5,9,16,45,80,144,720$. By Theorem~\ref{the:repnum}, the representation numbers of these genera are given explicitly in Table~\ref{tab:example1}.

	\begin{longtable}{@{}ll@{}}
		\caption{Example for $N=45$}\label{tab:example1} \\
		\toprule[1.5pt]
		\multicolumn{1}{c}{\textbf{Genus}} & 
		\multicolumn{1}{c}{\textbf{Representation numbers}} \\
		\midrule[1pt]
		\endfirsthead 
		
		\multicolumn{2}{@{}l@{}}{\small\textbf{Table \thetable\ (cont.):} Example for $N=45$} \\
		\toprule[1.5pt]
		\multicolumn{1}{c}{\textbf{Genus}} & 
		\multicolumn{1}{c}{\textbf{Representation numbers}} \\
		\midrule[1pt]
		\endhead
		
		\midrule[0.5pt]
		\multicolumn{2}{r@{}}{\footnotesize\textit{continued}} \\
		\endfoot
		\bottomrule[1.5pt]
		\endlastfoot
		$G_{180,32400,5}$&$2R_{(8, 23, 47, 2, 8, 4)}(n)+2R_{(20, 23, 27, -18, 0, -20)}(n)=H^{(5,9)}(n)$\\
		$G_{180,32400,45}$&$2R_{(7, 27, 52, -24, -4, -6)}(n)=H^{(45,1)}(n)$\\
		$G_{180,6480,5}$&$2R_{(4, 19, 27, -18, 0, -4)}(n)+2R_{(7, 16, 16, -4, -4, -4)}(n)=H^{(5,9)}(5n)$\\
		$G_{180,6480,45}$&$2R_{(8, 11, 23, 2, 4, 8)}(n)=H^{(45,1)}(5n)$\\
		$G_{180,3600,5}$&$2R_{(3, 7, 47, -6, -2, -2)}(n)+2R_{(7, 12, 12, 4, 4, 4)}(n)=H^{(5,9;9)}(n)$\\
		$G_{180,3600,45}$&$2R_{(7, 7, 23, -2, -2, -6)}(n)=H^{(45,1;9)}(n)$\\
		$G_{180,2025,5}$&$2R_{(2, 17, 17, -11, -1, -1)}(n)+2R_{(5, 8, 17, 7, 5, 5)}(n)=H^{(5,9)}(4n)$\\
		$G_{180,2025,45}$&$2R_{(7, 7, 13, 2, 7, 4)}(n)=H^{(45,1)}(4n)$\\
		$G_{180,720,5}$&$2R_{(3, 3, 23, 2, 2, 2)}(n)+2R_{(3, 7, 11, -6, -2, -2)}(n)=H^{(5,9;9)}(5n)$\\
		$G_{180,720,45}$&$2R_{(4, 7, 8, -4, -4, 0)}(n)=H^{(45,1;9)}(5n)$\\
		$G_{180,405,5}$&$2R_{(1, 7, 16, 5, 1, 1)}(n)+2R_{(4, 4, 7, -2, -2, -1)}(n)=H^{(5,9)}(20n)$\\
		$G_{180,405,45}$&$2R_{(2, 8, 8, 7, 1, 1)}(n)=H^{(45,1)}(20n)$\\
		$G_{180,225,5}$&$2R_{(2, 3, 12, 3, 1, 2)}(n)+2R_{(3, 3, 7, 2, 2, 1)}(n)=H^{(5,9;9)}(4n)$\\
		$G_{180,225,45}$&$2R_{(2, 5, 7, -5, -1, 0)}(n)=H^{(45,1;9)}(4n)$\\
		$G_{180,45,5}$&$2R_{(1, 2, 6, -1, -1, 0)}(n)+2R_{(2, 3, 3, 3, 1, 2)}(n)=H^{(5,9;9)}(20n)$\\
		$G_{180,45,45}$&$2R_{(1, 2, 7, 2, 1, 1)}(n)=H^{(45,1;9)}(20n)$
	\end{longtable}
	
	\subsection{Example for $2\parallel N$}
	Take $N=50=2\cdot5^2$. 
	By Theorem~\ref{the:genera},  the forms of level $200=4\cdot 50$ and discriminant $40000=16\cdot 50^2$ fall into $2^{2-1}=2$ genera, $G_{200,40000,2}$ and $G_{200,40000,50}$.
	By Theorem~\ref{the:classnum}, 
	\begin{equation*}
		|G_{200,40000,2}|=2,\quad |G_{200,40000,50}|=2.
	\end{equation*}
	 The possible values of $\delta $ are $1,25,32,800$. Since $2\parallel N$, the forms of level $100=2\cdot 50$ and discriminant $2500=50^2$ likewise fall into $2^{2-1}=2$ genera, $G_{100,2500,2}$ and $G_{100,2500,50}$, with
	\begin{equation*}
		|G_{100,2500,2}|=2,\quad |G_{100,2500,50}|=2,
	\end{equation*}
	and here $\delta '=1,25$. By Theorem~\ref{the:repnum}, the representation numbers are given explicitly in Table~\ref{tab:example2}.

	\begin{longtable}{@{}ll@{}}
		\caption{Example for $N=50$}\label{tab:example2}\\
		\toprule[1.5pt]
		\multicolumn{1}{c}{\textbf{Genus}} & 
		\multicolumn{1}{c}{\textbf{Representation numbers}} \\
		\midrule[1pt]
		\endfirsthead 
		
		\multicolumn{2}{@{}l@{}}{\small\textbf{Table \thetable\ (cont.):} Example for $N=50$} \\
		\toprule[1.5pt]
		\multicolumn{1}{c}{\textbf{Genus}} & 
		\multicolumn{1}{c}{\textbf{Representation numbers}} \\
		\midrule[1pt]
		\endhead
		
		\midrule[0.5pt]
		\multicolumn{2}{r@{}}{\footnotesize\textit{continued}} \\
		\endfoot
		\bottomrule[1.5pt]
		\endlastfoot
		$G_{200,40000,2}$&$R_{(4, 51, 51, 2, 4, 4)}(n)+4R_{(11, 19, 56, -16, -8, -6)}(n)=2H^{(2,25)}(n)$\\
		$G_{200,40000,50}$&$2R_{(3, 67, 67, -66, -2, -2)}(n)+3R_{(8, 27, 52, 4, 8, 8)}(n)=3H^{(50,1)}(n)$\\
		$G_{200,1600,2}$&$R_{(4, 4, 27, -4, -4, 0)}(n)+4R_{(3, 8, 19, -8, -2, 0)}(n)=2H^{(2,25;25)}(n)$\\
		$G_{200,1600,50}$&$2R_{(8, 8, 11, 8, 8, 8)}(n)+3R_{(3, 3, 51, -2, -2, -2)}(n)=3H^{(50,1;25)}(n)$\\
		$G_{200,1250,2}$&$R_{(2, 13, 13, 1, 2, 2)}(n)+4R_{(3, 7, 17, 1, 2, 3)}(n)=2H^{(2,25)}(8n)$\\
		$G_{200,1250,50}$&$2R_{(6, 9, 9, -7, -4, -4)}(n)+3R_{(1, 19, 19, 13, 1, 1)}(n)=3H^{(50,1)}(8n)$\\
		$G_{200,50,2}$&$R_{(1, 1, 13, -1, -1, 0)}(n)+4R_{(1, 3, 5, 2, 1, 1)}(n)=2H^{(2,25;25)}(8n)$\\
		$G_{200,50,50}$&$2R_{(1, 1, 17, 1, 1, 1)}(n)+3R_{(1, 2, 7, -2, -1, 0)}(n)=3H^{(50,1;25)}(8n)$\\
		$G_{100,2500,2}$&$R_{(1, 25, 25, 0, 0, 0)}(n)+4R_{(6, 9, 14, 2, 6, 4)}(n)=2H^{(2,25)}(4n)$\\
		$G_{100,2500,50}$&$2R_{(3, 17, 17, -16, -2, -2)}(n)+3R_{(2, 13, 25, 0, 0, -2)}(n)=3H^{(50,1)}(4n)$\\
		$G_{100,100,2}$&$R_{(1, 1, 25, 0, 0, 0)}(n)+4R_{(2, 3, 5, -2, -2, 0)}(n)=2H^{(2,25;25)}(4n)$\\
		$G_{100,100,50}$&$2R_{(2, 2, 9, 2, 2, 2)}(n)+3R_{(1, 2, 13, -2, 0, 0)}(n)=3H^{(50,1;25)}(4n)$
	\end{longtable}

	\subsection{Example for $4\mid N$}
	Take $N=60=2^{2}\cdot3\cdot5$. 
	By Theorem~\ref{the:genera},  the forms of level $240=4\cdot 60$ and discriminant $57600=16\cdot 60^2$ fall into $2^{3-1}=4$ genera, $G_{240,57600,3}$, $G_{240,57600,5}$, $G_{240,57600,12}$, and $G_{240,57600,20}$.
	By Theorem~\ref{the:classnum}, 
	\begin{equation*}
		|G_{240,57600,3}|=2,\, |G_{240,57600,5}|=3,\, |G_{240,57600,12}|=1,\, |G_{240,57600,20}|=2.
	\end{equation*}
	The possible values of $\delta $ are $1,3,5,15,64,192,320,960$. By Theorem~\ref{the:repnum}, the representation numbers are given explicitly in Table~\ref{tab:example3}.

	\begin{longtable}{@{}ll@{}}
		\caption{Example for $N=60$}\label{tab:example3}\\
		\toprule[1.5pt]
		\multicolumn{1}{c}{\textbf{Genus}} & 
		\multicolumn{1}{c}{\textbf{Representation numbers}} \\
		\midrule[1pt]
		\endfirsthead 
		
		\multicolumn{2}{@{}l@{}}{\small\textbf{Table \thetable\ (cont.):} Example for $N=60$} \\
		\toprule[1.5pt]
		\multicolumn{1}{c}{\textbf{Genus}} & 
		\multicolumn{1}{c}{\textbf{Representation numbers}} \\
		\midrule[1pt]
		\endhead
		
		\midrule[0.5pt]
		\multicolumn{2}{r@{}}{\footnotesize\textit{continued}} \\
		\endfoot
		\bottomrule[1.5pt]
		\endlastfoot
		$G_{240,57600,3}$&$2R_{(15, 16, 64, -16, 0, 0)}(n)+4R_{(16, 31, 31, 2, 8, 8)}(n)=H^{(3,20)}(n)$\\
		$G_{240,57600,5}$&$2R_{(12, 23, 63, 6, 12, 12)}(n)+2R_{(15, 32, 32, -16, 0, 0)}(n)+4R_{(23, 23, 32, -8, -8, -14)}(n)=H^{(5,12)}(n)$\\
		$G_{240,57600,12}$&$2R_{(16, 19, 64, 8, 16, 16)}(n)=H^{(12,5)}(n)$\\
		$G_{240,57600,20}$&$2R_{(3, 80, 80, -80, 0, 0)}(n)+6R_{(27, 27, 32, -24, -24, -6)}(n)=3H^{(20,3)}(n)$\\
		$G_{240,19200,3}$&$2R_{(5, 32, 32, -16, 0, 0)}(n)+4R_{(13, 13, 32, 8, 8, 6)}(n)=H^{(3,20)}(3n)$\\
		$G_{240,19200,5}$&$2R_{(4, 21, 61, 2, 4, 4)}(n)+2R_{(5, 16, 64, -16, 0, 0)}(n)+4R_{(16, 20, 21, -20, -8, 0)}(n)=H^{(5,12)}(3n)$\\
		$G_{240,19200,12}$&$2R_{(17, 17, 25, -10, -10, -14)}(n)=H^{(12,5)}(3n)$\\
		$G_{240,19200,20}$&$2R_{(1, 80, 80, -80, 0, 0)}(n)+6R_{(9, 9, 64, -8, -8, -2)}(n)=3H^{(20,3)}(3n)$\\
		$G_{240,11520,3}$&$2R_{(3, 32, 32, -16, 0, 0)}(n)+4R_{(11, 11, 32, 8, 8, 10)}(n)=H^{(3,20)}(5n)$\\
		$G_{240,11520,5}$&$2R_{(3, 16, 64, -16, 0, 0)}(n)+2R_{(16, 16, 19, 8, 16, 16)}(n)+4R_{(12, 16, 19, -8, -12, 0)}(n)=H^{(5,12)}(5n)$\\
		$G_{240,11520,12}$&$2R_{(12, 15, 23, 6, 12, 12)}(n)=H^{(12,5)}(5n)$\\
		$G_{240,11520,20}$&$2R_{(15, 16, 16, -16, 0, 0)}(n)+6R_{(7, 7, 63, -6, -6, -2)}(n)=3H^{(20,3)}(5n)$\\
		$G_{240,3840,3}$&$2R_{(1, 16, 64, -16, 0, 0)}(n)+4R_{(4, 16, 17, -8, -4, 0)}(n)=H^{(3,20)}(15n)$\\
		$G_{240,3840,5}$&$2R_{(1, 32, 32, -16, 0, 0)}(n)+2R_{(9, 9, 16, -8, -8, -2)}(n)+4R_{(4, 9, 32, -8, 0, -4)}(n)=H^{(5,12)}(15n)$\\
		$G_{240,3840,12}$&$2R_{(4, 5, 61, 2, 4, 4)}(n)=H^{(12,5)}(15n)$\\
		$G_{240,3840,20}$&$2R_{(5, 16, 16, -16, 0, 0)}(n)+6R_{(4, 13, 21, 2, 4, 4)}(n)=3H^{(20,3)}(15n)$\\
		$G_{240,900,3}$&$2R_{(1, 4, 60, 0, 0, -1)}(n)+4R_{(4, 4, 16, 1, 4, 2)}(n)=H^{(3,20;64)}(n)$\\
		$G_{240,900,5}$&$2R_{(2, 2, 60, 0, 0, -1)}(n)+2R_{(5, 8, 8, 4, 5, 5)}(n)+4R_{(2, 8, 17, 8, 1, 2)}(n)=H^{(5,12;64)}(n)$\\
		$G_{240,900,12}$&$2R_{(1, 15, 19, -15, -1, 0)}(n)=H^{(12,5;64)}(n)$\\
		$G_{240,900,20}$&$2R_{(5, 5, 12, 0, 0, -5)}(n)+6R_{(3, 5, 17, -5, -3, 0)}(n)=3H^{(20,3;64)}(n)$\\
		$G_{240,300,3}$&$2R_{(2, 2, 20, 0, 0, -1)}(n)+4R_{(2, 7, 7, -6, -1, -1)}(n)=H^{(3,20;64)}(3n)$\\
		$G_{240,300,12}$&$2R_{(3, 5, 7, -5, -3, 0)}(n)=H^{(12,5;64)}(3n)$\\
		$G_{240,300,5}$&$2R_{(1, 4, 20, 0, 0, -1)}(n)+2R_{(4, 5, 6, -5, -4, 0)}(n)+4R_{(4, 4, 6, 1, 4, 2)}(n)=H^{(5,12;64)}(3n)$\\
		$G_{240,300,20}$&$2R_{(4, 5, 5, -5, 0, 0)}(n)+6R_{(1, 9, 9, 3, 1, 1)}(n)=3H^{(20,3;64)}(3n)$\\
		$G_{240,180,3}$&$2R_{(2, 2, 12, 0, 0, -1)}(n)+4R_{(2, 5, 5, -2, -1, -1)}(n)=H^{(3,20;64)}(5n)$\\
		$G_{240,180,5}$&$2R_{(1, 3, 16, -3, -1, 0)}(n)+2R_{(1, 4, 12, 0, 0, -1)}(n)+4R_{(4, 4, 4, 1, 2, 4)}(n)=H^{(5,12;64)}(5n)$\\
		$G_{240,180,12}$&$2R_{(3, 5, 5, 5, 3, 3)}(n)=H^{(12,5;64)}(5n)$\\
		$G_{240,180,20}$&$2R_{(1, 1, 60, 0, 0, -1)}(n)+6R_{(1, 7, 7, 2, 1, 1)}(n)=3H^{(20,3;64)}(5n)$\\
		$G_{240,60,3}$&$2R_{(1, 4, 4, 0, 0, -1)}(n)+4R_{(2, 2, 4, -1, -1, 0)}(n)=H^{(3,20;64)}(15n)$\\
		$G_{240,60,5}$&$2R_{(1, 3, 6, -3, -1, 0)}(n)+2R_{(2, 2, 4, 0, 0, -1)}(n)+4R_{(2, 3, 3, 2, 1, 1)}(n)=H^{(5,12;64)}(15n)$\\
		$G_{240,60,12}$&$2R_{(1, 4, 5, -4, -1, 0)}(n)=H^{(12,5;64)}(15n)$\\
		$G_{240,60,20}$&$2R_{(1, 1, 20, 0, 0, -1)}(n)+6R_{(3, 3, 3, 2, 3, 3)}(n)=3H^{(20,3;64)}(15n)$\\
	\end{longtable}
	
	\newpage 
	\appendix
	
	\section{Table of one-class genera}\label{App:A}
	When $|G_{4N_1N_2,\,16(N_1N_2)^2/\delta ,\,N_1}|=1$, we can obtain  explicit formulas for the representation numbers of the corresponding  ternary quadratic forms. By Theorem~\ref{the:classnum},  the $39$ pairs $(N_1,N_2)$ with class number $1$ are listed in Table~\ref{tab:class1}.
	
	\begin{table}[H]
		\caption{$(N_1,N_2)$ with class number $1$}\label{tab:class1}
		\centering
		\renewcommand{\arraystretch}{1.3} %
		\begin{tabular}{|c|c|c|}
			\hline
			\multirow{2}{*}{\(e(N_1N_2)\)} & \multicolumn{2}{c|}{$(N_1,N_2)$} \\
			\cline{2-3}
			& \(2 \nmid N_1N_2\) & \(2 \mid N_1N_2\) \\
			\hline
			1 & \((3,1),(5,1),(7,1),(13,1)\) & \((2,1),(8,1)\) \\
			\hline
			\multirow{3}{*}{2} & \multirow{3}{*}{\((3,5),(3,11),(7,3),(45,1)\)} & \((2,3),(2,5),(2,7),(2,9),(2,11),(2,23)\) \\[4pt]
			& & \((3,2),(3,4),(3,8),(5,2),(5,4)\) \\[4pt]
			& & \((8,5),(12,1),(18,1),(20,1),(28,1),(52,1),(72,1)\) \\
			\hline
			\multirow{2}{*}{3} & & \multirow{2}{*}{\((2,15),(12,5),(12,11),(18,7),(28,3)\)} \\[4pt]
			& & \((30,1),(42,1),(70,1),(78,1),(168,1),(180,1)\) \\
			\hline
		\end{tabular}
	\end{table}
	By Theorem~\ref{the:repnum} we obtain explicit formulas for the representation numbers of $255$ ternary quadratic forms; 
	of these, $110$ are reported here for the first time, while the remaining $145$ are already known from earlier work in which $N_1N_2$ is squarefree. 
	The forms listed in Table~\ref{tab:genera-detail} were computed with SageMath and subsequently verified by hand.

	\begin{longtable}{@{}lll@{}}
		\caption{One-class Genera}\label{tab:genera-detail}\\
		\toprule[1.5pt]
		\multicolumn{1}{c}{\textbf{$N_f$}} & 
		\multicolumn{1}{c}{\textbf{$d_f$}} & 
		\multicolumn{1}{c}{\textbf{$R_f(n)$}} \\
		\midrule[1pt]
		\endfirsthead 
		
		\multicolumn{3}{@{}l@{}}{\small\textbf{Table \thetable\ (cont.):} One-class Genera} \\
		\toprule[1.5pt]
		\multicolumn{1}{c}{\textbf{$N_f$}} & 
		\multicolumn{1}{c}{\textbf{$d_f$}} & 
		\multicolumn{1}{c}{\textbf{$R_f(n)$}} \\
		\midrule[1pt]
		\endhead
		
		\midrule[0.5pt]
		\multicolumn{3}{r@{}}{\footnotesize\textit{to be continued}} \\
		\endfoot
		\bottomrule[1.5pt]
		\endlastfoot
		$4=2^2$&$4=2^2$&$R_{(1,1,1,0,0,0)}(n)=12H^{(2,1)}(4n)$\\
		$8=2^3$&$2=2^1$&$R_{(1,1,1,1,1,1)}(n)=12H^{(2,1)}(8n)$\\
		$8=2^3$&$64=2^6$&$R_{(3,3,3,-2,-2,-2)}(n)=12H^{(2,1)}(n)$\\
		$12=2^2\cdot3$&$3=3^1$&$R_{(1,1,1,0,0,-1)}(n)=6H^{(3,1)}(12n)$\\
		$12=2^2\cdot3$&$9=3^2$&$R_{(1,1,3,0,0,-1)}(n)=6H^{(3,1)}(4n)$\\
		$12=2^2\cdot3$&$12=2^2\cdot3$&$R_{(1,2,2,-2,0,0)}(n)=3H^{(2,3)}(12n)$\\
		$12=2^2\cdot3$&$12=2^2\cdot3$&$R_{(1,1,3,0,0,0)}(n)=2H^{(3,2)}(12n)$\\
		$12=2^2\cdot3$&$36=2^2\cdot3^2$&$R_{(2,2,3,0,0,-2)}(n)=3H^{(2,3)}(4n)$\\
		$12=2^2\cdot3$&$36=2^2\cdot3^2$&$R_{(1,3,3,0,0,0)}(n)=2H^{(3,2)}(4n)$\\
		$12=2^2\cdot3$&$48=2^4\cdot3$&$R_{(1,4,4,-4,0,0)}(n)=6H^{(3,1)}(3n)$\\
		$12=2^2\cdot3$&$144=2^4\cdot3^2$&$R_{(3,4,4,-4,0,0)}(n)=6H^{(3,1)}(n)$\\
		$20=2^2\cdot5$&$5=5^1$&$R_{(1,1,2,1,1,1)}(n)=3H^{(5,1)}(20n)$\\
		$20=2^2\cdot5$&$20=2^2\cdot5$&$R_{(1,1,5,0,0,0)}(n)=2H^{(2,5)}(20n)$\\
		$20=2^2\cdot5$&$20=2^2\cdot5$&$R_{(1,2,3,-2,0,0)}(n)=H^{(5,2)}(20n)$\\
		$20=2^2\cdot5$&$25=5^2$&$R_{(2,2,2,-1,-1,-1)}(n)=3H^{(5,1)}(4n)$\\
		$20=2^2\cdot5$&$80=2^4\cdot5$&$R_{(3,3,3,2,2,2)}(n)=3H^{(5,1)}(5n)$\\
		$20=2^2\cdot5$&$100=2^2\cdot5^2$&$R_{(1,5,5,0,0,0)}(n)=2H^{(2,5)}(4n)$\\
		$20=2^2\cdot5$&$100=2^2\cdot5^2$&$R_{(2,3,5,0,0,-2)}(n)=H^{(5,2)}(4n)$\\
		$20=2^2\cdot5$&$400=2^4\cdot5^2$&$R_{(3,7,7,-6,-2,-2)}(n)=3H^{(5,1)}(n)$\\
		$24=2^3\cdot3$&$6=2\cdot3$&$R_{(1, 1, 2, 0, 0, -1)}(n)=3H^{(2,3)}(24n)$\\
		$24=2^3\cdot3$&$6=2\cdot3$&$R_{(1, 1, 2, -1, -1, 0)}(n)=2H^{(3,2)}(24n)$\\
		$24=2^3\cdot3$&$18=2\cdot3^2$&$R_{(1, 1, 6, 0, 0, -1)}(n)=3H^{(2,3)}(8n)$\\
		$24=2^3\cdot3$&$18=2\cdot3^2$&$R_{(2, 2, 2, 1, 2, 2)}(n)=2H^{(3,2)}(8n)$\\
		$24=2^3\cdot3$&$192=2^6\cdot3$&$R_{(1, 8, 8, -8, 0, 0)}(n)=3H^{(2,3)}(3n)$\\
		$24=2^3\cdot3$&$192=2^6\cdot3$&$R_{(4, 4, 5, -4, -4, 0)}(n)=2H^{(3,2)}(3n)$\\
		$24=2^3\cdot3$&$576=2^6\cdot3^2$&$R_{(3, 8, 8, -8, 0, 0)}(n)=3H^{(2,3)}(n)$\\
		$24=2^3\cdot3$&$576=2^6\cdot3^2$&$R_{(4, 7, 7, 2, 4, 4)}(n)=2H^{(3,2)}(n)$\\
		$28=2^2\cdot7$&$7=7^1$&$R_{(1,1,2,0,-1,0)}(n)=2H^{(7,1)}(28n)$\\
		$28=2^2\cdot7$&$28=2^2\cdot7$&$2R_{(2,2,3,2,2,2)}(n)=3H^{(2,7)}(28n)$\\
		$28=2^2\cdot7$&$49=7^2$&$R_{(1,2,7,0,0,-1)}(n)=2H^{(7,1)}(4n)$\\
		$28=2^2\cdot7$&$112=2^4\cdot7$&$R_{(1,4,8,-4,0,0)}(n)=2H^{(7,1)}(7n)$\\
		$28=2^2\cdot7$&$196=2^2\cdot7^2$&$2R_{(3,5,5,-4,-2,-2)}(n)=3H^{(2,7)}(4n)$\\
		$28=2^2\cdot7$&$784=2^4\cdot7^2$&$R_{(4,7,8,0,-4,0)}(n)=2H^{(7,1)}(n)$\\
		$32=2^5$&$8=2^{3}$&$R_{(1, 1, 3, 1, 1, 1)}(n)=3H^{(8,1;128)}(n)$\\
		$32=2^5$&$1024=2^{10}$&$R_{(3, 11, 11, -10, -2, -2)}(n)=3H^{(8,1)}(n)$\\
		$36=2^2\cdot3^2$&$36=2^2\cdot3^2$&$R_{(1, 2, 5, -2, 0, 0)}(n)=H^{(2,9;9)}(4n)$\\
		$36=2^2\cdot3^2$&$36=2^2\cdot3^2$&$R_{(1, 1, 9, 0, 0, 0)}(n)=2H^{(18,1;9)}(4n)$\\
		$36=2^2\cdot3^2$&$324=2^2\cdot3^4$&$R_{(2, 5, 9, 0, 0, -2)}(n)=H^{(2,9)}(4n)$\\
		$36=2^2\cdot3^2$&$324=2^2\cdot3^4$&$R_{(1, 9, 9, 0, 0, 0)}(n)=2H^{(18,1)}(4n)$\\
		$40=2^3\cdot5$&$10=2\cdot5$&$R_{(1, 1, 3, -1, -1, 0)}(n)=2H^{(2,5)}(40n)$\\
		$40=2^3\cdot5$&$10=2\cdot5$&$R_{(1, 2, 2, 2, 1, 1)}(n)=H^{(5,2)}(40n)$\\
		$40=2^3\cdot5$&$50=2\cdot5^2$&$R_{(2, 3, 3, 1, 2, 2)}(n)=2H^{(2,5)}(8n)$\\
		$40=2^3\cdot5$&$50=2\cdot5^2$&$R_{(1, 4, 4, 3, 1, 1)}(n)=H^{(5,2)}(8n)$\\
		$40=2^3\cdot5$&$320=2^6\cdot5$&$R_{(4, 4, 7, -4, -4, 0)}(n)=2H^{(2,5)}(5n)$\\
		$40=2^3\cdot5$&$320=2^6\cdot5$&$R_{(3, 3, 11, -2, -2, -2)}(n)=H^{(5,2)}(5n)$\\
		$40=2^3\cdot5$&$1600=2^6\cdot5^2$&$R_{(4, 11, 11, 2, 4, 4)}(n)=2H^{(2,5)}(n)$\\
		$40=2^3\cdot5$&$1600=2^6\cdot5^2$&$R_{(7, 7, 12, -4, -4, -6)}(n)=H^{(5,2)}(n)$\\
		$44=2^2\cdot11$&$44=2^2\cdot11$&$R_{(1,2,6,-2,0,0)}(n)=H^{(2,11)}(44n)$\\
		$44=2^2\cdot11$&$484=2^2\cdot11^2$&$R_{(2,6,11,0,0,-2)}(n)=H^{(2,11)}(4n)$\\
		$48=2^4\cdot3$&$12=2^2\cdot3$&$R_{(1, 2, 2, 1, 1, 1)}(n)=H^{(3,4;64)}(3n)$\\
		$48=2^4\cdot3$&$12=2^2\cdot3$&$R_{(1, 1, 4, 0, 0, -1)}(n)=3H^{(12,1;64)}(3n)$\\
		$48=2^4\cdot3$&$36=2^2\cdot3^2$&$R_{(1, 1, 12, 0, 0, -1)}(n)=3H^{(12,1;64)}(n)$\\
		$48=2^4\cdot3$&$36=2^2\cdot3^2$&$R_{(1, 3, 4, -3, -1, 0)}(n)=H^{(3,4;64)}(n)$\\
		$48=2^4\cdot3$&$768=2^8\cdot3$&$R_{(4, 5, 13, 2, 4, 4)}(n)=H^{(3,4)}(3n)$\\
		$48=2^4\cdot3$&$768=2^8\cdot3$&$R_{(1, 16, 16, -16, 0, 0)}(n)=3H^{(12,1)}(3n)$\\
		$48=2^4\cdot3$&$2304=2^8\cdot3^2$&$R_{(7, 7, 15, -6, -6, -2)}(n)=H^{(3,4)}(n)$\\
		$48=2^4\cdot3$&$2304=2^8\cdot3^2$&$R_{(3, 16, 16, -16, 0, 0)}(n)=3H^{(12,1)}(n)$\\
		$52=2^2\cdot13$&$13=13^1$&$R_{(1,2,2,-1,0,-1)}(n)=H^{(13,1)}(52n)$\\
		$52=2^2\cdot13$&$169=13^2$&$R_{(2,5,5,-3,-1,-1)}(n)=H^{(13,1)}(4n)$\\
		$52=2^2\cdot13$&$208=2^4\cdot13$&$R_{(3,3,7,2,2,2)}(n)=H^{(13,1)}(13n)$\\
		$52=2^2\cdot13$&$2704=2^4\cdot13^2$&$R_{(7,8,15,8,2,4)}(n)=H^{(13,1)}(n)$\\
		$56=2^3\cdot7$&$14=2\cdot7$&$2R_{(1, 1, 5, 1, 1, 1)}(n)=3H^{(2,7)}(56n)$\\
		$56=2^3\cdot7$&$98=2\cdot7^2$&$2R_{(3, 3, 3, -1, -1, -1)}(n)=3H^{(2,7)}(8n)$\\
		$56=2^3\cdot7$&$448=2^6\cdot7$&$2R_{(5, 5, 5, 2, 2, 2)}(n)=3H^{(2,7)}(7n)$\\
		$56=2^3\cdot7$&$3136=2^6\cdot7^2$&$2R_{(3, 19, 19, -18, -2, -2)}(n)=3H^{(2,7)}(n)$\\
		$60=2^2\cdot3\cdot5$&$15=3\cdot5$&$R_{(1,1,4,0,-1,0)}(n)=H^{(3,5)}(60n)$\\
		$60=2^2\cdot3\cdot5$&$45=3^2\cdot5$&$R_{(2,2,3,0,0,-1)}(n)=H^{(3,5)}(20n)$\\
		$60=2^2\cdot3\cdot5$&$60=2^2\cdot3\cdot5$&$2R_{(2,3,3,0,0,-2)}(n)=H^{(2,15)}(60n)$\\
		$60=2^2\cdot3\cdot5$&$60=2^2\cdot3\cdot5$&$2R_{(2,2,5,0,0,-2)}(n)=3H^{(30,1)}(60n)$\\
		$60=2^2\cdot3\cdot5$&$75=3\cdot5^2$&$R_{(2,2,5,0,0,-1)}(n)=H^{(3,5)}(12n)$\\
		$60=2^2\cdot3\cdot5$&$180=2^2\cdot3^2\cdot5$&$2R_{(1,6,9,-6,0,0)}(n)=H^{(2,15)}(20n)$\\
		$60=2^2\cdot3\cdot5$&$180=2^2\cdot3^2\cdot5$&$2R_{(2,2,15,0,0,-2)}(n)=3H^{(30,1)}(20n)$\\
		$60=2^2\cdot3\cdot5$&$225=3^2\cdot5^2$&$R_{(1,4,15,0,0,-1)}(n)=H^{(3,5)}(4n)$\\
		$60=2^2\cdot3\cdot5$&$240=2^4\cdot3\cdot5$&$R_{(1,4,16,-4,0,0)}(n)=H^{(3,5)}(15n)$\\
		$60=2^2\cdot3\cdot5$&$300=2^2\cdot3\cdot5^2$&$2R_{(2,3,15,0,0,-2)}(n)=H^{(2,15)}(12n)$\\
		$60=2^2\cdot3\cdot5$&$300=2^2\cdot3\cdot5^2$&$2R_{(1,10,10,-10,0,0)}(n)=3H^{(30,1)}(12n)$\\
		$60=2^2\cdot3\cdot5$&$720=2^4\cdot3^2\cdot5$&$R_{(3,8,8,-4,0,0)}(n)=H^{(3,5)}(5n)$\\
		$60=2^2\cdot3\cdot5$&$900=2^2\cdot3^2\cdot5^2$&$2R_{(5,6,9,-6,0,0)}(n)=H^{(2,15)}(4n)$\\
		$60=2^2\cdot3\cdot5$&$900=2^2\cdot3^2\cdot5^2$&$2R_{(3,10,10,-10,0,0)}(n)=3H^{(30,1)}(4n)$\\
		$60=2^2\cdot3\cdot5$&$1200=2^4\cdot3\cdot5^2$&$R_{(5,8,8,-4,0,0)}(n)=H^{(3,5)}(3n)$\\
		$60=2^2\cdot3\cdot5$&$3600=2^4\cdot3^2\cdot5^2$&$R_{(4,15,16,0,-4,0)}(n)=H^{(3,5)}(n)$\\
		$72=2^3\cdot3^2$&$18=2\cdot3^2$&$R_{(1, 1, 5, -1, -1, 0)}(n)=2H^{(18,1;9)}(8n)$\\
		$72=2^3\cdot3^2$&$18=2\cdot3^2$&$R_{(1, 2, 3, -2, -1, 0)}(n)=H^{(2,9;9)}(8n)$\\
		$72=2^3\cdot3^2$&$162=2\cdot3^4$&$R_{(2, 5, 5, 1, 2, 2)}(n)=2H^{(18,1)}(8n)$\\
		$72=2^3\cdot3^2$&$162=2\cdot3^4$&$R_{(1, 7, 7, 5, 1, 1)}(n)=H^{(2,9)}(8n)$\\
		$72=2^3\cdot3^2$&$576=2^6\cdot3^2$&$R_{(3, 3, 19, -2, -2, -2)}(n)=H^{(2,9;9)}(n)$\\
		$72=2^3\cdot3^2$&$576=2^6\cdot3^2$&$R_{(4, 4, 11, -4, -4, 0)}(n)=2H^{(18,1;9)}(n)$\\
		$72=2^3\cdot3^2$&$5184=2^6\cdot3^4$&$R_{(8, 11, 20, 4, 8, 8)}(n)=H^{(2,9)}(n)$\\
		$72=2^3\cdot3^2$&$5184=2^6\cdot3^4$&$R_{(4, 19, 19, 2, 4, 4)}(n)=2H^{(18,1)}(n)$\\
		$80=2^4\cdot5$&$20=2^2\cdot5$&$2R_{(1, 1, 7, 1, 1, 1)}(n)=3H^{(20,1;64)}(5n)$\\
		$80=2^4\cdot5$&$20=2^2\cdot5$&$2R_{(1, 2, 3, -1, 0, -1)}(n)=H^{(5,4;64)}(5n)$\\
		$80=2^4\cdot5$&$100=2^2\cdot5^2$&$2R_{(2, 2, 7, -1, -1, -1)}(n)=H^{(5,4;64)}(n)$\\
		$80=2^4\cdot5$&$100=2^2\cdot5^2$&$2R_{(3, 3, 3, 1, 1, 1)}(n)=3H^{(20,1;64)}(n)$\\
		$80=2^4\cdot5$&$1280=2^8\cdot5$&$2R_{(3, 11, 11, 6, 2, 2)}(n)=H^{(5,4)}(5n)$\\
		$80=2^4\cdot5$&$1280=2^8\cdot5$&$2R_{(7, 7, 7, -2, -2, -2)}(n)=3H^{(20,1)}(5n)$\\
		$80=2^4\cdot5$&$6400=2^8\cdot5^2$&$2R_{(7, 12, 23, 12, 2, 4)}(n)=H^{(5,4)}(n)$\\
		$80=2^4\cdot5$&$6400=2^8\cdot5^2$&$2R_{(3, 27, 27, -26, -2, -2)}(n)=3H^{(20,1)}(n)$\\
		$84=2^2\cdot3\cdot7$&$21=3\cdot7$&$2R_{(1,2,3,-1,-1,0)}(n)=H^{(7,3)}(84n)$\\
		$84=2^2\cdot3\cdot7$&$63=3^2\cdot7$&$2R_{(2,2,5,2,2,1)}(n)=H^{(7,3)}(28n)$\\
		$84=2^2\cdot3\cdot7$&$84=2^2\cdot3\cdot7$&$R_{(1,1,21,0,0,0)}(n)=H^{(42,1)}(84n)$\\
		$84=2^2\cdot3\cdot7$&$147=3\cdot7^2$&$2R_{(3,3,5,-2,-2,-1)}(n)=H^{(7,3)}(12n)$\\
		$84=2^2\cdot3\cdot7$&$252=2^2\cdot3^2\cdot7$&$R_{(3,3,7,0,0,0)}(n)=H^{(42,1)}(28n)$\\
		$84=2^2\cdot3\cdot7$&$336=2^4\cdot3\cdot7$&$2R_{(3,3,11,2,2,2)}(n)=H^{(7,3)}(21n)$\\
		$84=2^2\cdot3\cdot7$&$441=3^2\cdot7^2$&$2R_{(2,8,8,-5,-1,-1)}(n)=H^{(7,3)}(4n)$\\
		$84=2^2\cdot3\cdot7$&$588=2^2\cdot3\cdot7^2$&$R_{(3,7,7,0,0,0)}(n)=H^{(42,1)}(12n)$\\
		$84=2^2\cdot3\cdot7$&$1008=2^4\cdot3^2\cdot7$&$2R_{(5,8,8,4,4,4)}(n)=H^{(7,3)}(7n)$\\
		$84=2^2\cdot3\cdot7$&$1764=2^2\cdot3^2\cdot7^2$&$R_{(1,21,21,0,0,0)}(n)=H^{(42,1)}(4n)$\\
		$84=2^2\cdot3\cdot7$&$2352=2^4\cdot3\cdot7^2$&$2R_{(5,12,12,-4,-4,-4)}(n)=H^{(7,3)}(3n)$\\
		$84=2^2\cdot3\cdot7$&$7056=2^4\cdot3^2\cdot7^2$&$2R_{(8,11,23,2,8,4)}(n)=H^{(7,3)}(n)$\\
		$88=2^3\cdot11$&$22=2\cdot11$&$R_{(1, 2, 3, 0, -1, 0)}(n)=H^{(2,11)}(88n)$\\
		$88=2^3\cdot11$&$242=2\cdot11^2$&$R_{(1, 3, 22, 0, 0, -1)}(n)=H^{(2,11)}(8n)$\\
		$88=2^3\cdot11$&$704=2^6\cdot11$&$R_{(1, 8, 24, -8, 0, 0)}(n)=H^{(2,11)}(11n)$\\
		$88=2^3\cdot11$&$7744=2^6\cdot11^2$&$R_{(8, 11, 24, 0, -8, 0)}(n)=H^{(2,11)}(n)$\\
		$92=2^2\cdot23$&$92=2^2\cdot23$&$2R_{(2,3,5,-2,0,-2)}(n)=H^{(2,23)}(92n)$\\
		$92=2^2\cdot23$&$2116=2^2\cdot23^2$&$2R_{(5,10,14,10,2,4)}(n)=H^{(2,23)}(4n)$\\
		$96=2^5\cdot3$&$24=2^{3}\cdot3$&$2R_{(1, 2, 4, 2, 1, 1)}(n)=H^{(3,8;128)}(3n)$\\
		$96=2^5\cdot3$&$72=2^{3}\cdot3^2$&$2R_{(2, 2, 5, 1, 1, 1)}(n)=H^{(3,8;128)}(n)$\\
		$96=2^5\cdot3$&$3072=2^{10}\cdot3$&$2R_{(5, 13, 13, -6, -2, -2)}(n)=H^{(3,8)}(3n)$\\
		$96=2^5\cdot3$&$9216=2^{10}\cdot3^2$&$2R_{(7, 15, 28, -12, -4, -6)}(n)=H^{(3,8)}(n)$\\
		$112=2^4\cdot7$&$28=2^2\cdot7$&$R_{(1,3, 3, 2, 1, 1)}(n)=H^{(28,1;64)}(7n)$\\
		$112=2^4\cdot7$&$196=2^2\cdot7^2$&$R_{(1, 7, 9, -7, -1, 0)}(n)=H^{(28,1;64)}(n)$\\
		$112=2^4\cdot7$&$1792=2^8\cdot7$&$R_{(4, 5, 29, 2, 4, 4)}(n)=H^{(28,1)}(7n)$\\
		$112=2^4\cdot7$&$12544=2^8\cdot7^2$&$R_{(11, 11, 32, -8, -8, -6)}(n)=H^{(28,1)}(n)$\\
		$120=2^3\cdot3\cdot5$&$30=2\cdot3\cdot5$&$2R_{(1,3,3,1,1,1)}(n)=H^{(2,15)}(120n)$\\
		$120=2^3\cdot3\cdot5$&$30=2\cdot3\cdot5$&$2R_{(1,1,10,0,0,-1)}(n)=3H^{(30,1)}(120n)$\\
		$120=2^3\cdot3\cdot5$&$90=2\cdot3^2\cdot5$&$2R_{(2,3,5,-3,2,0)}(n)=H^{(2,15)}(40n)$\\
		$120=2^3\cdot3\cdot5$&$90=2\cdot3^2\cdot5$&$2R_{(1,1,30,0,0,-1)}(n)=3H^{(30,1)}(40n)$\\
		$120=2^3\cdot3\cdot5$&$150=2\cdot3\cdot5^2$&$2R_{(1,5,9,-5,1,0)}(n)=H^{(2,15)}(24n)$\\
		$120=2^3\cdot3\cdot5$&$150=2\cdot3\cdot5^2$&$2R_{(2,5,5,-5,0,0)}(n)=3H^{(30,1)}(24n)$\\
		$120=2^3\cdot3\cdot5$&$450=2\cdot3^2\cdot5^2$&$2R_{(3,7,7,4,3,3)}(n)=H^{(2,15)}(8n)$\\
		$120=2^3\cdot3\cdot5$&$450=2\cdot3^2\cdot5^2$&$2R_{(5,5,6,0,0,-5)}(n)=3H^{(30,1)}(8n)$\\
		$120=2^3\cdot3\cdot5$&$960=2^6\cdot3\cdot5$&$2R_{(5,5,12,-4,-4,-2)}(n)=H^{(2,15)}(15n)$\\
		$120=2^3\cdot3\cdot5$&$960=2^6\cdot3\cdot5$&$2R_{(5,8,8,-8,0,0)}(n)=3H^{(30,1)}(15n)$\\
		$120=2^3\cdot3\cdot5$&$2880=2^6\cdot3^2\cdot5$&$2R_{(4,7,31,2,4,4)}(n)=H^{(2,15)}(5n)$\\
		$120=2^3\cdot3\cdot5$&$2880=2^6\cdot3^2\cdot5$&$2R_{(8,8,15,0,0,-8)}(n)=3H^{(30,1)}(5n)$\\
		$120=2^3\cdot3\cdot5$&$4800=2^6\cdot3\cdot5^2$&$2R_{(8,12,17,4,8,8)}(n)=H^{(2,15)}(3n)$\\
		$120=2^3\cdot3\cdot5$&$4800=2^6\cdot3\cdot5^2$&$2R_{(1,40,40,-40,0,0)}(n)=3H^{(30,1)}(3n)$\\
		$120=2^3\cdot3\cdot5$&$14400=2^6\cdot3^2\cdot5^2$&$2R_{(11,11,35,-10,-10,-2)}(n)=H^{(2,15)}(n)$\\
		$120=2^3\cdot3\cdot5$&$14400=2^6\cdot3^2\cdot5^2$&$2R_{(3,40,40,-40,0,0)}(n)=3H^{(30,1)}(n)$\\
		$132=2^2\cdot3\cdot11$&$33=3\cdot11$&$2R_{(1,2,5,1,1,1)}(n)=H^{(3,11)}(132n)$\\
		$132=2^2\cdot3\cdot11$&$99=3^2\cdot11$&$2R_{(2,3,5,-3,-1,0)}(n)=H^{(3,11)}(44n)$\\
		$132=2^2\cdot3\cdot11$&$363=3\cdot11^2$&$2R_{(2,7,7,3,1,1)}(n)=H^{(3,11)}(12n)$\\
		$132=2^2\cdot3\cdot11$&$528=2^4\cdot3\cdot11$&$2R_{(4,7,7,-6,0,-4)}(n)=H^{(3,11)}(33n)$\\
		$132=2^2\cdot3\cdot11$&$1089=3^2\cdot11^2$&$2R_{(6,7,10,7,3,6)}(n)=H^{(3,11)}(4n)$\\
		$132=2^2\cdot3\cdot11$&$1584=2^4\cdot3^2\cdot11$&$2R_{(5,5,17,-2,-2,-2)}(n)=H^{(3,11)}(11n)$\\
		$132=2^2\cdot3\cdot11$&$5808=2^4\cdot3\cdot11^2$&$2R_{(8,13,17,2,4,8)}(n)=H^{(3,11)}(3n)$\\
		$132=2^2\cdot3\cdot11$&$17424=2^4\cdot3^2\cdot11^2$&$2R_{(7,19,39,-18,-6,-2)}(n)=H^{(3,11)}(n)$\\
		$140=2^2\cdot5\cdot7$&$140=2^2\cdot5\cdot7$&$2R_{(1,2,18,-2,0,0)}(n)=H^{(70,1)}(140n)$\\
		$140=2^2\cdot5\cdot7$&$700=2^2\cdot5^2\cdot7$&$2R_{(5,6,6,-2,0,0)}(n)=H^{(70,1)}(28n)$\\
		$140=2^2\cdot5\cdot7$&$980=2^2\cdot5\cdot7^2$&$2R_{(6,6,7,0,0,-2)}(n)=H^{(70,1)}(20n)$\\
		$140=2^2\cdot5\cdot7$&$4900=2^2\cdot5^2\cdot7^2$&$2R_{(2,18,35,0,0,-2)}(n)=H^{(70,1)}(4n)$\\
		$156=2^2\cdot3\cdot13$&$156=2^2\cdot3\cdot13$&$2R_{(2,3,7,0,-2,0)}(n)=H^{(78,1)}(156n)$\\
		$156=2^2\cdot3\cdot13$&$468=2^2\cdot3^2\cdot13$&$2R_{(1,6,21,-6,0,0)}(n)=H^{(78,1)}(52n)$\\
		$156=2^2\cdot3\cdot13$&$2028=2^2\cdot3\cdot13^2$&$2R_{(2,7,39,0,0,-2)}(n)=H^{(78,1)}(12n)$\\
		$156=2^2\cdot3\cdot13$&$6084=2^2\cdot3^2\cdot13^2$&$2R_{(6,13,21,0,-6,0)}(n)=H^{(78,1)}(4n)$\\
		$160=2^5\cdot5$&$40=2^{3}\cdot5$&$2R_{(1, 3, 4, -2, 0, -1)}(n)=H^{(8,5;128)}(5n)$\\
		$160=2^5\cdot5$&$200=2^{3}\cdot5^2$&$2R_{(3, 3, 7, 3, 3, 1)}(n)=H^{(8,5;128)}(n)$\\
		$160=2^5\cdot5$&$5120=2^{10}\cdot5$&$2R_{(7, 15, 15, -2, -6, -6)}(n)=H^{(8,5)}(5n)$\\
		$160=2^5\cdot5$&$25600=2^{10}\cdot5^2$&$2R_{(11, 16, 44, 16, 4, 8)}(n)=H^{(8,5)}(n)$\\
		$168=2^3\cdot3\cdot7$&$42=2\cdot3\cdot7$&$R_{(1,1,11,-1,-1,0)}(n)=H^{(42,1)}(168n)$\\
		$168=2^3\cdot3\cdot7$&$126=2\cdot3^2\cdot7$&$R_{(3,3,5,-3,-3,0)}(n)=H^{(42,1)}(56n)$\\
		$168=2^3\cdot3\cdot7$&$294=2\cdot3\cdot7^2$&$R_{(5,5,5,-3,-3,-4)}(n)=H^{(42,1)}(24n)$\\
		$168=2^3\cdot3\cdot7$&$882=2\cdot3^2\cdot7^2$&$R_{(2,11,11,1,2,2)}(n)=H^{(42,1)}(8n)$\\
		$168=2^3\cdot3\cdot7$&$1344=2^6\cdot3\cdot7$&$R_{(4,4,23,-4,-4,0)}(n)=H^{(42,1)}(21n)$\\
		$168=2^3\cdot3\cdot7$&$4032=2^6\cdot3^2\cdot7$&$R_{(12,12,13,-12,-12,0)}(n)=H^{(42,1)}(7n)$\\
		$168=2^3\cdot3\cdot7$&$9408=2^6\cdot3\cdot7^2$&$R_{(12,17,17,6,12,12)}(n)=H^{(42,1)}(3n)$\\
		$168=2^3\cdot3\cdot7$&$28224=2^6\cdot3^2\cdot7^2$&$R_{(4,43,43,2,4,4)}(n)=H^{(42,1)}(n)$\\
		$180=2^2\cdot3^2\cdot5$&$45=3^2\cdot5$&$2R_{(1, 2, 7, 2, 1, 1)}(n)=H^{(45,1;9)}(20n)$\\
		$180=2^2\cdot3^2\cdot5$&$225=3^2\cdot5^2$&$2R_{(2, 5, 7, -5, -1, 0)}(n)=H^{(45,1;9)}(4n)$\\
		$180=2^2\cdot3^2\cdot5$&$405=3^4\cdot5$&$2R_{(2, 8, 8, 7, 1, 1)}(n)=H^{(45,1)}(20n)$\\
		$180=2^2\cdot3^2\cdot5$&$720=2^4\cdot3^2\cdot5$&$2R_{(4, 7, 8, -4, -4, 0)}(n)=H^{(45,1;9)}(5n)$\\
		$180=2^2\cdot3^2\cdot5$&$2025=3^4\cdot5^2$&$2R_{(7, 7, 13, 2, 7, 4)}(n)=H^{(45,1)}(4n)$\\
		$180=2^2\cdot3^2\cdot5$&$3600=2^4\cdot3^2\cdot5^2$&$2R_{(7, 7, 23, -2, -2, -6)}(n)=H^{(45,1;9)}(n)$\\
		$180=2^2\cdot3^2\cdot5$&$6480=2^4\cdot3^4\cdot5$&$2R_{(8, 11, 23, 2, 4, 8)}(n)=H^{(45,1)}(5n)$\\
		$180=2^2\cdot3^2\cdot5$&$32400=2^4\cdot3^4\cdot5^2$&$2R_{(7, 27, 52, -24, -4, -6)}(n)=H^{(45,1)}(n)$\\
		$184=2^3\cdot23$&$46=2\cdot23$&$2R_{(1, 3, 5, 3, 1, 1)}(n)=H^{(2,23)}(184n)$\\
		$184=2^3\cdot23$&$1058=2\cdot23^2$&$2R_{(5, 7, 10, 2, 4, 5)}(n)=H^{(2,23)}(8n)$\\
		$184=2^3\cdot23$&$1472=2^6\cdot23$&$2R_{(5, 8, 12, -8, -4, 0)}(n)=H^{(2,23)}(23n)$\\
		$184=2^3\cdot23$&$33856=2^6\cdot23^2$&$2R_{(11, 19, 51, -14, -6, -10)}(n)=H^{(2,23)}(n)$\\
		$208=2^4\cdot13$&$52=2^2\cdot13$&$2R_{(1, 3, 5, 1, 1, 1)}(n)=H^{(52,1;64)}(13n)$\\
		$208=2^4\cdot13$&$676=2^2\cdot13^2$&$2R_{(5, 5, 8, -2, -2, -3)}(n)=H^{(52,1;64)}(n)$\\
		$208=2^4\cdot13$&$3328=2^8\cdot13$&$2R_{(7, 12, 12, 8, 4, 4)}(n)=H^{(52,1)}(13n)$\\
		$208=2^4\cdot13$&$43264=2^8\cdot13^2$&$2R_{(11, 19, 59, -18, -10, -2)}(n)=H^{(52,1)}(n)$\\
		$240=2^4\cdot3\cdot5$&$60=2^2\cdot3\cdot5$&$2R_{(1, 4, 5, -4, -1, 0)}(n)=H^{(12,5;64)}(15n)$\\
		$240=2^4\cdot3\cdot5$&$180=2^2\cdot3^2\cdot5$&$2R_{(3, 5, 5, 5, 3, 3)}(n)=H^{(12,5;64)}(5n)$\\
		$240=2^4\cdot3\cdot5$&$300=2^2\cdot3\cdot5^2$&$2R_{(3, 5, 7, -5, -3, 0)}(n)=H^{(12,5;64)}(3n)$\\
		$240=2^4\cdot3\cdot5$&$900=2^2\cdot3^2\cdot5^2$&$2R_{(1, 15, 19, -15, -1, 0)}(n)=H^{(12,5;64)}(n)$\\
		$240=2^4\cdot3\cdot5$&$3840=2^8\cdot3\cdot5$&$2R_{(4, 5, 61, 2, 4, 4)}(n)=H^{(12,5)}(15n)$\\
		$240=2^4\cdot3\cdot5$&$11520=2^8\cdot3^2\cdot5$&$2R_{(12, 15, 23, 6, 12, 12)}(n)=H^{(12,5)}(5n)$\\
		$240=2^4\cdot3\cdot5$&$19200=2^8\cdot5^2\cdot3$&$2R_{(17, 17, 25, -10, -10, -14)}(n)=H^{(12,5)}(3n)$\\
		$240=2^4\cdot3\cdot5$&$57600=2^8\cdot3^2\cdot5^2$&$2R_{(16, 19, 64, 8, 16, 16)}(n)=H^{(12,5)}(n)$\\
		$252=2^2\cdot3^2\cdot7$&$252=2^2\cdot3^2\cdot7$&$4R_{(2, 5, 7, -2, -2, 0)}(n)=H^{(18,7;9)}(28n)$\\
		$252=2^2\cdot3^2\cdot7$&$1764=2^2\cdot3^2\cdot7^2$&$4R_{(5, 10, 10, 6, 2, 2)}(n)=H^{(18,7;9)}(4n)$\\
		$252=2^2\cdot3^2\cdot7$&$2268=2^2\cdot3^4\cdot7$&$4R_{(7, 7, 13, -2, -2, -4)}(n)=H^{(18,7)}(28n)$\\
		$252=2^2\cdot3^2\cdot7$&$15876=2^2\cdot3^4\cdot7^2$&$4R_{(10, 13, 34, 2, 10, 4)}(n)=H^{(18,7)}(4n)$\\
		$280=2^3\cdot5\cdot7$&$70=2\cdot5\cdot7$&$2R_{(1, 2, 9, 0, -1, 0)}(n)=H^{(70,1)}(280n)$\\
		$280=2^3\cdot5\cdot7$&$350=2\cdot5^2\cdot7$&$2R_{(3, 3, 10, 0, 0, -1)}(n)=H^{(70,1)}(56n)$\\
		$280=2^3\cdot5\cdot7$&$490=2\cdot5\cdot7^2$&$2R_{(3, 3, 14, 0, 0, -1)}(n)=H^{(70,1)}(40n)$\\
		$280=2^3\cdot5\cdot7$&$2240=2^6\cdot5\cdot7$&$2R_{(1, 8, 72, -8, 0, 0)}(n)=H^{(70,1)}(35n)$\\
		$280=2^3\cdot5\cdot7$&$2450=2\cdot5^2\cdot7^2$&$2R_{(1, 9, 70, 0, 0, -1)}(n)=H^{(70,1)}(8n)$\\
		$280=2^3\cdot5\cdot7$&$11200=2^6\cdot5^2\cdot7$&$2R_{(5, 24, 24, -8, 0, 0)}(n)=H^{(70,1)}(7n)$\\
		$280=2^3\cdot5\cdot7$&$15680=2^6\cdot5\cdot7^2$&$2R_{(7, 24, 24, -8, 0, 0)}(n)=H^{(70,1)}(5n)$\\
		$280=2^3\cdot5\cdot7$&$78400=2^6\cdot5^2\cdot7^2$&$2R_{(8, 35, 72, 0, -8, 0)}(n)=H^{(70,1)}(n)$\\
		$288=2^5\cdot3^2$&$72=2^{3}\cdot3^2$&$2R_{(1, 4, 5, -2, 1, 0)}(n)=H^{(72,1;1152)}(n)$\\
		$288=2^5\cdot3^2$&$648=2^{3}\cdot3^4$&$2R_{(5, 5, 8, 4, 4, 1)}(n)=H^{(72,1;128)}(n)$\\
		$288=2^5\cdot3^2$&$9216=2^{10}\cdot3^2$&$2R_{(11, 16, 16, 0, -8, -8)}(n)=H^{(72,1;9)}(n)$\\
		$288=2^5\cdot3^2$&$82944=2^{10}\cdot3^4$&$2R_{(16, 19, 76, 4, 16, 8)}(n)=H^{(72,1)}(n)$\\
		$312=2^3\cdot3\cdot13$&$78=2\cdot3\cdot13$&$2R_{(1, 5, 5, 4, 1, 1)}(n)=H^{(78,1)}(312n)$\\
		$312=2^3\cdot3\cdot13$&$234=2\cdot3^2\cdot13$&$2R_{(2, 3, 11, -3, -2, 0)}(n)=H^{(78,1)}(104n)$\\
		$312=2^3\cdot3\cdot13$&$1014=2\cdot3\cdot13^2$&$2R_{(1, 13, 23, -13, -1, 0)}(n)=H^{(78,1)}(24n)$\\
		$312=2^3\cdot3\cdot13$&$2496=2^6\cdot3\cdot13$&$2R_{(5, 5, 28, -4, -4, -2)}(n)=H^{(78,1)}(39n)$\\
		$312=2^3\cdot3\cdot13$&$3042=2\cdot3^2\cdot13^2$&$2R_{(3, 17, 17, 8, 3, 3)}(n)=H^{(78,1)}(8n)$\\
		$312=2^3\cdot3\cdot13$&$7488=2^6\cdot3^2\cdot13$&$2R_{(4, 7, 79, 2, 4, 4)}(n)=H^{(78,1)}(13n)$\\
		$312=2^3\cdot3\cdot13$&$32448=2^6\cdot3\cdot13^2$&$2R_{(8, 28, 41, 4, 8, 8)}(n)=H^{(78,1)}(3n)$\\
		$312=2^3\cdot3\cdot13$&$97344=2^6\cdot3^2\cdot13^2$&$2R_{(19, 19, 84, -12, -12, -14)}(n)=H^{(78,1)}(n)$\\
		$336=2^4\cdot3\cdot7$&$84=2^2\cdot3\cdot7$&$4R_{(1, 3, 8, -2, 0, -1)}(n)=H^{(28,3;64)}(21n)$\\
		$336=2^4\cdot3\cdot7$&$252=2^2\cdot3^2\cdot7$&$4R_{(3, 5, 5, -1, 0, -3)}(n)=H^{(28,3;64)}(7n)$\\
		$336=2^4\cdot3\cdot7$&$588=2^2\cdot3\cdot7^2$&$4R_{(3, 3, 17, -1, -1, -1)}(n)=H^{(28,3;64)}(3n)$\\
		$336=2^4\cdot3\cdot7$&$1764=2^2\cdot3^2\cdot7^2$&$4R_{(8, 9, 9, -3, -6, -6)}(n)=H^{(28,3;64)}(n)$\\
		$336=2^4\cdot3\cdot7$&$5376=2^8\cdot3\cdot7$&$4R_{(12, 12, 15, 12, 12, 8)}(n)=H^{(28,3)}(21n)$\\
		$336=2^4\cdot3\cdot7$&$16128=2^8\cdot3^2\cdot7$&$4R_{(5, 29, 29, 10, 2, 2)}(n)=H^{(28,3)}(7n)$\\
		$336=2^4\cdot3\cdot7$&$37632=2^8\cdot3\cdot7^2$&$4R_{(17, 20, 33, 20, 2, 4)}(n)=H^{(28,3)}(3n)$\\
		$336=2^4\cdot3\cdot7$&$112896=2^8\cdot3^2\cdot7^2$&$4R_{(11, 32, 92, 32, 4, 8)}(n)=H^{(28,3)}(n)$\\
		$504=2^3\cdot3^2\cdot7$&$126=2\cdot3^2\cdot7$&$4R_{(1, 5, 7, 2, 1, 1)}(n)=H^{(18,7;9)}(56n)$\\
		$504=2^3\cdot3^2\cdot7$&$882=2\cdot3^2\cdot7^2$&$4R_{(5, 5, 10, 2, 2, 3)}(n)=H^{(18,7;9)}(8n)$\\
		$504=2^3\cdot3^2\cdot7$&$1134=2\cdot3^4\cdot7$&$4R_{(5, 9, 9, 9, 3, 3)}(n)=H^{(18,7)}(56n)$\\
		$504=2^3\cdot3^2\cdot7$&$4032=2^6\cdot3^2\cdot7$&$4R_{(5, 8, 28, -8, -4, 0)}(n)=H^{(18,7;9)}(7n)$\\
		$504=2^3\cdot3^2\cdot7$&$7938=2\cdot3^4\cdot7^2$&$4R_{(5, 17, 26, 2, 4, 5)}(n)=H^{(18,7)}(8n)$\\
		$504=2^3\cdot3^2\cdot7$&$28224=2^6\cdot3^2\cdot7^2$&$4R_{(19, 19, 27, -6, -6, -18)}(n)=H^{(18,7;9)}(n)$\\
		$504=2^3\cdot3^2\cdot7$&$36288=2^6\cdot3^4\cdot7$&$4R_{(13, 28, 28, -16, -4, -4)}(n)=H^{(18,7)}(7n)$\\
		$504=2^3\cdot3^2\cdot7$&$254016=2^6\cdot3^4\cdot7^2$&$4R_{(19, 27, 136, -24, -16, -6)}(n)=H^{(18,7)}(n)$\\
		$528=2^4\cdot3\cdot11$&$132=2^2\cdot3\cdot11$&$4R_{(1, 5, 7, -1, 0, -1)}(n)=H^{(12,11;64)}(33n)$\\
		$528=2^4\cdot3\cdot11$&$396=2^2\cdot3^2\cdot11$&$4R_{(3, 5, 8, -2, 0, -3)}(n)=H^{(12,11;64)}(11n)$\\
		$528=2^4\cdot3\cdot11$&$1452=2^2\cdot3\cdot11^2$&$4R_{(7, 7, 8, 2, 2, 3)}(n)=H^{(12,11;64)}(3n)$\\
		$528=2^4\cdot3\cdot11$&$4356=2^2\cdot3^2\cdot11^2$&$4R_{(7, 13, 13, -7, -1, -1)}(n)=H^{(12,11;64)}(n)$\\
		$528=2^4\cdot3\cdot11$&$8448=2^8\cdot3\cdot11$&$4R_{(11, 11, 19, 2, 2, 6)}(n)=H^{(12,11)}(33n)$\\
		$528=2^4\cdot3\cdot11$&$25344=2^8\cdot3^2\cdot11$&$4R_{(17, 20, 20, -8, -4, -4)}(n)=H^{(12,11)}(11n)$\\
		$528=2^4\cdot3\cdot11$&$92928=2^8\cdot3\cdot11^2$&$4R_{(17, 32, 52, 32, 4, 8)}(n)=H^{(12,11)}(3n)$\\
		$528=2^4\cdot3\cdot11$&$278784=2^8\cdot3^2\cdot11^2$&$4R_{(19, 28, 139, 28, 2, 4)}(n)=H^{(12,11)}(n)$\\
		$672=2^5\cdot3\cdot7$&$168=2^{3}\cdot3\cdot7$&$4R_{(1, 4, 11, -2, -1, 0)}(n)=H^{(168,1;128)}(21n)$\\
		$672=2^5\cdot3\cdot7$&$504=2^{3}\cdot3^2\cdot7$&$4R_{(3, 5, 11, 4, 3, 3)}(n)=H^{(168,1;128)}(7n)$\\
		$672=2^5\cdot3\cdot7$&$1176=2^{3}\cdot3\cdot7^2$&$4R_{(5, 5, 13, 1, 1, 3)}(n)=H^{(168,1;128)}(3n)$\\
		$672=2^5\cdot3\cdot7$&$3528=2^{3}\cdot3^2\cdot7^2$&$4R_{(8, 11, 11, 1, 4, 4)}(n)=H^{(168,1;128)}(n)$\\
		$672=2^5\cdot3\cdot7$&$21504=2^{10}\cdot3\cdot7$&$4R_{(16, 16, 23, -8, -8, 0)}(n)=H^{(168,1)}(21n)$\\
		$672=2^5\cdot3\cdot7$&$64512=2^{10}\cdot3^2\cdot7$&$4R_{(13, 37, 37, -22, -2, -2)}(n)=H^{(168,1)}(7n)$\\
		$672=2^5\cdot3\cdot7$&$150528=2^{10}\cdot3\cdot7^2$&$4R_{(17, 41, 68, -36, -12, -10)}(n)=H^{(168,1)}(3n)$\\
		$672=2^5\cdot3\cdot7$&$451584=2^{10}\cdot3^2\cdot7^2$&$4R_{(16, 43, 172, 4, 16, 8)}(n)=H^{(168,1)}(n)$\\
		$720=2^4\cdot3^2\cdot5$&$180=2^{2}\cdot3^2\cdot5$&$4R_{(1, 7, 7, -3, 0, -1)}(n)=H^{(180,1;576)}(5n)$\\
		$720=2^4\cdot3^2\cdot5$&$900=2^{2}\cdot3^2\cdot5^2$&$4R_{(5, 7, 8, -2, 0, -5)}(n)=H^{(180,1;576)}(n)$\\
		$720=2^4\cdot3^2\cdot5$&$1620=2^{2}\cdot3^4\cdot5$&$4R_{(8, 9, 9, 9, 6, 6)}(n)=H^{(180,1;64)}(5n)$\\
		$720=2^4\cdot3^2\cdot5$&$8100=2^{2}\cdot3^4\cdot5^2$&$4R_{(7, 13, 27, 3, 6, 7)}(n)=H^{(180,1;64)}(n)$\\
		$720=2^4\cdot3^2\cdot5$&$11520=2^{8}\cdot3^2\cdot5$&$4R_{(7, 16, 31, -16, -6, 0)}(n)=H^{(180,1;9)}(5n)$\\
		$720=2^4\cdot3^2\cdot5$&$57600=2^{8}\cdot3^2\cdot5^2$&$4R_{(27, 28, 28, -24, -12, -12)}(n)=H^{(180,1;9)}(n)$\\
		$720=2^4\cdot3^2\cdot5$&$103680=2^{8}\cdot3^4\cdot5$&$4R_{(23, 32, 44, 32, 4, 8)}(n)=H^{(180,1)}(5n)$\\
		$720=2^4\cdot3^2\cdot5$&$518400=2^{8}\cdot3^4\cdot5^2$&$4R_{(27, 28, 187, 28, 6, 12)}(n)=H^{(180,1)}(n)$\\
	\end{longtable}
	
	\bibliographystyle{plain}
	\bibliography{TQF}
	
\end{document}